\documentclass[12pt]{amsart}
\usepackage[russian]{babel}
\usepackage{amsmath,amsfonts,amssymb}

\theoremstyle{plain}

\newtheorem{sct}{╤ртўєъ}[section]
\newtheorem{Lemma}[sct]{╦хььр}

\newtheorem{Theorem}[sct]{╥хюЁхьр}
\newtheorem{Corollary}[sct]{╤ыхфёЄтшх}

\newtheorem{Note}[sct]{╟рьхўрэшх}

\begin{document}

{\Large {\bf \centerline{╬яхЁрЄюЁ ─шЁрър}} {\bf \centerline{  ё
ъюьяыхъёэючэрўэ√ь ёєььшЁєхь√ь яюЄхэЎшрыюь.}}}

\bigskip

{\centerline{\textbf{╤ртўєъ └.~╠., ╪ърышъют └.~└.}\footnote{╧хЁт√щ ртЄюЁ яюффхЁцрэ уЁрэЄюь ╨╘╘╚, No 13-01-12476 ╬╘╚\_ь2, тЄюЁющ --- уЁрэЄюь ╨═╘ No
14-11-00754.}}}

{\centerline{\textit{╠юёъютёъшщ ├юёєфрЁёЄтхээ√щ ╙эштхЁёшЄхЄ, ╠юёътр, ╨юёёш }}}

{\centerline{\textit{─хърсЁ№ 2014}\footnote{╤ЄрЄ№  юяєсышъютрэр т Math. Notes \textbf{96} (5), 777--810 (2014).}}}

\bigskip

 ┬ ЁрсюЄх шчєўрхЄё  юяхЁрЄюЁ ─шЁрър, яюЁюцфхээ√щ эр юЄЁхчъх $[0,\pi]$ фшЇЇхЁхэЎшры№э√ь т√Ёрцхэшхь $-B\mathbf{y}'+Q(x)\mathbf{y}$, уфх
\begin{equation*}
B = \begin{pmatrix} 0 & 1 \\ -1 & 0 \end{pmatrix},
        \qquad
Q(x) = \begin{pmatrix} q_1(x) & q_2(x) \\ q_3(x) & q_4(x)        \end{pmatrix},
\end{equation*}
р ЇєэъЎшш $q_j(x)$ яЁшэрфыхцрЄ яЁюёЄЁрэёЄтє $L_p[0,\pi]$ фы  эхъюЄюЁюую $p\ge1$. ╩ырёё Ёхуєы Ёэ√ї ш ёшы№эю Ёхуєы Ёэ√ї юяхЁрЄюЁют Єръюую тшфр
юяЁхфхы хЄё  ёююЄтхЄёЄтє■∙шьш ъЁрхт√ьш єёыютш ьш. ┬ ЁрсюЄх яюыєўхэ√ рёшьяЄюЄшўхёъшх ЇюЁьєы√ фы  ёюсёЄтхээ√ї чэрўхэшщ ш ёюсёЄтхээ√ї ЇєэъЎшщ Єръшї
юяхЁрЄюЁют ё юЎхэърьш юёЄрЄъют, чртшё ∙шьш юЄ $p$. ─юърчрэю, ўЄю ёшёЄхьр ёюсёЄтхээ√ї ш яЁшёюхфшэхээ√ї ЇєэъЎшщ яЁюшчтюы№эюую Ёхуєы Ёэюую юяхЁрЄюЁр
юсЁрчєхЄ срчшё ╨шёёр ёю ёъюсърьш т яЁюёЄЁрэёЄтх $(L_2[0,\pi])^2$ ш юс√ўэ√щ срчшё ╨шёёр, яЁш єёыютшш ёшы№эющ Ёхуєы ЁэюёЄш юяхЁрЄюЁр.

\medskip
╩ы■ўхт√х ёыютр: \textit{╬яхЁрЄюЁ ─шЁрър, Ёхуєы Ёэ√х ъЁрхт√х єёыютш , рёшьяЄюЄшўхёъшх ЇюЁьєы√ фы  ёюсёЄтхээ√ї чэрўхэшщ ш ёюсёЄтхээ√ї ЇєэъЎшщ,
срчшё ╨шёёр.}

\vskip 1cm

\section*{┬тхфхэшх ш юсючэрўхэш .}

┬ ЁрсюЄх шчєўрхЄё  юяхЁрЄюЁ $L_Q$, яюЁюцфхээ√щ фшЇЇхЁхэЎшры№э√ь т√Ёрцхэшхь
\begin{equation}\label{eq:0.lQ}
l_Q(\mathbf{y})=-B\mathbf{y}'+Q\mathbf{y}
\end{equation}
т яЁюёЄЁрэёЄтх $\mathbb H=L_2[0,\pi]\oplus L_2[0,\pi]\ni \mathbf{y}$, уфх
\begin{equation}\label{eq:0.BQ}
    B = \begin{pmatrix} 0 & 1 \\ -1 & 0
        \end{pmatrix},
        \qquad
    Q(x) = \begin{pmatrix} q_1(x) & q_2(x) \\ q_3(x) & q_4(x)
        \end{pmatrix},
        \qquad
     \mathbf{y}(x)=\begin{pmatrix}y_1(x)\\ y_2(x)\end{pmatrix}.
\end{equation}
╘єэъЎшш $q_j$, $j=1,2,3,4$, яЁхфяюырур■Єё  ёєььшЁєхь√ьш эр юЄЁхчъх $[0,\pi]$ ш ъюьяыхъёэючэрўэ√ьш.  ┼ёЄхёЄтхээю, фы  ъюЁЁхъЄэюую юяЁхфхыхэш  ¤Єюую
юяхЁрЄюЁр  эєцэю фюяюыэшЄхы№эю чрфрЄ№ х∙х ъЁрхт√х єёыютш ,  эю юс ¤Єюь ь√ сєфхь уютюЁшЄ№ эшцх. ┬ ¤Єющ ЁрсюЄх ь√ ёЇюЁьєышЁєхь ш фюърцхь Ёхчєы№ЄрЄ√ ю
ётющёЄтрї юяхЁрЄюЁр ─шЁрър фы  ёыєўр  $q_j \in L_p[0, \pi]$  яЁш эхъюЄюЁюь  $p\geqslant 1$. ═ю юфэютЁхьхээю ь√ яюфуюЄютшь срчє  фы   шёёыхфютрэш 
¤Єюую юяхЁрЄюЁр ё яюЄхэЎшрырьш шч  фЁєушї ЇєэъЎшюэры№э√ї яЁюёЄЁрэёЄт, т ўрёЄэюёЄш фы  $q_j$ шч яЁюёЄЁрэёЄт ╤юсюыхтр $ W_p^\theta[0,\pi],
\theta\geqslant 0,\ \, p\geqslant 1$. ╩юэхўэю, яЁш $\theta >0$ яюыєўхээ√х чфхё№ Ёхчєы№ЄрЄ√ юс рёшьяЄюЄшўхёъюь яютхфхэшш ёюсёЄтхээ√ї ЇєэъЎшщ ш
ёюсёЄтхээ√ї чэрўхэшщ  ьюуєЄ с√Є№ ёє∙хёЄтхээю єЄюўэхэ√.

╤шёЄхьр \eqref{eq:0.lQ} с√ыр ттхфхэр т ЁрёёьюЄЁхэшх ╧.─шЁръюь т 1929 уюфє т ёт чш ё шчєўхэшхь Ёхы ЄштшёЄёъющ ьюфхыш ¤тюы■Ўшш ёяшэ-$1/2$ ўрёЄшЎ√ т
¤ыхъЄЁюьруэшЄэюь яюых. ╤рь ─шЁръ ЁрёёьрЄЁштры ёыєўрщ $q_2=q_3=0$, $q_1=V(x)-m$, $q_4(x)=V(x)+m$, уфх ЇєэъЎш  $V(x)$ юяшё√трхЄ яюЄхэЎшры яюы , р $m$
--- ьрёёр ўрёЄшЎ√ (ёь. \cite{Th}).   ╬яхЁрЄюЁ $L_Q$  шчєўрыё  тю ьэюушї ЁрсюЄрї, эю фю эхфртэхую тЁхьхэш т юёэютэюь  т ёыєўрх ёшььхЄЁшўхёъющ ьрЄЁшЎ√ $Q$  ё эхяЁхЁ√тэ√ьш
ЇєэъЎш ьш $q_j$  (ёь, эряЁшьхЁ, ьюэюуЁрЇш■ ╦хтшЄрэр ш ╤рЁуё эр \cite{LS} ш ёё√ыъш т эхщ). ╬яхЁрЄюЁє $L_Q$ ё яхЁшюфшўхёъшьш ш рэЄшяхЁшюфшўхёъшьш
єёыютш ьш  яюёт ∙хэр ёхЁш  ёЄрЄхщ ─№ ъютр ш ╠шЄ ушэр (ёь., эряЁшьхЁ, ЁрсюЄє \cite{DM1} ш ышЄхЁрЄєЁє т эхщ). ┬ ЁрсюЄрї \cite{DM2}---\cite{DM4} Єхї цх
ртЄюЁют шчєўхэ юяхЁрЄюЁ ─шЁрър ё яюЄхэЎшрыюь $Q\in L_2$ ш Ёхуєы Ёэ√ьш ъЁрхт√ьш єёыютш ьш. ┬ ўрёЄэюёЄш, т ёыєўрх $Q\in L_2$ фы  ёшы№эю Ёхуєы Ёэ√ї
єёыютшщ т ЁрсюЄх \cite{DM3} с√ыр фюърчрэр срчшёэюёЄ№ ╨шёёр ш ЄхюЁхьр юс рёшьяЄюЄшўхёъюь яютхфхэшш ёюсёЄтхээ√ї чэрўхэшщ (юэр рэрыюушўэр яхЁтющ ўрёЄш
╥хюЁхь√ \ref{tm:4.1} ¤Єющ ЁрсюЄ√). ╧Ёш юЄёєЄёЄтшш ёшы№эющ Ёхуєы ЁэюёЄш с√ыр фюърчрэр срчшёэюёЄ№ ╨шёёр шч яюфяЁюёЄЁрэёЄт. ─ы  ёыєўр  $q_j\in
L_1[0,\pi]$ юЄьхЄшь ЁрсюЄє └ы№схтхЁшю, ├Ёшэштр ш ╠шъшЄ■ър \cite{AHM}, т ъюЄюЁющ ьхЄюфюь юяхЁрЄюЁют яЁхюсЁрчютрэш  шчєўхэр юсЁрЄэр  чрфрўр
тюёёЄрэютыхэш  тх∙хёЄтхээючэрўэющ ьрЄЁшЎ√ $Q, q_2 =q_3=0$, яю фтєь ёяхъЄЁрь юяхЁрЄюЁр $L_Q$ ё ъЁрхт√ьш єёыютш ьш ─шЁшїых $u_1(0)=u_1(\pi)=0$ ш
─шЁшїых--═хщьрэр $u_1(0)=u_2(\pi)=0$.

╩ ёшёЄхьх ─шЁрър яЁшьхэшь√ Ёхчєы№ЄрЄ√, яюыєўхээ√х фы  юс∙шї ёшёЄхь яхЁтюую яюЁ фър ё ьрЄЁшЎхщ $n\times n$. ▌Єш шёёыхфютрэш  схЁєЄ ётюх эрўрыю ё ЁрсюЄ
┴шЁъуюЇр \cite{B1, B2}, ╥рьрЁъшэр \cite{Ta1, Ta2}, ┴шЁъуюЇр ш ╦рэухЁр \cite{BirLan}. ┬ ўрёЄэюёЄш, т ¤Єшї ЁрсюЄрї с√ыш ттхфхэ√ яюэ Єш  Ёхуєы Ёэ√ї ш
ёшы№эю Ёхуєы Ёэ√ї ъЁрхт√ї єёыютшщ, яЁш ъюЄюЁ√ї юяхЁрЄюЁ $L_Q$ шчєўрхЄё  т ¤Єющ ЁрсюЄх. ╨рсюЄ√ ─рэЇюЁфр \cite{Du}, ╠шїрщыютр \cite{Mikh} ш
├.~╩хёхы№ьрэр \cite{Kes} с√ыш яхЁт√ьш, т ъюЄюЁ√ї с√ыю яЁютхфхэю  шчєўхэшх
 ётющёЄт срчшёэюёЄш ёшёЄхь√ ъюЁэхт√ї тхъЄюЁют фы  ёшы№эю Ёхуєы Ёэ√ї юс√ъэютхээ√ї
фшЇЇхЁхэЎшры№э√ї юяхЁрЄюЁют. ╬сюс∙хэшх ¤Єшї Ёхчєы№ЄрЄют ш шї ЁрчтшЄшх с√ыю яЁютхфхэю т ЁрсюЄрї ╪ърышъютр \cite{S1, S2}. ╠хЄюф√ шч \cite{S1}  ь√
шёяюы№чєхь ш т ¤Єющ ЁрсюЄх. ╬с∙шх Ёхчєы№ЄрЄ√ ю срчшёэюёЄш фы  тючьє∙хэшщ ёрьюёюяЁ цхээ√ї юяхЁрЄюЁют т ушы№схЁЄютюь яЁюёЄЁрэёЄтх с√ыш яюыєўхэ√
╠рЁъєёюь \cite{M1} ш ╩рЎэхы№ёюэюь \cite{Kaz}, р яючфэхх с√ыш юсю∙хэ√ ╠рЁъєёюь ш ╠рЎрхт√ь \cite{MM} (ёь. Єръцх \cite{M2}). └уЁрэютшў \cite{A1, A2}
ёфхыры трцэ√х фюсртыхэш  ъ ¤Єшь ЄхюЁхьрь. ═ют√х шфхш т ¤Єющ Єхьх фы  юяхЁрЄюЁют ўрёЄэюую тшфр с√ыш яЁхфыюцхэ√  └фєўўш ш ╠шЄ ушэ√ь \cite{AdM}, р фы 
юс∙шї юяхЁрЄюЁют  ╪ърышъют√ь \cite{S3}.
 ╚ч яюёыхфэшї ЁрсюЄ юс юс√ъэютхэээ√ї фшЇЇхЁхэЎшры№э√ї юяхЁрЄюЁрї, яюЁюцфрхь√ї
 ьрЄЁшЎхщ $n\times n$, юЄьхЄшь ЁрсюЄ√ ╠рырьєфр, ╬ЁшфюЁюуш ш ╦єэхтр \cite{MalOri, ML1}. ─ы  Єръшї юяхЁрЄюЁют т эшї
фюърчрэ√ ЄхюЁхь√ ю яюыэюЄх, ьшэшьры№эюёЄш ш срчшёэюёЄш
ёшёЄхь√ ъюЁэхт√ї тхъЄюЁют (фы  ёыєўр  ёшёЄхь√ ─шЁрър ЁрёёьюЄЁхэ яюЄхэЎшры $P\in L_\infty$ ш рэюэёшЁютрэ Ёхчєы№ЄрЄ фы  $P\in L_2$). ╬ЄьхЄшь Єръцх  ЁрсюЄ√ └ьшЁютр ш ├єёхщэютр \cite{AG}, ┴рёъръютр, ─хЁсє°хтр ш ┘хЁсръютр \cite{Ba12},
 ╩юЁэхтр ш ╒Ёюьютр \cite{Kh10},  ╥Ёє°шэр ш ▀ьрьюЄю  \cite{TrYa}, т ъюЄюЁ√ї ўшЄрЄхы№ ьюцхЄ эрщЄш шэЄхЁхёэ√х Ёхчєы№ЄрЄ√,
 ёт чрээ√х ё юяхЁрЄюЁюь ─шЁрър .

╬ёэютэр  Ўхы№ эрёЄю ∙хщ ЁрсюЄ√ --- эрщЄш рёшьяЄюЄшъш ёюсёЄтхээ√ї чэрўхэшщ ш ёюсёЄтхээ√ї ЇєэъЎшщ  Ёхуєы Ёэюую юяхЁрЄюЁр  $L_Q$ яЁш $q_j\in L_p[0,\pi],
p\geqslant 1$, р Єръцх фюърчрЄ№ срчшёэюёЄ№ ╨шёёр ёюсёЄтхээ√ї ЇєэъЎшщ  (ёю ёъюсърьш шыш схч эшї)  т яЁюёЄЁрэёЄтх $\mathbb H$ яЁш ы■сюь $p\geqslant 1$.
▌Єш Ёхчєы№ЄрЄ√ с√ыш яюыєўхэ√ ртЄюЁрьш х∙х т 2011 уюфє (т ўрёЄэюёЄш, ╥хюЁхьр \ref{tm:2.9} с√ыр рэюэёшЁютрэр  т \cite{Sav1}). ╬фэръю, эрёЄю ∙р  тхЁёш 
ЁрсюЄ√ с√ыр чряшёрэр Єюы№ъю т ртуєёЄх  2013 уюфр. ╥юуфр цх юэр с√ыр яюёырэр эхъюЄюЁ√ь ёяхЎшрышёЄрь фы  ючэръюьыхэш . ┬ ўрёЄэюёЄш, ь√ шёъЁхээх
сыруюфрЁшь яЁюЇхёёюЁр ┴.~╤.~╠шЄ ушэр,  ёфхырт°хую  эхёъюы№ъю Ўхээ√ї чрьхўрэшщ. ╧єсышърЎш■ ЁрсюЄ√  ртЄюЁ√ юЄъырф√трыш, эрфх ё№ эрщЄш сюыхх яЁюёЄюх
фюърчрЄхы№ёЄтю Єхїэшўхёъш ёыюцэющ ╦хьь√ \ref{lem:2.1} (юЄьхЄшь, ўЄю ¤Єр ыхььр  ты хЄё  ъы■ўхтющ т ЁрсюЄх). ╤фхырЄ№ ¤Єю эх єфрыюё№ ш ртЄюЁ√ Ёх°шыш
юяєсышъютрЄ№ шёїюфэє■ тхЁёш■ фюърчрЄхы№ёЄтр. ╬фэръю, яЁш яюшёъх сюыхх яЁюёЄюую фюърчрЄхы№ёЄтр ╦хьь√ \ref{lem:2.1}, ртЄюЁрь єфрыюё№ эрщЄш фЁєующ
яюфїюф ъ яюыєўхэш■ эєцэ√ї рёшьяЄюЄшъ, ъюЄюЁ√щ юёэютрэ эр фЁєушї шфх ї ш сєфхЄ яЁхфёЄртыхэ т фЁєушї ЁрсюЄрї. ═хфртэю ртЄюЁрь ёЄрыю шчтхёЄэю ю ЁрсюЄх
╦єэхтр ш ╠рырьєфр \cite{ML2}, т ъюЄюЁющ юэш рэюэёшЁютрыш Ёхчєы№ЄрЄ ю срчшёэюёЄш ╨шёёр ёшcЄхь√ ъюЁэхт√ї тхъЄюЁют юяхЁрЄюЁр ─шЁрър, яюЁюцфхээюую ёшы№эю
Ёхуєы Ёэ√ьш ъЁрхт√ьш  єёыютш ьш,  ё яюЄхэЎшрыюь $Q\in L_1[0,\pi]$. ╚фх  фюърчрЄхы№ёЄтр т \cite{ML2} юЄышўэр юЄ эр°хщ.  ╩Ёюьх Єюую, ь√ фюърч√трхь
сюыхх юс∙шщ Ёхчєы№ЄрЄ, эх яЁхфяюырур  ёшы№эє■ Ёхуєы ЁюёЄ№, р Єюы№ъю Ёхуєы ЁэюёЄ№. ╬ЄьхЄшь Єръцх, ўЄю яЁш фюърчрЄхы№ёЄтх рёшьяЄюЄшўхёъшї ЇюЁьєы фы 
ёюсёЄтхээ√ї чэрўхэшщ ш ёюсёЄтхээ√ї ЇєэъЎшщ, эрьш яюыєўхэю ёє∙хёЄтхээю сюыхх Єюўэ√х юЎхэъш юёЄрЄъют, т Єю тЁхь  ъръ фы  фюърчрЄхы№ёЄтр срчшёэюёЄш
╨шёёр фюёЄрЄюўэю юЎхэюъ тшфр $o(1)$ яЁш $n\to\infty$.

═р яЁюЄ цхэшш тёхщ ёЄрЄ№ш ь√ сєфхь шёяюы№чютрЄ№ ёыхфє■∙шх юсючэрўхэш . ╫хЁхч $\mathbf{y}(x)=(y_1(x),y_2(x))^t$ сєфхь юсючэрўрЄ№ тхъЄюЁ-ЇєэъЎшш эр
$[0,\pi]$, ўхЁхч $\langle\cdot,\,\cdot\rangle$
--- ёъры Ёэюх яЁюшчтхфхэшх т яЁюёЄЁрэёЄтх $\mathbb H$. ╫Єюс√ эх
єёыюцэ Є№ чряшё№ ь√ ўрёЄю сєфхь   яшёрЄ№ ${\bold f \in L_p}$, яЁхфяюырур , ўЄю  ${\bold f} \in L_p\times L_p$, шыш $Q\in L_p$, яЁхфяюырур , ўЄю тёх
ъюьяюэхэЄ√ ьрЄЁшЎ√ шч $L_p$. ═юЁьє т $L_p$ шыш т $L_p\times L_P$  сєфхь юсючэрўрЄ№ $\|\cdot\|_p$.  ┼ёыш шч ъюэЄхъёЄр  ёэю, т ъръюь яЁюёЄЁрэёЄтх
ЇєэъЎшщ (тхъЄюЁэ√ї шыш ёъры Ёэ√ї) ь√ ЁрсюЄрхь, ь√ эх фхырхь Ёрчышўш  т юсючэрўхэш ї ьхцфє $L_p$ ш $L_p\times L_p$.

\vskip 0,3cm

\section{┬ёяюьюурЄхы№э√х Ёхчєы№ЄрЄ√.}

\vskip 0,3cm

╤ фшЇЇхЁхэЎшры№э√ь т√Ёрцхэшхь $l_Q$  ёт цхь ьръёшьры№э√щ юяхЁрЄюЁ
$$
L_{Q,M}\, \mathbf{y}:=l_Q(\mathbf{y});\qquad\frak{D}(L_{Q,M})=\{\mathbf{y}\in AC[0,\pi]\vert\,\ \, l(\mathbf{y})\in\mathbb H\}
$$
ш ьшэшьры№э√щ юяхЁрЄюЁ $L_{Q,m}$,  ты ■∙шщё  ёєцхэшхь ьръёшьры№эюую юяхЁрЄюЁр эр юсырёЄ№
$$
\frak{D}(L_m)=\{\mathbf{y}\in\frak{D}(L_M)\vert\,\mathbf{y}(0)=\mathbf{y}(\pi)=0\}.
$$
╟фхё№ $AC[0,\pi] =W^1_1[0,\pi]$ --- яЁюёЄЁрэёЄтю рсёюы■Єэю эхяЁхЁ√тэ√ї ЇєэъЎшщ,   р ¤ыхьхэЄ√ ьрЄЁшЎ√ $Q$ ---  яЁюшчтюы№э√х ёєььшЁєхь√х
ъюьяыхъёэючэрўэ√х ЇєэъЎшш. ┬ ¤Єюь ёыєўрх юср ёырурхь√ї фшЇЇхЁхэЎшры№эюую т√Ёрцхэш  $l_Q(\mathbf{y})$  ъюЁЁхъЄэю юяЁхфхыхэ√, ъръ ЇєэъЎшш шч $L_1
\times L_1$. ═ю т юсырёЄ№ юяЁхфхыхэш  тїюф Є Єюы№ъю Єх ЇєэъЎшш $\mathbf{y}$,  фы  ъюЄюЁ√ї ёєььр ¤Єшї ёырурхь√ї яЁшэрфыхцшЄ $\mathbb H$.  ╫хЁхч
$L_{\overline Q , M}$ ш $L_{\overline Q , m}$ сєфхь юсючэрўрЄ№ ьръёшьры№э√щ ш ьшэшьры№э√щ юяхЁрЄюЁ√, яюЁюцфхээ√х ёюяЁ цхээ√ь фшЇЇхЁхэЎшры№э√ь
т√Ёрцхэшхь
$$
l_{\overline{Q}}(\mathbf{y}):=-B\mathbf{y}'+\overline{Q}\mathbf{y},\qquad\text{уфх}\ \overline{Q}=\begin{pmatrix}\overline{q_1} & \overline{q_3}\\
\overline{q_2} & \overline{q_4}\end{pmatrix}.
$$

\begin{Lemma}\label{lem:1.0} ╧єёЄ№  $\varphi$ --- яЁюшчтюы№эр  рсёюы■Єэю эхяЁхЁ√тэр  эр $[0,\pi]$ ЇєэъЎш . ╥юуфр
юяхЁрЄюЁ $L_Q$ єэшЄрЁэю ¤ътштрыхэЄхэ юяхЁрЄюЁє $L_{\widetilde Q}$, уфх
\begin{gather*}
\widetilde Q=\begin{pmatrix}\widetilde{q_1}&\widetilde{q_2}\\ \widetilde{q_3}&\widetilde{q_4}\end{pmatrix},\\
\widetilde{q_1}=-\dot{\varphi}+q_1\cos^2\varphi-(q_2+q_3)\sin\varphi\cos\varphi+q_4\sin^2\varphi,\\
\widetilde{q_2}=q_2\cos^2\varphi+(q_1-q_4)\sin\varphi\cos\varphi-q_3\sin^2\varphi,\\
\widetilde{q_3}=q_3\cos^2\varphi+(q_1-q_4)\sin\varphi\cos\varphi-q_2\sin^2\varphi,\\
\widetilde{q_4}=-\dot{\varphi}+q_4\cos^2\varphi+(q_2+q_3)\sin\varphi\cos\varphi+q_1\sin^2\varphi.
\end{gather*}
\end{Lemma}
\begin{proof}
▌Єю єЄтхЁцфхэшх їюЁю°ю шчтхёЄэю; ёь. \cite[├ы. 1]{LS}. ╨рёёьюЄЁшь т яЁюёЄЁрэёЄтх $\mathbb H$ єэшЄрЁэ√х юяхЁрЄюЁ√
$$
H=\begin{pmatrix}\cos\varphi & \sin\varphi\\ -\sin\varphi & \cos\varphi\end{pmatrix}, \qquad
H^{-1}=\begin{pmatrix}\cos\varphi & -\sin\varphi\\ \sin\varphi & \cos\varphi\end{pmatrix}.
$$
╥юуфр
$$
H^{-1}L_QH=-H^{-1}BH\frac{d}{dx}+(H^{-1}QH-H^{-1}B\dot{H})=-B\frac{d}{dx}+(H^{-1}QH-\dot{\varphi}I),
$$
уфх $I$ --- ЄюцфхёЄтхээ√щ юяхЁрЄюЁ т $\mathbb H$ ═хяюёЁхфёЄтхээ√ь яюфёўхЄюь яюыєўрхь юЄё■фр єЄтхЁцфхэшх ыхьь√.
\end{proof}
\begin{Note}
╟рьхЄшь ЄхяхЁ№, ўЄю $\widetilde{q_1}+\widetilde{q_4}=-2\dot{\varphi}+q_1+q_4$, Єръ ўЄю т√сЁрт
$$
\varphi=\frac12\int\limits_0^x(q_1+q_4)d\xi,
$$
яюыєўшь ЁртхэёЄтю $\widetilde{q_4}=-\widetilde{q_1}$. ─рыхх, трцэю чрьхЄшЄ№, ўЄю яЁш Єръюь т√сюЁх $\varphi$ ъырёё уырфъюёЄш яюЄхэЎшрыр
$\widetilde{Q}$ ёютярфрхЄ ё ъырёёюь яюЄхэЎшрыр $Q$, эряЁшьхЁ $Q\in W_p^\theta$ Єюуфр ш Єюы№ъю Єюуфр, ъюуфр $\widetilde{Q}\in W_p^\theta$. ═ръюэхЎ,
чрьхЄшь (¤Єю трцэю фы  юяшёрэш  ёрьюёюяЁ цхээ√ї Ёрё°шЁхэшщ юяхЁрЄюЁр $L_{Q, m}$), ўЄю яЁш юЄюсЁрцхэшш $\mathbf{y}\mapsto H\mathbf{y}=\mathbf{z}$
яЁюёЄЁрэёЄтр $\mathbf{y}\in AC[0,\pi]$ т√яюыэхэю:
$$
\mathbf{z}(0)=\mathbf{y}(0), \qquad\mathbf{z}(\pi)=\begin{pmatrix}\cos\varphi(\pi) & \sin\varphi(\pi)\\-\sin\varphi(\pi) &
\cos\varphi(\pi)\end{pmatrix}\mathbf{y}(\pi).
$$
╧ю¤Єюьє, ёфхырт, хёыш эхюсїюфшью,  чрьхэє ёяхъЄЁры№эюую ярЁрьхЄЁр $\widetilde{\lambda}=\lambda+c$ ь√ тёхуфр ьюцхь ёўшЄрЄ№, ўЄю
$$
\varphi(\pi)=\frac12\int_0^\pi (q_1+q_4)d\xi=0,
$$
Є.х. ўЄю яюёых чрьхэ√ $ \mathbf{z} = H\mathbf{y}$ т√яюыэ хЄё  ЁртхэёЄтю $\mathbf{z}(\pi)=\mathbf{y}(\pi)$.
\end{Note}
┬ё■фє фрыхх т ЁрсюЄх ь√ ёўшЄрхь, ўЄю ЇєэъЎшш $q_j$ єфютыхЄтюЁ ■Є ёыхфє■∙шь єёыютш ь\\
(i)\quad ╘єэъЎшш $q_2$ ш $q_3$ ёєььшЁєхь√ эр $[0,\pi]$;\\
(ii)\quad ╘єэъЎшш $q_1$ ш $q_2+q_3$ яЁшэрфыхцрЄ °рЁє $B(0,R)$ яЁюёЄЁрэёЄтр $L_p[0,\pi]$ фы  эхъюЄюЁюую $p\ge 1$ ш $R>0$.\\
(iii)\quad $q_4=-q_1 $ (ъръ с√ыю юЄьхўхэю т√°х, ¤Єю єёыютшх эх юуЁрэшўштрхЄ юс∙эюёЄш).

\begin{Lemma}[╘юЁьєыр ╦руЁрэцр]\label{lem:1.1}
─ы  ЇєэъЎшщ \(\mathbf{f}\in\frak{D}(L_{Q,M})\), \(\mathbf{g}\in\frak{D}(L_{\overline Q, M})\) ёяЁртхфыштю ЄюцфхёЄтю
\[
    \int_0^\pi L_{Q,M}(\mathbf{f})\ \overline{\mathbf{g}}\, dx=\int_0^\pi\mathbf{f}\ \overline{L_{\overline Q,
M}(\mathbf{g})}\,dx+[\mathbf{f},\mathbf{g}]^{\pi}_0,
\]
уфх \([\mathbf{f},\mathbf{g}]_0^{\pi}=-\left.f_2(x)\overline{g_1(x)}\right|_0^{\pi}- \left.\overline{g_2(x)}f_1(x)\right|_0^{\pi}\).
\end{Lemma}
╚ч ¤Єющ ЇюЁьєы√, т ўрёЄэюёЄш, яюыєўрхь
\begin{align*}
    \langle L_{Q,M}(\mathbf{f}),\mathbf{g}\rangle&=\langle\mathbf{f},L_{\overline Q,m}(\mathbf{g})\rangle,&\mathbf{f}&\in
    \frak{D}(L_{Q,M}),\;\mathbf{g}\in\frak{D}(L_{\overline Q,
m}),
\end{align*}
Є.~х. юяхЁрЄюЁ√ \(L_{Q,M}\) ш \(L_{\overline Q,m}\) тчршьэю ёюяЁ цхэ√.

\vskip 0,2cm

┬ фры№эхщ°хь трцэє■ Ёюы№ шуЁрхЄ ёыхфє■∙шщ Ёхчєы№ЄрЄ.

\begin{Theorem}\label{tm:0} ╧єёЄ№ $\mathbf{A}(x)$ --- ьрЄЁшЎр ЁрчьхЁр $n\times n$,
¤ыхьхэЄ√ ъюЄюЁющ  ты ■Єё  ЇєэъЎш ьш яЁюёЄЁрэёЄтр $L_1[0,\pi]$, р $\mathbf{f}\in \big [ L_1[0,\pi]\big ]^n$ --- тхъЄюЁ-ЇєэъЎш . ╥юуфр яЁш ы■сюь
$c\in[0,\pi]$ єЁртэхэшх
$$
\mathbf{y}'=\mathbf{A}(x)\mathbf{y}+\mathbf{f}, \mathbf{y}(c)=\mathbf{\xi}\in\mathbb C^n,
$$
шьххЄ хфшэёЄтхээюх Ёх°хэшх $\mathbf{y}(x)$, яЁшўхь $\mathbf{y}(x)$ --- рсёюы■Єэю эхяЁхЁ√тэр  эр $[0,\pi]$ тхъЄюЁ--ЇєэъЎш . ┼ёыш яюёыхфютрЄхы№эюёЄ№
ьрЄЁшЎ $\mathbf{A}_\varepsilon(x)$ ё ¤ыхьхэЄрьш шч $L_1[0,\pi]$ Єръютр, ўЄю $\|\mathbf{A}_\varepsilon(x)-\mathbf{A}(x)\|_{L_1}\to 0$ яЁш
$\varepsilon\to 0$, Єю Ёх°хэш  єЁртэхэшщ
\begin{align*}
    \mathbf{y}_\varepsilon'&=\mathbf{A}_\varepsilon(x)\mathbf{y}_\varepsilon+\mathbf{f},&
    \mathbf{y}_\varepsilon(c)=\mathbf{\xi},
\end{align*}
ёїюф Єё  ъ $\mathbf{y}(x)$ ЁртэюьхЁэю эр $[0,\pi]$ (ш фрцх т ьхЄЁшъх яЁюёЄЁрэёЄтр $W_1^1[0,\pi]$).  ╩Ёюьх Єюую, ёяЁртхфыштр юЎхэър
\begin{equation*}
    \|\mathbf{y}(x)-\mathbf{y}_\varepsilon(x)\|_{1,1}\le
    C(\|\mathbf{f}\|_{L_1}+\|\xi\|)\|\mathbf{A}(x)-\mathbf{A}_\varepsilon(x)\|_{L_1}
\end{equation*}
ё яюёЄю ээющ $C$, эх чртшё ∙хщ юЄ $\mathbf{f}$ ш $\varepsilon$.
\end{Theorem}
\begin{proof} ╧хЁтюх єЄтхЁцфхэшх ЄхюЁхь√ їюЁю°ю шчтхёЄэю (ёь., эряЁшьхЁ, \cite [├ы. 8] {CL}). ─юърчрЄхы№ёЄтю тЄюЁюую єЄтхЁцфхэш  ёюфхЁцшЄё  т ЁрсюЄх
\cite{SS99}.
\end{proof}

═ряюьэшь, ўЄю юяхЁрЄюЁ \(F\), фхщёЄтє■∙шщ т ушы№схЁЄютюь (шыш срэрїютюь) яЁюёЄЁрэёЄтх \(\mathfrak H\), эрч√трхЄё  ЇЁхфуюы№ьют√ь, хёыш хую юсырёЄ№
юяЁхфхыхэш  \(\mathfrak D(F)\) яыюЄэр т \(\mathfrak H\), юсЁрч чрьъэєЄ, р фхЇхъЄэ√х ўшёыр \(\{\alpha,\beta\}\), Ёртэ√х ЁрчьхЁэюёЄ ь  фЁр ш ъю фЁр,
ъюэхўэ√.

\begin{Theorem}\label{tm:1.2} ╧Ёш ы■сюь \(\lambda\in\mathbb C\) юяхЁрЄюЁ√
\(L_{Q, M}-\lambda\) ш \( L_{\overline Q, m} -\overline{\lambda}\) ЇЁхфуюы№ьют√,  ты ■Єё  тчршьэю ёюяЁ цхээ√ьш, р шї фхЇхъЄэ√х ўшёыр Ёртэ√
\(\{0,2\}\) ш \(\{2,0\}\), ёююЄтхЄёЄтхээю.
\end{Theorem}
\begin{proof} ─юърчрЄхы№ёЄтю ¤Єющ ЄхюЁхь√ ыхуъю яюыєўрхЄё  ё яюью∙№■
яЁхф√фє∙хщ ЄхюЁхь√ Єюўэю Єръцх, ъръ ¤Єю ёфхырэю т ЁрсюЄх
\cite{SS99}.
\end{proof}

╧хЁхщфхь ъ юяшёрэш■ Ёрё°шЁхэшщ юяхЁрЄюЁр $L_{Q,m}$. ═рё шэЄхЁхёє■Є юяхЁрЄюЁ√ $L_Q$ ё эєыхт√ьш фхЇхъЄэ√ьш ўшёырьш, ъюЄюЁ√х  ты ■Єё  юфэютЁхьхээю
эхЄЁштшры№э√ь Ёрё°шЁхэшхь юяхЁрЄюЁр $L_m$ ш эхЄЁштшры№э√ь ёєцхэшхь юяхЁрЄюЁр $L_M$, Є.х. $L_m\subset L\subset L_M$. ╚ч ╥хюЁхь√ \ref{tm:1.2} ш шч
юяЁхфхыхэшщ ьръёшьры№эюую ш ьшэшьры№эюую юяхЁрЄюЁют ёыхфєхЄ, ўЄю ы■сющ Єръющ юяхЁрЄюЁ $L$ шьххЄ юсырёЄ№ юяЁхфхыхэш 
\begin{equation} \label{dom}
    \mathfrak D(L)=\left\{\mathbf{y}\;\vline\;\mathbf{y}\in\mathfrak D(L_M),\;U(\mathbf{y})=0\right\},
\end{equation}
 уфх \(U(\mathbf{y})\) --- эхъюЄюЁр  эхэєыхтр  ышэхщэр  ЇюЁьр юЄ
яхЁхьхээ√ї \(\mathbf{y}(0)\) ш \(\mathbf{y}(\pi)\). ▀ёэю, ўЄю ы■ср  Єрър  ЇюЁьр ьюцхЄ с√Є№ чряшёрэр т тшфх
\begin{equation}\label{matrU}
U(\mathbf{y}(0),\mathbf{y}(\pi))=A\mathbf{y}(0)+B\mathbf{y}(\pi), \text{уфх}\quad A=\begin{pmatrix}u_{11} & u_{12}\\ u_{21} & u_{22}\end{pmatrix},
B=\begin{pmatrix}u_{13} & u_{14}\\ u_{23} & u_{24}\end{pmatrix}.
\end{equation}
╬сючэрўшь ўхЁхч \(J_{\alpha\beta}\) юяЁхфхышЄхы№, ёюёЄртыхээ√щ шч
\(\alpha\)-ую ш \(\beta\)-ую ёЄюысЎр ьрЄЁшЎ√
\[
    \begin{pmatrix}u_{11}&u_{12}&u_{13}&u_{14}\\
    u_{21}&u_{22}&u_{23}&u_{24}\end{pmatrix}.
\]
ышэхщэющ ЇюЁь√ $U$,
ттхфхээющ т \eqref{matrU}
╩Ёрхтюх єёыютшх, юяЁхфхыхээюх ЇюЁьющ $U$ эрчютхь {\it Ёхуєы Ёэ√ь},
хёыш эх  Ёртэ√ эєы■ фтр  ўшёыр
\begin{equation}\label{regular}
J_{14}+J_{32}+i(J_{42}-J_{13})\quad\text{ш}\
J_{14}+J_{32}-i(J_{42}-J_{13}).
\end{equation}
╩юэхўэю, хёыш  тёх ъю¤ЇЇшЎшхэЄ√ $u_{kj}$, єўрёЄтє■∙шх т ышэхщэющ ЇюЁьх $U$, тх∙хёЄтхээ√, Єю фюёЄрЄюўэю юЄышўш  юЄ эєы  Єюы№ъю юфэюую шч ¤Єшї ўшёхы.
╬яхЁрЄюЁ ─шЁрър, яюЁюцфхээ√щ Ёхуєы Ёэ√ь ъЁрхт√ь єёыютшхь $U$ (Є.х. юяхЁрЄюЁ $L_Q$  ё юсырёЄ№■ \eqref{dom}) сєфхь эрч√трЄ№ {\it  Ёхуєы Ёэ√ь}. ╧ючцх ь√
т√фхышь Єръцх {\it ёшы№эю Ёхуєы Ёэ√х} юяхЁрЄюЁ√ ─шЁрър. ╟фхё№ Єюы№ъю юЄьхЄшь, ўЄю фртрЄ№ юяЁхфхыхэшх ёшы№эю Ёхуєы Ёэюую ъЁрхтюую єёыютш  ъюЁЁхъЄэю
Єюы№ъю яЁш фюяюыэшЄхы№эюь єёыютшш $q_2=q_3$. ┬ яЁюЄштэюь ёыєўрх, яЁш юфэюь ш Єюь цх ъЁрхтюь єёыютшш юяхЁрЄюЁ ─шЁрър ё яюЄхэЎшрыюь $Q$  ьюцхЄ с√Є№
ёшы№эю Ёхуєы Ёэ√ь, р ё фЁєушь яюЄхэЎшрыюь $Q_1$  Єръшь ётющёЄтюь эх юсырфрЄ№.

╩Ёрхт√х єёыютш  эрчютхь {\it эхт√Ёюцфхээ√ьш}, хёыш шч ЄЁхї ўшёхы
$J_{12}+J_{34}$ ш $J_{14}+J_{32}\pm i(J_{42}-J_{13})$ їюЄ  с√ фтр
ўшёыр юЄышўэ√ юЄ эєы . ▀ёэю, ўЄю т ёыєўрх тх∙хёЄтхээ√ї
ъю¤ЇЇшЎшхэЄют $u_{kj}$  яюэ Єшх Ёхуєы ЁэюёЄш ш эхт√ЁюцфхээюёЄш
ёютярфр■Є. ╠√ т√ ёэшь, ўЄю юяхЁрЄюЁ $L_Q$ ё юсырёЄ№■ \eqref{dom},
яюЁюцфхээ√щ эхт√Ёюцфхээ√ь ъЁрхт√ь єёыютшхь шьххЄ фшёъЁхЄэ√щ
ёяхъЄЁ, эю яюыэюх шчєўхэшх ь√ сєфхь яЁютюфшЄ№ Єюы№ъю фы 
Ёхуєы Ёэ√ї юяхЁрЄюЁют.

\section{└ёшьяЄюЄшўхёъшх ЇюЁьєы√ фы  ЇєэфрьхэЄры№эющ ёшёЄхь√ Ёх°хэшщ єЁртэхэш  $l_Q(y)=\lambda y$.}

╟фхё№ ь√ яюыєўшь  Ёхчєы№ЄрЄ√ юс рёшьяЄюЄшўхёъюь яютхфхэшш ЇєэфрьхэЄры№эющ ёшёЄхь√ Ёх°хэшщ ёшёЄхь√ $L_Q\mathbf{y}=\lambda\mathbf{y}$, ъюЄюЁ√х
ёє∙хёЄтхээю сєфєЄ шёяюы№чютрЄ№ё  фрыхх. ╬сючэрўшь ўхЁхч $\mathbf{s}(x,\lambda)$ Ёх°хэшх ё эрўры№э√ьш єёыютш ьш $s_1(0,\lambda)=0$, $s_2(0,\lambda)=1$
ш ўхЁхч $\mathbf{c}(x,\lambda)$
--- Ёх°хэшх ё эрўры№э√ьш єёыютш ьш $c_1(0,\lambda)=1$, $c_2(0,\lambda)=0$.
─юърчрЄхы№ёЄтю  рёшьяЄюЄшўхёъшї ЇюЁьєы фы  ¤Єшї Ёх°хэшщ яюЄЁхсєхЄ ёхЁ№хчэющ Єхїэшўхёъющ ЁрсюЄ√. ─рыхх ь√ сєфхь яЁютюфшЄ№ юЎхэъш т ъюьяыхъёэющ
$\lambda$--яыюёъюёЄш тэєЄЁш яюыюё
$$ \Pi_\alpha=\{\lambda\in\mathbb C\,\vert\ |\text{Im} \lambda|<\alpha\}.
$$
─ы  юЎхэъш юёЄрЄъют ь√ ттхфхь ЇєэъЎшш
\begin{gather*}
\upsilon_1 (x,\lambda) = \int\limits_0^xq_1(t)\sin(2\lambda t)\, dt, \qquad
\upsilon_2 (x,\lambda) = \int\limits_0^xq_1(t)\cos(2\lambda t)\, dt, \\
\notag \upsilon_3 (x,\lambda) =\frac12 \int\limits_0^x(q_2(t)+q_3(t))\sin(2\lambda t)\, dt, \qquad \upsilon_4(x,\lambda) = \frac12
\int\limits_0^x(q_2(t)+q_3(t))\cos(2\lambda t)\, dt.
\end{gather*}
═юЁьє ЇєэъЎшш $\upsilon_j(x,\lambda)$ яЁш ЇшъёшЁютрээюь $\lambda$ яю яхЁхьхээющ $x$ т яЁюёЄЁрэёЄтх $L_\nu[0,\pi]$ юсючэрўшь $v_{j,\nu}(\lambda)$.
╚эфхъё√ $p$ ш $p'$ тхчфх фрыхх сєфхь ёўшЄрЄ№ ЇшъёшЁютрээ√ьш: $p$ юяЁхфхы хЄё  єёыютшхь (ii), р $p'$  ты хЄё  ёюяЁ цхээ√ь шэфхъёюь ъ $p$, Є.х.
$1/p'+1/p=1$, яЁшўхь $p'=\infty$ яЁш $p=1$. ╬сючэрўшь Єръцх
\begin{equation}\label{gamma}
 \Upsilon(x,\lambda) =
\sum_{j=1}^4|\upsilon_j(x,\lambda)|,\qquad \Upsilon_\nu(\lambda)=\sum_{j=1}^4v_{j,\nu}(\lambda).
\end{equation}
╧Ёш ¤Єюь шэфхъё ЇєэъЎшш $\Upsilon_\infty(\lambda)$ фы  ъЁрЄъюёЄш сєфхь юяєёърЄ№, юсючэрўр  хх $\Upsilon(\lambda)$. ╟рьхЄшь ёЁрчє, ўЄю
$$
\Upsilon_\nu(\lambda)\le\pi^{1/\nu}\Upsilon(\lambda)\quad \text{ш} \quad\sup_{x\in[0,\pi]}|\Upsilon(x,\lambda)|\le\Upsilon(\lambda).
$$
╠√ т√сЁрыш шьхээю Єръшх
ЇєэъЎшш фы  юЎхэъш юёЄрЄюўэ√ї ўыхэют, яюёъюы№ъє  тэ√щ тшф ¤Єшї ЇєэъЎшщ яючтюышЄ эрь т фры№эхщ°хь яюыєўрЄ№ рёшьяЄюЄшўхёъшх ЇюЁьєы√ фы  ёюсёЄтхээ√ї
чэрўхэшщ ш ёюсёЄтхээ√ї ЇєэъЎшщ юяхЁрЄюЁют ─шЁрър т чртшёшьюёЄш юЄ яЁшэрфыхцэюёЄш яюЄхэЎшрыр $Q$ ъ Єюьє шыш шэюьє ъырёёє ЇєэъЎшщ. ╦хуъю тшфхЄ№ (¤Єю
ёЁрчє цх ёыхфєхЄ шч ыхьь√ ╨шьрэр--╦хсхур), ўЄю яЁш ърцфюь ЇшъёшЁютрээюь $x\in[0,\pi]$ ЇєэъЎш  $\Upsilon(x,\lambda)$ ёЄЁхьшЄё  ъ эєы■ яЁш
$\lambda\to\pm\infty$. ─юърцхь, ўЄю ¤Єю ёЄЁхьыхэшх ъ эєы■ ЁртэюьхЁэю яю $x\in[0,\pi]$. ┬ фры№эхщ°шї юЎхэърї сєфхь шёяюы№чютрЄ№ ёыхфє■∙шх юўхтшфэ√х
эхЁртхэёЄтр. ╧єёЄ№ яЁюшчтюы№э√х Єюўъш $z$, $z_1$ ш $z_2$ ыхцрЄ тэєЄЁш ъЁєур $|z|<\xi\le1$, р Єюўър $w$ ыхцшЄ тэєЄЁш яюыюё√ $\Pi_\alpha$ т ъюьяыхъёэющ
яыюёъюёЄш. ╥юуфр
\begin{gather}
 |\sin w|\le\ch\alpha,\qquad|\cos w|\le\ch\alpha, \qquad |\sin
 z|\le|z|\ch\xi,\notag
 \\ |\cos z-1|\le\frac12|z|^2\ch\xi,\qquad|\sin
z-z|\le\frac16|z|^2\ch\xi,\label{complexineq}\\ |\cos z_2-\cos z_1|\le|z_2-z_1|\ch\xi, \qquad |(\sin z_2-z_2)-(\sin
z_1-z_1)|\le\frac12|z_2-z_1|\xi^2\ch\xi,\notag\\
|\exp(z)-1|\le|z|\exp(\xi),\qquad|\exp(z)-1-z|\le\frac12|z|^2\exp(\xi).\notag
\end{gather}
\begin{Lemma}\label{lem:Ups1}
╧єёЄ№ $Q\in L_1[0,\pi]$ ш т√яюыэхэ√ єёыютш  (i), (ii), (iii), т√яшёрээ√х т√°х. ╧єёЄ№ $\alpha>0$ --- яЁюшчтюы№эюх ЇшъёшЁютрээюх ўшёыю,
$\Pi_\alpha=\{\lambda\in\mathbb C\,\vert\ |\text{Im} \lambda|<\alpha\}$.  ╥юуфр $\Upsilon(\lambda)\to0$ яЁш $|\lambda|\to\infty$ ЁртэюьхЁэю т яюыюёх
$\Pi_\alpha$.
\end{Lemma}
\begin{proof}
▀ёэю, ўЄю фюърчрЄхы№ёЄтю фюёЄрЄюўэю яЁютхёЄш юЄфхы№эю фы  ърцфющ шч ЇєэъЎшщ $\upsilon_j$. ╠√ ЁрёёьюЄЁшь Єюы№ъю ЇєэъЎш■ $\upsilon_1$
--- юёЄры№э√х ёырурхь√х юЎхэштр■Єё  рэрыюушўэю. ╧юфсхЁхь
эхяЁхЁ√тэю фшЇЇхЁхэЎшЁєхьє■ ЇєэъЎш■ $p_1(x)$ Єръ, ўЄю
$$
\|q_1-p_1\|_{L_1[0,\pi]}<\varepsilon/(2\ch(2\pi\alpha)).
$$
╥юуфр
$$
\left|\int_0^x(q_1(t)-p_1(t))\sin(2\lambda t)dt\right|\le\ch(2\pi\alpha)\|q_1-p_1\|_{L_1[0,\pi]}<\varepsilon/2.
$$
─рыхх, яЁш $\lambda\to\infty$, $\lambda\in\Pi_\alpha$,
$$
\int_0^xp_1(t)\sin(2\lambda t)dt=\frac{p_1(0)-p_1(x)\cos(2\lambda x)}{2\lambda}+\frac1{2\lambda}\int_0^xp_1'(t)\cos(2\lambda
t)dt=O\left(\frac1{\lambda}\right),
$$
уфх юЎхэъш т√яюыэхэ√ ЁртэюьхЁэю яю $x\in[0,\pi]$. ═ю Єюуфр $|\upsilon_1(x,\lambda)|<\varepsilon$ яЁш $\lambda\in\Pi_\alpha$, $|\lambda|>\mu$, яЁш
фюёЄрЄюўэю сюы№°юь $\mu$.
\end{proof}
╬ЄьхЄшь Єръцх трцэє■ фы  фры№эхщ°хую юЎхэъє
\begin{Lemma}
\begin{equation}\label{UpsUps}
\sup\limits_{x\in[0,\pi]}|\Upsilon(x,\lambda)|^2\le 8R\ch(2\pi\alpha)\Upsilon_\nu(\lambda).
\end{equation}
\end{Lemma}
\begin{proof}
╘єэъЎш  $\Upsilon(x,\lambda)^2$ ёюёЄюшЄ шч ёєьь√ °хёЄэрфЎрЄш яюярЁэ√ї яЁюшчтхфхэшщ тшфр $\upsilon_k(x,\lambda)\upsilon_j(x,\lambda)$,
$k,\,j=1,2,3,4$. ╧ю¤Єюьє фюёЄрЄюўэю юЎхэшЄ№ юфэю Єръюх ёырурхьюх. ╧Ёшьхэшт эхЁртхэёЄтю ├хы№фхЁр, яюыєўшь
$$
|\upsilon_k(x,\lambda)\upsilon_j(x,\lambda)|=\left|\int\limits_0^x(\upsilon_k'(t,\lambda)\upsilon_j(t,\lambda)+
\upsilon_k(t,\lambda)\upsilon_j'(t,\lambda)),dt\right|\le R\ch(2\pi\alpha)(\upsilon_{j,\nu}(\lambda)+\upsilon_{k,\nu}(\lambda)).
$$
╤ъырф√тр  ¤Єш эхЁртхэёЄтр ш єўшЄ√тр , ўЄю ёрьр ЇєэъЎш  $\Upsilon_\nu(\lambda)$ ёюёЄюшЄ шч ўхЄ√Ёхї ёырурхь√ї, яюыєўрхь юЎхэъє \eqref{UpsUps}.
\end{proof}

╠√ эрўэхь ё рёшьяЄюЄшъш ЇєэъЎшш $\mathbf{s}(x,\lambda)$. ┬ ёшёЄхьх $l_Q\mathbf{y}=\lambda\mathbf{y}$ ёфхырхь чрьхэє ╧Ё■ЇхЁр (ёь. \cite{Ha})
\begin{equation}\label{Pr}
y_1(x,\lambda)=-r(x,\lambda)\sin\theta(x,\lambda),\ y_2(x,\lambda)=r(x,\lambda)\cos\theta(x,\lambda)
\end{equation}
(хх ьюцэю ЄЁръЄютрЄ№ ъръ яхЁхїюф ъ яюы Ёэ√ь ъююЁфшэрЄрь). ╥юуфр ёшёЄхьє ьюцэю чряшёрЄ№ т тшфх
\begin{equation}
\begin{array}{ccc}
r'\cos\theta-r\theta'\sin\theta-q_1r\sin\theta+q_2r\cos\theta &=&-\lambda r\sin\theta,\\
r'\sin\theta+r\theta'\cos\theta-q_3r\sin\theta-q_1r\cos\theta &=&\lambda r\cos\theta,
\end{array}
\label{polar}
\end{equation}
уфх $r=r(x,\lambda)$, $\theta=\theta(x,\lambda)$, $q_j=q_j(x)$, $j=1,2,3$, р яЁюшчтюфэ√х ЇєэъЎшщ $r$ ш $\theta$ схЁєЄё  яю яхЁхьхээющ $x$. ╙ьэюцшь
яхЁтюх єЁртэхэшх т~\eqref{polar} эр $-\sin\theta$ ш яЁшсртшь тЄюЁюх єЁртэхэшх, єьэюцхээюх эр $\cos\theta$. ┬ Ёхчєы№ЄрЄх яюыєўшь єЁртэхэшх фы  ЇєэъЎшш
$\theta(x,\lambda)$
\begin{equation}
\theta'(x,\lambda)-q_1(x)\cos2\theta(x,\lambda)- \frac12(q_2(x)+q_3(x))\sin2\theta(x,\lambda)=\lambda. \label{th}
\end{equation}
═рўры№эюх єёыютшх яюырурхь $\theta(0,\lambda)=0$, Єръ ъръ $s_1(0,\lambda)=0$. ┼ёыш ь√ єьэюцшь яхЁтюх єЁртэхэшх т~\eqref{polar} эр $\cos\theta$ ш
яЁшсртшь тЄюЁюх єЁртэхэшх, єьэюцхээюх эр $\sin\theta$, Єю яюыєўшь єЁртэхэшх эр ЇєэъЎш■ $r(x,\lambda)$
\begin{equation}
r'(x,\lambda)=r(x,\lambda)\left[\frac12(q_3(x)-q_2(x))+q_1(x)\sin 2\theta(x,\lambda)-\frac12(q_2(x)+q_3(x)) \cos 2\theta(x,\lambda)\right] \label{r}
\end{equation}
ё эрўры№э√ь єёыютшхь $r(0,\lambda)=1$. ╥ръшь юсЁрчюь юЄ ёшёЄхь√ $l_Q(y)=\lambda y$ ь√ яхЁх°ыш ъ ёшёЄхьх фтєї єЁртэхэшщ~\eqref{th}, \eqref{r}. ▌Єш
єЁртэхэш  эх  ты ■Єё  ышэхщэ√ьш. ═ю фюёЄюшэёЄтю эютющ ёшёЄхь√ ёюёЄюшЄ т Єюь, ўЄю єЁртэхэшх~\eqref{th} эх ёюфхЁцшЄ эхшчтхёЄэющ ЇєэъЎшш $r(x,\lambda)$
ш  ты хЄё  эхчртшёшь√ь фшЇЇхЁхэЎшры№э√ь єЁртэхэшхь эр ЇєэъЎш■ $\theta(x,\lambda)$. ═р°р сышцрщ°р  Ўхы№
--- эрщЄш рёшьяЄюЄшўхёъшх яЁхфёЄртыхэш  фы  ЇєэъЎшщ
$\theta(x,\lambda)$ ш $r(x,\lambda)$. ▌Єш ЇюЁьєы√ сєфєЄ яюыєўхэ√ т эшцхёыхфє■∙шї ыхььрї. ╧хЁтр  шч ¤Єшї ыхьь  ты хЄё  ъы■ўхтющ. ═р эхщ срчшЁє■Єё 
юёэютэ√х Ёхчєы№ЄрЄ√ ЁрсюЄ√. ─рыхх схч эряюьшэрэшщ шёяюы№чєхь ттхфхээ√х т√°х юсючэрўхэш .
\begin{Lemma}\label{lem:2.1}
╧єёЄ№ $Q\in L_1[0,\pi]$ ш т√яюыэхэ√ єёыютш  (i), (ii), (iii), т√яшёрээ√х т√°х. ╧єёЄ№ $\alpha>0$ --- яЁюшчтюы№эюх ЇшъёшЁютрээюх ўшёыю,
$\Pi_\alpha=\{\lambda\in\mathbb C\,\vert\ |\text{Im} \lambda|<\alpha\}$.  ╥юуфр яЁш ы■сюь $\lambda$ шч юсырёЄш
$$
D_{Q,\alpha}=\left\{\lambda\in\mathbb C: \ |Im \lambda|<\alpha,\,\Upsilon(\lambda)<\frac1{8 k^4}\right\}, \quad \text{уфх} \ \,
k=2+12R\ch(2\pi\alpha+1),
$$
єЁртэхэшх \eqref{th} шьххЄ хфшэёЄтхээюх Ёх°хэшх $\theta(x,\lambda)$, юяЁхфхыхээюх яЁш тёхї $0\le x\le\pi$ ш єфютыхЄтюЁ ■∙хх эрўры№эюьє єёыютш■
$\theta(0,\lambda)=0$. ▌Єю Ёх°хэшх фюяєёърхЄ яЁхфёЄртыхэшх
\begin{equation}\label{null}
\theta(x,\lambda)=\lambda x+\eta(x,\lambda),
\end{equation}
уфх
$$
|\eta(x,\lambda)|\le\Upsilon(x,\lambda)+k^2\Upsilon_{p'}(\lambda)\qquad\text{ш}\quad \|\eta(x,\lambda)\|_{p'}\le(1+\pi k^2)\Upsilon_{p'}(\lambda).
$$
┴юыхх Єюую,
\begin{equation}\label{first}
\eta(x,\lambda)=\int\limits_0^xq_1(t)\cos(2\lambda t)\, dt+\frac12\int\limits_0^x(q_2(t)+q_3(t))\sin(2\lambda t)dt+\zeta(x,\lambda),
\end{equation}
уфх
$$
\|\zeta(x,\lambda)\|_C\le k^2\Upsilon_{p'}(\lambda).
$$
╧Ёш ¤Єюь
\begin{equation}\label{null1}
|\eta'(x,\lambda)|\le\ch(2\pi\alpha+1)\left(|q_1(x)|+\frac12|q_2(x)+q_3(x)|\right)
\end{equation}
яюўЄш тё■фє эр $x\in[0,\pi]$.
\end{Lemma}
\begin{proof}  \textit{╪ру 1.} ╧Ёхцфх тёхую чрьхЄшь, ўЄю юсырёЄ№ $\Pi_\alpha \setminus D_{Q,\alpha}$ юуЁрэшўхэр, Є.х.
ёє∙хёЄтєхЄ  ўшёыю $\mu=\mu(Q, \alpha)$,  Єръюх ўЄю тёх ўшёыр $\lambda\in \Pi_\alpha$ ыхцрЄ т $D_{Q,\alpha}$,  хёыш Єюы№ъю $|\text{Re} \lambda|
\geqslant \mu$. ▌Єю єЄтхЁцфхэшх т√ЄхърхЄ шч ╦хьь√ \ref{lem:Ups1}.

╧хЁхяш°хь єЁртэхэшх \eqref{th} т шэЄхуЁры№эюь тшфх
\begin{equation*}
\theta(x,\lambda)=\lambda x+\int\limits_0^x q_1(t)\cos2\theta(t,\lambda)dt+\frac12\int\limits_0^x(q_2(t)+q_3(t))\sin2\theta(t,\lambda)dt.
\end{equation*}
ш ёфхырхь чрьхэє $\eta(x,\lambda)=\theta(x,\lambda)-\lambda x$. ╥юуфр єЁртэхэшх яЁшьхЄ тшф
\begin{equation}
\label{integr} \eta(x,\lambda)=\int\limits_0^xq_1(t)\cos(2\eta(t,\lambda)+2\lambda
t)dt+\frac12\int\limits_0^x(q_2(t)+q_3(t))\sin(2\eta(t,\lambda)+2\lambda t)dt.
\end{equation}
╧Ёртє■ ўрёЄ№ ¤Єюую єЁртэхэш  ь√ юсючэрўшь ўхЁхч $F_{Q,\lambda}(\eta)$, уфх $F_{Q,\lambda}$ --- эхышэхщэюх юЄюсЁрцхэшх яЁюёЄЁрэёЄтр $C[0,\pi]$ т ёхс ,
$$
F_{Q,\lambda}(\eta)=\int\limits_0^xq_1(t)\cos(2\eta(t,\lambda)+2\lambda t)dt+\frac12\int\limits_0^x(q_2(t)+q_3(t))\sin(2\eta(t,\lambda)+2\lambda
t)dt.
$$
═ряюьэшь ЇюЁьєышЁютъє  яЁшэЎшяр  ёцшьр■∙шї юЄюсЁрцхэшщ, ъюЄюЁ√ь ь√ эрьхЁхтрхьё  тюёяюы№чютрЄ№ё . {\sl ╧єёЄ№ юЄюсЁрцхэшх $\Phi$ яхЁхтюфшЄ чрьъэєЄ√щ
°рЁ $B(y,r)$ срэрїютр яЁюёЄЁрэёЄтр $X$ т ёхс  ш фы  эхъюЄюЁюую $q\in[0,1)$ ш $\forall y_1,\,y_2\in B(y,r)$ т√яюыэхэю $\|\Phi(y_2)-\Phi(y_1)\|\le
q\|y_2-y_1\|$. ╥юуфр т °рЁх $B(y,r)$ ёє∙хёЄтєхЄ хфшэёЄтхээр  эхяюфтшцэр  Єюўър $x$ (Є.х. $\Phi(x)=x$). ╧Ёш ¤Єюь $x=\lim x_n$, уфх $x_0\in B(y,r)$ ---
яЁюшчтюы№эю, р $x_{n}=\Phi(x_{n-1})$ фы  тёхї $n=1,\,2,\,3,\dots$.}

┬ ърўхёЄтх ЎхэЄЁр °рЁр $B(y,r)$ яюыюцшь
$$
f_0=\int\limits_0^xq_1(t)\cos(2\lambda t)dt+\frac12\int\limits_0^x(q_2(t)+q_3(t))\sin(2\lambda t)dt,
$$
р Ёрфшєё т√схЁхь Ёртэ√ь
$$
r=k^2\Upsilon_{p'}(\lambda)\le k^2\pi\Upsilon(\lambda)<\pi/(8k^2)<1/10.
$$

\textit{╪ру 2.} ╬Ўхэшь ЁрчэюёЄ№ $\|F(f_0)-f_0\|_{C[0,\pi]}$.
\\╧ЁхфёЄртшь юЄюсЁрцхэшх $F$ ёыхфє■∙шь юсЁрчюь:
\begin{gather*}
F(\eta)=f_0+F_1\eta+\mathcal F(\eta),\quad\text{уфх}\\
F_1\eta=-2\int\limits_0^xq_1(t)\sin(2\lambda t)\eta(t, \lambda)dt+\int\limits_0^x(q_2(t)+q_3(t))\cos(2\lambda t)\eta(t, \lambda)dt,\\
\mathcal F(\eta)=\int\limits_0^xq_1(t)\left[\cos(2\eta(t, \lambda)+2\lambda t)-\cos(2\lambda t)+2\sin(2\lambda t)\eta(t,
\lambda)\right]dt+\\+\frac12\int\limits_0^x(q_2(t)+q_3(t))\left[\sin(2\eta(t, \lambda)+2\lambda t)-\sin(2\lambda t)-2\cos(2\lambda t)\eta(t,
\lambda)\right]dt.
\end{gather*}
╥юуфр $F(f_0)=f_0+F_1(f_0)+\mathcal F(f_0)$. ┬эрўрых юЄьхЄшь, ўЄю
\begin{equation}
\label{F_0est} |f_0(x)|\le\Upsilon(x,\lambda),\qquad\|f_0\|_{p'}\le\Upsilon_{p'}(\lambda),
\qquad\|f_0\|_C\le\Upsilon(\lambda)<\frac1{8k^4}<\frac1{128}.
\end{equation}
╬Ўхэшь ЄхяхЁ№ эюЁьє тхъЄюЁр
\begin{multline*}
F_1(f_0)=-2\int\limits_0^xq_1(t)\sin(2\lambda t)\int\limits_0^tq_1(s)\cos(2\lambda s)dsdt-\\-\int\limits_0^xq_1(t)\sin(2\lambda
t)\int\limits_0^t(q_2(s)+q_3(s))\sin(2\lambda s)dsdt+\label{F1F0}\\
+\int\limits_0^x(q_2(t)+q_3(t))\cos(2\lambda t)\int\limits_0^tq_1(s)\cos(2\lambda s)dsdt+\\+\frac12\int\limits_0^x(q_2(t)+q_3(t))\cos(2\lambda
t)\int\limits_0^t(q_2(s)+q_3(s))\sin(2\lambda s)dsdt.
\end{multline*}
╧Ёютхфхь юЎхэъє Єюы№ъю яхЁтюую шч ўхЄ√Ёхї яютЄюЁэ√ї шэЄхуЁрыют --- юёЄры№э√х ЄЁш юЎхэштр■Єё  рэрыюушўэю
$$
\left|\int\limits_0^xq_1(t)\sin(2\lambda t)\int\limits_0^tq_1(s)\cos(2\lambda s)dsdt\right|\le \|q_1(t)\sin(2\lambda
t)\|_p\cdot\upsilon_{2,p'}(\lambda)\le R\ch(2\pi\alpha)\cdot\upsilon_{2,p'}(\lambda),
$$
уфх ЇєэъЎш  $\upsilon_{2,p'}$ юяЁхфхыхэр т эрўрых ярЁруЁрЇр. ┬ Ёхчєы№ЄрЄх яюыєўшь юЎхэъє
$$
\|F_1(f_0)\|_C\le 2R\ch(2\pi\alpha)\Upsilon_{p'}(\lambda).
$$

╧хЁхщфхь ъ юЎхэъх эюЁь√ тхъЄюЁр $\mathcal F(f_0)$. ╧Ёхцфх тёхую, чрьхЄшь, ўЄю яЁш єёыютшш $|f_0(t)|<1/{128}$ шч эхЁртхэёЄт \eqref{complexineq} ш
ЄЁшуюэюьхЄЁшўхёъшї ЇюЁьєы ёыхфєхЄ
\begin{gather}
|\cos(2f_0(t)+2\lambda t)-\cos(2\lambda t)+2\sin(2\lambda t)f_0(t)|\le\notag\\
\le\frac83\ch(2\pi\alpha)\left(\ch\frac1{64}\right)|f_0(t)|^2\le\frac{\ch(2\pi\alpha)}2|f_0(t)|,\label{calF}\\
|\sin(2f_0(t)+2\lambda t)-\sin(2\lambda t)-2\cos(2\lambda t)f_0(t)|\le \frac{\ch(2\pi\alpha)}2|f_0(t)|.\notag
\end{gather}
╬Єё■фр т√ЄхърхЄ юЎхэър
$$
\|\mathcal F(f_0)\|_C\le\left(\|q_1\|_p+\frac12\|q_2+q_3\|_p\right)\frac12\ch(2\pi\alpha)\|f_0\|_{p'}\le R\ch(2\pi\alpha)\Upsilon_{p'}(\lambda).
$$
╚Єръ, ь√ яюърчрыш, ўЄю
\begin{equation}\label{St1res}
\|F(f_0)-f_0\|_C\le 3R\ch(2\pi\alpha)\Upsilon_{p'}(\lambda)<\frac{k}4\Upsilon_{p'}(\lambda).
\end{equation}
╪ру 2 чртхЁ°хэ.

\textit{╪ру 3.} ═рь эх єфрыюё№ фюърчрЄ№, ўЄю юЄюсЁрцхэшх $F$ ёцшьрхЄ. ─ры№эхщ°р  шфх  фюърчрЄхы№ёЄтр ёюёЄюшЄ т ёыхфєъ■∙хь: ь√ яюърцхь, ўЄю
юЄюсЁрцхэшх $\Phi$, юяЁхфхыхээюх т \eqref{Phi} єцх  ты хЄё  ёцшьр■∙шь. ─ы  ¤Єюую ь√ Ёрчыюцшь юЄюсЁрцхэшх $F$ т Єюўъх $f_0$:
\begin{gather}\label{Fdecomp}
F(\eta)=F(f_0)+G_1(\zeta)+\mathcal G(\zeta),\qquad\text{уфх}\
\zeta=\eta-f_0,\\
G_1(\zeta)=-2\int\limits_0^xq_1(t)\sin(2\lambda t+2f_0(t))\zeta(t)dt+\int\limits_0^x(q_2(t)+q_3(t))\cos(2\lambda t+2f_0(t))\zeta(t)dt,\notag\\
\mathcal G(\zeta)=\int\limits_0^xq_1(t)\left[\cos(2\zeta(t)+2\lambda t+2f_0(t))-\cos(2\lambda t+2f_0(t))+2\sin(2\lambda
t+2f_0(t))\zeta(t)\right]dt+\notag\\+\frac12\int\limits_0^x(q_2(t)+q_3(t))\left[\sin(2\zeta(t)+2\lambda t+2f_0(t))-\sin(2\lambda
t+2f_0(t))-2\cos(2\lambda t+2f_0(t))\zeta(t)\right]dt.\notag
\end{gather}
╬Ўхэшь эюЁь√ юЄюсЁрцхэшщ $G_1$, $G_1^2$ ш $\mathcal G$. ╦хуъю тшфхЄ№, ўЄю яЁш єёыютшш $\|f_0\|_C\le1/2$, ъюЄюЁюх т√яюыэхэю т $D_{Q,\alpha}$ т ёшыє
юЎхэъш \eqref{F_0est}, шьххь
$$
\|G_1\|_C\le3\pi R\ch(2\pi\alpha+1).
$$
╩юэхўэю, эюЁьр $G_1^2$ юЎхэштрхЄё  ътрфЁрЄюь яЁртющ ўрёЄш ¤Єюую эхЁртхэёЄтр, эю эрь яюЄЁхсєхЄё  сюыхх Єюўэр  юЎхэър. ╦хуъю тшфхЄ№, ўЄю
\begin{gather*}
G_1^2\zeta=4\int\limits_0^xq_1(t)\sin(2\lambda t+2f_0(t))\int\limits_0^tq_1(s)\sin(2\lambda s+2f_0(s))\zeta(s)\,ds\,dt-\\
-2\int\limits_0^xq_1(t)\sin(2\lambda
t+2f_0(t))\int\limits_0^t(q_2(s)+q_3(s))\cos(2\lambda s+2f_0(s))\zeta(s)\,ds\,dt-\notag\\
-2\int\limits_0^x(q_2(t)+q_3(t))\cos(2\lambda t+2f_0(t))\int\limits_0^tq_1(s)\sin(2\lambda s+2f_0(s))\zeta(s)\,ds\,dt+\notag\\
+\int\limits_0^x(q_2(t)+q_3(t))\cos(2\lambda t+2f_0(t))\int\limits_0^t(q_2(s)+q_3(s))\cos(2\lambda s+2f_0(s))\zeta(s)\,ds\,dt.\notag
\end{gather*}
─ы  юяЁхфхыхээюёЄш юЎхэшь яхЁт√щ шэЄхуЁры (юёЄры№э√х юЎхэштр■Єё  рэрыюушўэю)
\begin{gather}
\left|\int\limits_0^xq_1(t)\sin(2\lambda t+2f_0(t))\int\limits_0^tq_1(s)\sin(2\lambda
s+2f_0(s))\zeta(s)\,ds\,dt\right|=\notag\\
=\left|\int\limits_0^xq_1(s)\sin(2\lambda s+2f_0(s))\zeta(s)\int\limits_s^xq_1(t)\sin(2\lambda
t+2f_0(t))\,dt\,ds\right|\le\label{G1sq2}\\
\|\zeta\|_CR\ch(2\pi\alpha+1)\left(\left\|\int\limits_s^xq_1(t)\cos(2f_0(t))\sin(2\lambda
t)\,dt\right\|_{p'}+\left\|\int\limits_s^xq_1(t)\sin(2f_0(t))\cos(2\lambda t)\,dt\right\|_{p'}\right)\notag.
\end{gather}
╧хЁт√щ шэЄхуЁры т яЁртющ ўрёЄш яюёыхфэхую эхЁртхэёЄтр яЁхфёЄртшь т тшфх
$$
\int\limits_s^xq_1(t)\sin(2\lambda t)\,dt+\int\limits_s^xq_1(t)(\cos(2f_0(t))-1)\sin(2\lambda t)dt.
$$
╥юуфр $\|\cdot\|_{p'}$ эюЁьр яхЁтюую ёырурхьюую эх яЁхтюёїюфшЄ $2\upsilon_{1,p'}(\lambda)$, р ьюфєы№ тЄюЁюую шэЄхуЁрыр (ёь. \eqref{complexineq}) эх
яЁхтюёїюфшЄ
$$
2R\ch(2\pi\alpha)\ch(1/64)\||f_0|^2\|_{p'}\le(1/32)R\ch(2\pi\alpha)\Upsilon_{p'}(\lambda).
$$
┬ЄюЁющ шэЄхуЁры т яЁртющ ўрёЄш \eqref{G1sq2} рэрыюушўэю ЁрчсштрхЄё  т ёєььє
$$
2\int\limits_s^xq_1(t)\cos(2\lambda t)f_0(t)\,dt+\int\limits_s^xq_1(t)(\sin(2f_0(t))-2f_0(t))\cos(2\lambda t)dt
$$
ш хую $\|\cdot\|_{p'}$ эюЁьр эх яЁхтюёїюфшЄ
$$
2R\ch(2\pi\alpha)\Upsilon_{p'}(\lambda)+\frac1{128}R\ch(2\pi\alpha)\Upsilon_{p'}(\lambda).
$$
╤тюф  тьхёЄх тёх ¤Єш юЎхэъш яюыєўшь
\begin{gather*}
\left|\int\limits_0^xq_1(t)\sin(2\lambda t+2f_0(t))\int\limits_0^tq_1(s)\sin(2\lambda
s+2f_0(s))\zeta(s)\,ds\,dt\right|\le\\
\le\|\zeta\|_CR\ch(2\pi\alpha+1)(2\upsilon_{1,p'}(\lambda)+3R\ch(2\pi\alpha)\Upsilon_{p'}(\lambda)).
\end{gather*}
╬Єё■фр ёыхфєхЄ юЎхэър
$$
\|G_1^2\|_C\le R\ch(2\pi\alpha+1)(12+27R\ch(2\pi\alpha))\Upsilon_{p'}(\lambda)<(k^2/2)\Upsilon_{p'}(\lambda).
$$
╧хЁхщфхь ъ юЎхэърь юЄюсЁрцхэш  $\mathcal G$. └эрыюушўэю эхЁртхэёЄтрь \eqref{calF}, т °рЁх $\|\zeta\|_C\le r$ юЎхэшь
\begin{gather*}
\left|\cos(2\zeta(t)+2\lambda t+2f_0(t))-\cos(2\lambda t+2f_0(t))+2\sin(2\lambda
t+2f_0(t))\zeta(t)\right|\ \text{ш}\\
\left|\sin(2\zeta(t)+2\lambda t+2f_0(t))-\sin(2\lambda t+2f_0(t))-2\cos(2\lambda t+2f_0(t))\zeta(t)\right|
\end{gather*}
тхышўшэющ $(8/3)\ch(2\pi\alpha+1)\ch(2r)|\zeta(t)|^2$. ╥юуфр
$$
\|\mathcal G(\zeta)\|_C\le\frac{16\pi}3R\ch(2\pi\alpha+1)\ch(2r)\|\zeta\|_C^2.
$$
╧юёъюы№ъє $r<1/10$, Єю ўшёыю $16\pi\ch(2r)/3$ эх яЁхтюёїюфшЄ $18$ ш Єюуфр
\begin{equation}\label{calG1}
\|\mathcal G(\zeta)\|_C\le\frac{3k}2\|\zeta\|_C^2.
\end{equation}
═рь яюЄЁхсєхЄё  Єръцх юЎхэър эюЁь√ ЁрчэюёЄш $\|\mathcal G(\zeta_1)-\mathcal G(\zeta_2)\|_C$. ┬эют№ юсЁрЄшьё  ъ эхЁртхэёЄтрь \eqref{complexineq}:
\begin{gather*}
\left|\cos(2\zeta_1(t)+2\lambda t+2f_0(t))-\cos(2\zeta_2+2\lambda t+2f_0(t))+2\sin(2\lambda
t+2f_0(t))(\zeta_1(t)-\zeta_2(t))\right|\le\\
\le|\cos(2\lambda t+2f_0(t))|\,|\cos(2\zeta_1(t))-\cos(2\zeta_2(t))|+\\
+|\sin(2\lambda t+2f_0(t))|\,\left|\left(\sin(2\zeta_1(t)-2\zeta_1(t)\right)-
\left(\sin(2\zeta_2(t)-2\zeta_2(t)\right)\right|\le\\
\le\frac52\ch(2\pi\alpha+1)r\ch r|\zeta_2(t)-\zeta_1(t)|.
\end{gather*}
╬Єё■фр яюыєўрхь юЎхэъє
$$
\|\mathcal G(\zeta_1)-\mathcal G(\zeta_2)\|_C\le5\pi R\ch(2\pi\alpha+1)r\ch r\|\zeta_2-\zeta_1\|_C,
$$
ёяЁртхфыштє■ т °рЁх $B(f_0,r)$. ╧юёъюы№ъє $r<1/10$, Єю ўшёыю $5\pi\ch r$ эх яЁхтюёїюфшЄ $16$ ш Єюуфр
\begin{equation*}
\|\mathcal G(\zeta_1)-\mathcal G(\zeta_2)\|_C\le 2kr\|\zeta_2-\zeta_1\|_C.
\end{equation*}
╪ру 3 чртхЁ°хэ.

\textit{╪ру 4. ─юърчрЄхы№ёЄтю юЎхэюъ \eqref{null}, \eqref{null1} ш \eqref{first}.} ╚ёяюы№чє  Ёрчыюцхэшх \eqref{Fdecomp}, чряш°хь єЁртэхэшх
\eqref{integr} т тшфх
\begin{equation*}
\zeta=F(f_0)-f_0+G_1(\zeta)+\mathcal G(\zeta),\qquad \text{уфх}\ \zeta=\eta-f_0,
\end{equation*}
ъюЄюЁюх ь√ т ётю■ юўхЁхф№ яЁхюсЁрчєхь ъ тшфє
\begin{equation}\label{Phi}
\eta=f_0+(I-G_1)^{-1}(F(f_0)-f_0+\mathcal G(\zeta)):=\Phi(\eta).
\end{equation}
╧юърцхь, ўЄю юЄюсЁрцхэшх $\Phi =\Phi_{Q,\lambda}$ (чфхё№ $\lambda\in D_{Q,\alpha}$) ъюЁЁхъЄэю юяЁхфхыхэю,  ты хЄё  ёцшьр■∙шь т чрьъэєЄюь °рЁх
$B(f_0,r)$ яЁюёЄЁрэёЄтр $C[0,\pi]$ ш яхЁхтюфшЄ хую т ёхс . ┬ ╪рух 3 ь√ юЎхэшыш эюЁьє юяхЁрЄюЁют $G_1$ ш $G_1^2$
$$
\|G_1\|_C\le3\pi R\ch(2\pi\alpha+1),\qquad \|G_1^2\|_C\le\frac{k^2}2\Upsilon_{p'}(\lambda)<\frac1{4k^2}.
$$
╥юуфр яЁш Ўхы√ї $s\geqslant 1$ ш $\lambda \in D_{Q,\alpha}$ ёяЁртхфышт√ юЎхэъш
$$
\|G_1^{2s}\|_C\le  \left(\frac 1{4k^2}\right)^s
 <\left(\frac 1{16}\right)^s, \quad \|G_1^{2s+1}\|_C\le
 \frac{\pi k}4 \left(\frac 1{4k^2}\right)^s<
 \frac{\pi}2\left( \frac1{16}\right)^s.
 $$
╤ыхфютрЄхы№эю, юяхЁрЄюЁ $(I-G_1)^{-1}= I+G_1+G_1^2+ \dots$ ёє∙хёЄтєхЄ, яЁшўхь хую эюЁьр фюяєёърхЄ юЎхэъє
\begin{equation} \label{norm}
\|(I-G_1)^{-1}\|_C \leqslant \sum_{j=0}^\infty\|G_1^j\|_C\le 1+\|G_1\|_C +1\le 2 +3\pi R\ch(2\pi\alpha+1)<k.
\end{equation}
╧хЁхщфхь ъ юЎхэърь юЄюсЁрцхэш  $\Phi$, юяЁхфхыхээюую т \eqref{Phi}. ─ы  ы■сюую $\eta=f_0+\zeta$, уфх $\|\zeta\|_C\le r$ шч эхЁртхэёЄт \eqref{St1res},
\eqref{calG1} ш \eqref{norm} ёыхфєхЄ
\begin{equation*}
\|\Phi(\eta)-f_0\|_C\le\|(I-G_1)^{-1}\|_C\left(\|F(f_0)-f_0\|_C+\|\mathcal G(\zeta)\|_C\right) \le
k\left(\frac{k}4\Upsilon_{p'}(\lambda)+\frac{3k}2r^2\right)\le r,
\end{equation*}
яюёъюы№ъє $k^2\Upsilon_{p'}=r$, р
$$
k^2r^2=k^4\Upsilon_{p'}(\lambda)r\le\pi k^4\Upsilon(\lambda)r\le(\pi/8)r<r/2.
$$
╥ръшь юсЁрчюь, ь√ фюърчрыш, ўЄю юЄюсЁрцхэшх $\Phi$ яхЁхтюфшЄ чрьъэєЄ√щ °рЁ $B(f_0,r)$ т ёхс . ─рыхх
\begin{equation*}
\|\Phi(\zeta_2)-\Phi(\zeta_1)\|_C\le\|(I-G_1)^{-1}\|_C\|\mathcal G(\zeta_2)-\mathcal G(\zeta_1)\|_C \le 2k^2r\|\zeta_2-\zeta_1\|_C.
\end{equation*}
╧юёъюы№ъє
$$
2k^2r=2k^4\Upsilon_{p'}(\lambda)<\pi/4<1,
$$
юЄюсЁрцхэшх $\Phi$  ты хЄё  ёцшьр■∙шь т °рЁх $B(f_0,r)$ ё ъю¤ЇЇшЎшхэЄюь ёцрЄш  $q=2k^4\Upsilon_{p'}(\lambda)$. ╧Ёшьхэ   ъ $\Phi$ яЁшэЎшя ёцшьр■∙шї
юЄюсЁрцхэшщ, яюыєўшь ёє∙хёЄтютрэшх Ёх°хэш  єЁртэхэш  \eqref{th} тшфр $\theta(x,\lambda)=\lambda x+\eta(x,\lambda)$, уфх ЇєэъЎш  $\eta$ шьххЄ тшф
$\eta=f_0+\zeta$, р $\zeta$ фюяєёърхЄ юЎхэъє $\|\zeta\|_C\le r=k^2\Upsilon_{p'}(\lambda)$. ╥ръшь юсЁрчюь, ь√ фюърчрыш рёшьяЄюЄшўхёъюх яЁхфёЄртыхэшх
\eqref{first}. ╧ЁхфёЄртыхэшх \eqref{null} ёыхфєхЄ шч \eqref{first} ш \eqref{F_0est}. ╬Ўхэър \eqref{null1} яЁюшчтюфэющ ЇєэъЎшш $\eta$ т√тюфшЄё  шч
єЁртэхэш  \eqref{integr}. ─хщёЄтшЄхы№эю, фшЇЇхЁхэЎшЁє  хую, яюыєўшь
\begin{gather*}
|\eta'(x,\lambda)|\le\ch(2\pi\alpha+2|\eta(x,\lambda)|) \left(|q_1(x)|+\frac12(|q_2(x)|+|q_3(x)|)\right)\le\\\le\ch(2\pi\alpha+1)
\left(|q_1(x)|+\frac12(|q_2(x)|+|q_3(x)|)\right),
\end{gather*}
яюёъюы№ъє
$$
|\eta(x,\lambda)|\le\Upsilon(x,\lambda)+k^2\Upsilon_{p'}(\lambda)\le (1+\pi k^2)\Upsilon(\lambda)<1/2.
$$
╪ру 4 чртхЁ°хэ.
\end{proof}

\begin{Lemma}\label{lem:2.2}
╧єёЄ№ $Q\in L_1[0,\pi]$ ш т√яюыэхэ√ єёыютш  (i), (ii), (iii). ╧єёЄ№ $\lambda\in D_{Q,\alpha}$ (ёь. юяЁхфхыхэшх т ╦хььх \ref{lem:2.1}), р
$\theta(x,\lambda)$
--- Ёх°хэшх єЁртэхэш ~\eqref{th} ё эрўры№э√ь єёыютшхь $\theta(0,\lambda)=0$. ╥юуфр Ёх°хэшх $r(x,\lambda)$ єЁртэхэш ~\eqref{r} ё
эрўры№э√ь єёыютшхь $r(0,\lambda)=1$ фюяєёърхЄ яЁхфёЄртыхэшх
\begin{equation}
r(x,\lambda)=\exp\left\{\frac12\int\limits_0^x(q_3(t)-q_2(t))dt\right\}\left(1+\rho(x,\lambda)\right), \label{asr}
\end{equation}
уфх юёЄрЄюъ $\rho(x,\lambda)$ яюфўшэхэ юЎхэъх
\begin{equation}
|\rho(x,\lambda)|\le\Upsilon(x,\lambda)+2k(2+\pi k^2)\Upsilon_{p'}(\lambda),\qquad \|\rho(x,\lambda)\|_{p'}\le (1+4\pi
k+2\pi^2k^3)\Upsilon_{p'}(\lambda). \label{est}
\end{equation}
┴юыхх Єюую,
\begin{equation}\label{rfirst}
\rho(x,\lambda)=\int\limits_0^xq_1(t)\sin(2\lambda t)\,dt- \frac12\int\limits_0^x(q_2(t)+q_3(t))\cos(2\lambda t)\,dt+\sigma(x,\lambda),
\end{equation}
уфх
$$
\|\sigma(x,\lambda)\|_C\le 2k(2+\pi k^2)\Upsilon_{p'}(\lambda).
$$
╧Ёш ¤Єюь
\begin{equation}\label{rho'}
|\rho'(x,\lambda)|\le2\ch(2\pi\alpha+1)\left(|q_1(x)|+\frac12|q_2(x)+q_3(x)|\right)
\end{equation}
яюўЄш тё■фє яЁш $x\in[0,\pi]$.
\end{Lemma}
\begin{proof}
╙Ёртэхэшх \eqref{r} Ёх°рхЄё   тэю
\begin{gather*}
r(x,\lambda)=\exp\left(\frac12\int\limits_0^x(q_3(t)-q_2(t))dt+\int\limits_0^xq_1(t)\sin(2\theta(t,\lambda))dt-\right.\\\left.-
\frac12\int\limits_0^x(q_2(t)+q_3(t))\cos(2\theta(t,\lambda))dt\right)=\exp\left(\frac12\int\limits_0^x(q_3(t)-q_2(t))dt\right)\exp(H(x,\lambda)),
\end{gather*}
уфх ь√ ттхыш юсючэрўхэшх
$$
H(x,\lambda)=\int\limits_0^xq_1(t)\sin(2\theta(t,\lambda))dt- \frac12\int\limits_0^x(q_2(t)+q_3(t))\cos(2\theta(t,\lambda))dt.
$$
┬эрўрых ь√ юЎхэшь ЁрчэюёЄ№
\begin{gather}\label{rint1}
\left|\int\limits_0^xq_1(t)\sin(2\theta(t,\lambda))dt-\int\limits_0^xq_1(t)\sin(2\lambda t)dt\right|\le\notag\\
\le\left|\int\limits_0^xq_1(t)\sin(2\lambda t)(\cos(2\eta(t,\lambda))-1)dt\right|+\left|\int\limits_0^xq_1(t)\cos(2\lambda t)\sin(2\eta(t))dt\right|,
\end{gather}
тюёяюы№чютрт°шё№ яЁхфёЄртыхэшхь \eqref{null}. ┬ ёшыє эхЁртхэёЄт \eqref{complexineq}
$$
|\cos(2\eta(t,\lambda))-1|\le2|\eta(t,\lambda)|^2\ch(2|\eta(t,\lambda)|).
$$
╧юёъюы№ъє $|\eta(t,\lambda)|\le(1+k^2)\Upsilon(t,\lambda)\le1/2$, яюыєўрхь
$$
\|\cos(2\eta(t,\lambda))-1\|_{p'}\le(1+\pi k^2)(\ch1)\Upsilon_{p'}(\lambda).
$$
╥ръшь юсЁрчюь, яхЁт√щ шэЄхуЁры т яЁртющ ўрёЄш \eqref{rint1} эх яЁхтюёїюфшЄ $R\ch(2\pi\alpha)(\ch1)(1+\pi k^2)\Upsilon_{p'}(\lambda)$. ┬ЄюЁющ шэЄхуЁры
т яЁртющ ўрёЄш \eqref{rint1} юЎхэшь ёыхфє■∙шь юсЁрчюь: Єръ ъръ $|\sin(2\eta(t,\lambda))|\le2(\ch1)|\eta(t,\lambda)|$, Єю
$$
\|\sin(2\eta(t,\lambda)\|_{p'}\le2(\ch1)(1+\pi k^2)\Upsilon_{p'}(\lambda),
$$
р ёрь шэЄхуЁры эх яЁхтюёїюфшЄ $2R\ch(2\pi\alpha)(\ch1)(1+\pi k^2)\Upsilon_{p'}(\lambda)$. ╚Єръ,
$$
\left|\int\limits_0^xq_1(t)\sin(2\theta(t,\lambda))dt-\int\limits_0^xq_1(t)\sin(2\lambda t)dt\right|\le 5R\ch(2\pi\alpha)(1+\pi
k^2)\Upsilon_{p'}(\lambda).
$$
└эрыюушўэю
$$
\left|\frac12\int\limits_0^x(q_2(t)+q_3(t))\cos(2\theta(t))dt-\frac12\int\limits_0^x(q_2(t)+q_3(t))\cos(2\lambda
t)dt\right|\le5R\ch(2\pi\alpha)(1+\pi k^2)\Upsilon_{p'}(\lambda).
$$
╧юыюцшь
\begin{equation}\label{H0}
H_0(x,\lambda)=\int\limits_0^xq_1(t)\sin(2\lambda t)dt-\frac12\int\limits_0^x(q_2(t)+q_3(t))\cos(2\lambda t)dt,
\end{equation}
╥юуфр $|H_0(x,\lambda)|\le\Upsilon(x,\lambda)$. ╟рьхЄшь, ўЄю т юсырёЄш $D_{Q,\alpha}$ т√яюыэхэю
$$
|H(x,\lambda)|\le5R\ch(2\pi\alpha)(1+\pi k^2)\Upsilon_{p'}(\lambda)+\Upsilon(x,\lambda)\le(5\pi R\ch(2\pi\alpha)(1+\pi
k^2)+1)\Upsilon(\lambda)\le\frac12.
$$
╥юуфр, т ёшыє эхЁртхэёЄт \eqref{complexineq},
$$
\exp(H(x,\lambda))=\exp(H_0(x,\lambda))(1+\sigma_1(x,\lambda)),
$$
уфх
$$
\|\sigma_1(x,\lambda)\|_C\le10R\ch(2\pi\alpha)(1+\pi k^2)\Upsilon_{p'}(\lambda).
$$
╬сЁрЄшьё  ъ т√Ёрцхэш■ $\exp(H_0(x,\lambda)$. ┬эют№ яЁшьхэ   эхЁртхэёЄтр \eqref{complexineq}, яюыєўшь
$$
\left|\exp(H_0(x,\lambda))-1-H_0(x,\lambda)\right|\le|H_0(x,\lambda)|^2\le|\Upsilon(x,\lambda)|^2\le8R\ch(2\pi\alpha)\Upsilon_{p'}(\lambda),
$$
уфх ь√ тюёяюы№чютрышё№ юЎхэъющ \eqref{UpsUps}. ╚Єръ,
$$
\exp(H(x,\lambda))=(1+H_0(x,\lambda)+\sigma_0(x,\lambda))(1+\sigma_1(x,\lambda)),
$$
уфх
\begin{gather*}
\|\sigma_1(x,\lambda)\|_C\le10R\ch(2\pi\alpha)(1+\pi k^2)\Upsilon_{p'}(\lambda)\le1,\\
\|\sigma_0(x,\lambda)\|_C\le8R\ch(2\pi\alpha)\Upsilon_{p'}(\lambda)\le(1/16).
\end{gather*}
╨рёъЁ√тр  ёъюсъш т яюёыхфэхщ ЇюЁьєых, яюыєўшь яЁхфёЄртыхэшх \eqref{rfirst}. ╬Ўхэъш \eqref{est} ёыхфє■Є шч \eqref{rfirst} ш \eqref{H0}. ╬Ўхэър
\eqref{rho'} ёыхфєхЄ шч єЁртэхэш  \eqref{r} ш яЁхфёЄртыхэш  \eqref{asr}.
\end{proof}
\begin{Theorem}\label{tm:2.4}
╧єёЄ№ $Q\in L_1[0,\pi]$ ш т√яюыэхэ√ єёыютш  (i), (ii), (iii). ╧єёЄ№ $\alpha>0$ --- яЁюшчтюы№эюх ЇшъёшЁютрээюх ўшёыю, $\Pi_\alpha=\{\lambda\in\mathbb
C|\,|\text{Im}\lambda|<\alpha\}$, р
$$
D_{Q,\alpha}=\left\{\lambda\in\mathbb C: \ |Im \lambda|<\alpha,\,\Upsilon(\lambda)<\frac1{8 k^4}\right\}, \quad \text{уфх} \ \,
k=2+12R\ch(2\pi\alpha).
$$
╥юуфр фы  ы■сюую $\lambda\in D_{Q,\alpha}$ Ёх°хэшх $\mathbf{s}(x,\lambda)=(s_1(x,\lambda), s_2(x,\lambda))^t$ фюяєёърхЄ яЁхфёЄртыхэшх
\begin{gather}
s_1(x,\lambda)=-\exp\left\{\frac12\int\limits_0^x(q_3(t)-q_2(t))dt\right\}\sin(\lambda x)+\rho_1(x,\lambda),\notag\\
s_2(x,\lambda)=\exp\left\{\frac12\int\limits_0^x(q_3(t)-q_2(t))dt\right\}\cos(\lambda x)+\rho_2(x,\lambda),\label{sinas}
\end{gather}
уфх $|\rho_j(x,\lambda)|\le M(\Upsilon(x,\lambda)+\Upsilon_{p'}(\lambda))$, $M=M(R,\alpha)$. ┴юыхх Єюую,
$$
\rho_j(x,\lambda)=\rho_{j,1}(x,\lambda)\cos(\lambda x)+\rho_{j,2}(x,\lambda)\sin(\lambda x),
$$
уфх $|\rho'_{j,k}(x)|\le M(\alpha)\left(|q_1(x)|+|q_2(x)+q_3(x)|/2\right)$.
\end{Theorem}
\begin{proof}
┬ ёююЄтхЄёЄтшш ё чрьхэющ \eqref{Pr}  шьххь яЁхфёЄртыхэш 
$$\mathbf{s}_1(x,\lambda)=-r(x,\lambda)\sin\theta(x,\lambda),\qquad
\mathbf{s}_2(x,\lambda)=r(x,\lambda)\cos\theta(x,\lambda).
$$
╧юфёЄрты   Ёрчыюцхэш , яюыєўхээ√х т ╦хььрї \ref{lem:2.1} ш \ref{lem:2.2}, эряЁшьхЁ фы  $\mathbf{s}_1$ яюыєўрхь
$$
\mathbf{s}_1(x,\lambda)=-\exp\left\{\frac12\int\limits_0^x(q_3(t)-q_2(t))dt\right\}(1+\rho(x,\lambda))(\sin(\lambda
x)\cos(\eta(x,\lambda))+\cos(\lambda x)\sin(\eta(x,\lambda)).
$$
┬ юўхЁхфэющ Ёрч шёяюы№чє  юЎхэъш \eqref{complexineq}, р Єръцх \eqref{null} ш \eqref{asr}, яюёых ЁрёъЁ√Єш  ёъюсюъ, яюыєўшь рёшьяЄюЄшўхёъє■ ЇюЁьєыє фы 
$\mathbf{s}_1(x,\lambda)$. ╤ыєўрщ $\mathbf{s}_2(x,\lambda)$ рэрыюушўхэ.
\end{proof}
\begin{Theorem}\label{tm:2.3}
╧єёЄ№ $Q\in L_1[0,\pi]$ ш т√яюыэхэ√ єёыютш  (i), (ii), (iii). ╧єёЄ№ $\alpha>0$ --- яЁюшчтюы№эюх ЇшъёшЁютрээюх ўшёыю, $\Pi_\alpha=\{\lambda\in\mathbb
C|\,|\text{Im}\lambda|<\alpha\}$, р
$$
D_{Q,\alpha}=\left\{\lambda\in\mathbb C: \ |Im \lambda|<\alpha,\,\Upsilon(\lambda)<\frac1{8 k^4}\right\}, \quad \text{уфх} \ \,
k=2+12R\ch(2\pi\alpha).
$$
╥юуфр фы  ы■сюую $\lambda\in D_{Q,\alpha}$ Ёх°хэшх $\mathbf{c}(x,\lambda) =(c_1(x,\lambda), c_2(x,\lambda))^t$ фюяєёърхЄ яЁхфёЄртыхэшх
\begin{gather*}
c_1(x,\lambda)=\exp\left\{\frac12\int\limits_0^x(q_3(t)-q_2(t))dt\right\}\cos(\lambda x)+\rho_1(x,\lambda),\notag\\
c_2(x,\lambda)=\exp\left\{\frac12\int\limits_0^x(q_3(t)-q_2(t))dt\right\}\sin(\lambda x)+\rho_2(x,\lambda),\label{cosas}
\end{gather*}
уфх $|\rho_j(x,\lambda)|\le M(\Upsilon(x,\lambda)+\Upsilon_{p'}(\lambda))$, $M=M(R,\alpha)$. ┴юыхх Єюую,
$$
\rho_j(x,\lambda)=\rho_{j,1}(x,\lambda)\cos(\lambda x)+\rho_{j,2}(x,\lambda)\sin(\lambda x),
$$
уфх $|\rho'_{j,k}(x)|\le M(\alpha)\left(|q_1(x)|+|q_2(x)+q_3(x)|/2\right)$.
\end{Theorem}
\begin{proof} ╙ърцхь, ъръшх шчьхэхэш  эєцэю яЁютхёЄш, яютЄюЁ   фюърчрЄхы№ёЄтю яЁхф√фє∙хщ ЄхЁхь√.
┬ ёшёЄхьх $l_Q\mathbf{y}=\lambda\mathbf{y}$ ёфхырхь чрьхэє
\begin{equation*}
y_1(x,\lambda)=r(x,\lambda)\cos\theta(x,\lambda),\quad y_2(x,\lambda)=r(x,\lambda)\sin\theta(x,\lambda).
\end{equation*}
╥юуфр ёшёЄхьє ьюцэю чряшёрЄ№ т тшфх
\begin{equation*}
\begin{array}{ccc}
r'\sin\theta+r\theta'\cos\theta &+q_1r\cos\theta+q_2r\sin\theta &=\lambda r\cos\theta,\\
-r'\cos\theta+r\theta'\sin\theta &+q_3r\cos\theta-q_1r\sin\theta &=\lambda r\sin\theta.
\end{array}
\end{equation*}
╙ьэюцшь яхЁтюх єЁртэхэшх эр $\cos\theta$ ш яЁшсртшь тЄюЁюх єЁртэхэшх, єьэюцхээюх эр $\sin\theta$. ┬ Ёхчєы№ЄрЄх яюыєўшь єЁртэхэшх фы  ЇєэъЎшш
$\theta(x,\lambda)$
\begin{equation}\label{th1}
\theta'(x,\lambda)+q_1(x)\cos2\theta(x,\lambda)+ \frac12(q_2(x)+q_3(x))\sin2\theta(x,\lambda)=\lambda x,
\end{equation}
ё эрўры№э√ь єёыютшхь $\theta(0,\lambda)=0$. ┼ёыш ь√ єьэюцшь яхЁтюх єЁртэхэшх эр $\sin\theta$ ш т√ўЄхь тЄюЁюх єЁртэхэшх, єьэюцхээюх эр $\cos\theta$,
Єю яюыєўшь єЁртэхэшх эр ЇєэъЎш■ $r(x,\lambda)$
\begin{equation}\label{r1}
r'(x,\lambda)=r(x,\lambda)\left[\frac12(q_3(x)-q_2(x))-q_1(x)\sin 2\theta(x,\lambda)+\frac12(q_2(x)+q_3(x)) \cos 2\theta(x,\lambda)\right]
\end{equation}
ё эрўры№э√ь єёыютшхь $r(0,\lambda)=1$. ╙Ёртэхэшх \eqref{th1} юЄышўрхЄё  юЄ \eqref{th} чрьхэющ чэръют є ЇєэъЎшщ $q_j$, Єръ ўЄю Ёхчєы№ЄрЄ ╦хьь√
\ref{lem:2.1} ёюїЁрэ хЄё  ё юўхтшфэ√ьш шчьхэхэш ьш. └эрыюушўэю юсёЄюшЄ фхыю ш ё єЁртэхэшхь \eqref{r1}. ╧юфёЄрты   Ёрчыюцхэш , яюыєўхээ√х т фтєї
ыхььрї, яюыєўрхь єЄтхЁцфхэшх ЄхюЁхь√.
\end{proof}

\section{╬яхЁрЄюЁ ─шЁрър ё эєыхт√ь яюЄхэЎшрыюь}

{\bf 3.1.}  ╟фхё№ ь√ яЁштхфхь  эєцэ√х эрь ЇръЄ√  ю яЁюёЄхщ°хь юяхЁрЄюЁх ─шЁрър $L_{0,U}$, яюЁюцфхээ√ь фшЇЇхЁхэЎшры№э√ь т√Ёрцхэшхь
\begin{equation}\label{Dirac0}
l_0(\mathbf{y})=-B\mathbf{y}', \quad\text{уфх}\ \ B = \begin{pmatrix} 0 & 1 \\ -1 & 0        \end{pmatrix},
\end{equation}
ш Ёхуєы Ёэ√ь ъЁрхт√ь єёыютшхь $U(y)=0$ тшфр \eqref{matrU}. ┴юы№°р  ўрёЄ№ ¤Єшї єЄтхЁцфхэшщ їюЁю°ю шчтхёЄэ√ (яю ъЁрщэхщ ьхЁх, фы  ъюэъЁхЄэ√ї ъЁрхт√ї
єёыютшщ: яхЁшюфшўхёъшї, єёыютшщ ─шЁшїых шыш ─шЁшїых-═хщьрэр).\footnote{╘ръЄ√ юс юяхЁрЄюЁх $L_{0,U}$ ё яЁюшчтюы№э√ьш Ёхуєы Ёэ√ьш ъЁрхт√ьш єёыютш ьш,
яЁштхфхээ√х т ╥хюЁхьх \ref{tm:spL0} ш ╟рьхўрэшш \ref{no:reg}, сюы№°хщ ўрёЄ№■ ёюфхЁцрЄё  т ёЄрЄ№х \cite{DM3} ─цръютр ш ╠шЄ ушэр. ╬фэръю, фы  єфюсёЄтр
ўшЄрЄхы , ь√ яЁштхфхь шї ё фюърчрЄхы№ёЄтюь (юэю юЄышўрхЄё  юЄ фюърчрЄхы№ёЄтр т \cite{DM3}).} ┬ Єю цх тЁхь , ¤Єш ЇръЄ√ эхюсїюфшь√ фы  яюэшьрэш 
юёэютэ√ї Ёхчєы№ЄрЄют ёыхфє■∙хую ярЁруЁрЇр. ╬сючэрўшь ўхЁхч
$$
\mathcal U =\begin{pmatrix} u_{11} &  u_{12} & u_{13} & u_{14} \\
u_{21} & u_{22} & u_{23} & u_{24}\end{pmatrix}
$$
Ёрё°шЁхээє■ ьрЄЁшЎє, ёюёЄртыхээє■ шч ъю¤ЇЇшЎшхэЄют ъЁрхт√ї єёыютшщ. ╩ръ яЁхцфх ўхЁхч $J_{sj}$ юсючэрўрхь юяЁхфхышЄхыш, ёюёЄртыхээ√х шч $s$-ую ш
$j$-ую ёЄюысЎют ¤Єющ ьрЄЁшЎ√.

\begin{Theorem}\label{tm:spL0}
╤яхъЄЁ юяхЁрЄюЁр $L_{0,U}$ ёюёЄртыхэ шч ёюсёЄтхээ√ї чэрўхэшщ,
ъюЄюЁ√х ьюцэю чряшёрЄ№ фтєь  ёхЁш ьш
\begin{equation}\label{lan0}
\left[\begin{array}{l}\lambda_n=-\frac{i}{\pi}\ln z_0+2n,\quad n\in\mathbb Z\\
\lambda_n=-\frac{i}{\pi}\ln z_1+2n,\quad n\in\mathbb Z.\end{array}\right.
\end{equation}
╫шёыр $z_0,\,z_1$ хёЄ№ ъюЁэш ътрфЁрЄэюую єЁртэхэш 
\begin{equation}\label{sqeq}
[J_{14}-J_{23}-i(J_{13}+J_{24})]z^2+2[J_{12}+J_{34}]z+[J_{14}-J_{23}+i(J_{13}+J_{24})]=0,
\end{equation}
р тхЄт№ ыюурЁшЇьр фы  юяЁхфхыхээюёЄш ЇшъёшЁєхь т яюыюёх $Im\,
z\in(-\pi,\pi]$. ┬ ёыєўрх, хёыш фшёъЁшьшэрэЄ ътрфЁрЄэюую єЁртэхэш 
\eqref{sqeq} Ёртхэ эєы■, шьххь $z_0=z_1$, ш Єюуфр  тёх ёюсёЄтхээ√х
чэрўхэш  юяхЁрЄюЁр $L_{0,U}$ фтєъЁрЄэ√.
\end{Theorem}
─ы  єфюсёЄтр ь√ сєфхь т фры№эхщ°хь эєьхЁютрЄ№ ўшёыр $\lambda_n$ юфэшь шэфхъёюь $n\in\mathbb Z$, юс·хфшэ   фтх ёхЁшш \eqref{lan0} т юфэє
$\lambda_n=\varkappa_j+n$, уфх
$$
\begin{array}{ll}
\varkappa_j=\varkappa_{j(n)}=-i\pi^{-1}\ln z_0\qquad & \text{фы  ўхЄэ√ї } n;\\
\varkappa_j=\varkappa_{j(n)}=1-i\pi^{-1}\ln z_1\qquad & \text{фы  эхўхЄэ√ї } n.
\end{array}
$$
\begin{proof}
╦хуъю тшфхЄ№, ўЄю Ёх°хэшхь єЁртэхэш  $l_0(\mathbf{y})=\lambda\mathbf{y}$ ё эрўры№э√ьш єёыютш ьш $y_1(0)=1$, $y_2(0)=0$  ты хЄё  ЇєэъЎш 
$$
\mathbf{c}(x,\lambda)=\begin{pmatrix}\cos(\lambda x)\\-\sin(\lambda x)\end{pmatrix}.
$$
└эрыюушўэю, Ёх°хэшхь Єюую цх єЁртэхэш  ё эрўры№э√ьш єёыютш ьш $y_1(0)=0$, $y_2(0)=1$  ты хЄё  ЇєэъЎш 
$$
\mathbf{s}(x,\lambda)=\begin{pmatrix}\sin(\lambda x)\\ \cos(\lambda x)\end{pmatrix}.
$$
╬с∙хх Ёх°хэшх єЁртэхэш  $l_0(\mathbf{y})=\lambda\mathbf{y}$ шьххЄ тшф $\mathbf{y}=\gamma_1\mathbf{c}+\gamma_2\mathbf{s}$. ╧юфёЄрты   ¤Єю т√Ёрцхэшх т
ъЁрхт√х єёыютш  яюыєўрхь ёшёЄхьє
\begin{equation}\label{eigeneq}
\begin{cases}[u_{12}+u_{13}\sin\pi\lambda+u_{14}\cos\pi\lambda]\gamma_1+[u_{11}+u_{13}\cos\pi\lambda-u_{14}\sin\pi\lambda]\gamma_2=0,\\
[u_{22}+u_{23}\sin\pi\lambda+u_{24}\cos\pi\lambda]\gamma_1+[u_{21}+u_{23}\cos\pi\lambda-u_{24}\sin\pi\lambda]\gamma_2=0.\end{cases}
\end{equation}
╬сючэрўшь ьрЄЁшЎє ¤Єющ ёшёЄхь√ ўхЁхч $M(\lambda)$. ╫шёыю $\lambda\in\mathbb C$  ты хЄё  ёюсёЄтхээ√ь чэрўхэшхь юяхЁрЄюЁр $L_{0,U}$ Єюуфр ш Єюы№ъю
Єюуфр, ъюуфр юяЁхфхышЄхы№ $\Delta_0(\lambda): = det M(\lambda)=0$. ╙ўшЄ√тр  ъюэъЁхЄэ√щ тшф Ёх°хэшщ $\mathbf{c}(x,\lambda)$ ш $\mathbf{s}(x,\lambda)$,
яюыєўрхь ёыхфє■∙хх т√Ёрцхэшх фы  їрЁръЄхЁшёЄшўхёъюую юяЁхфхышЄхы 
\begin{equation}\label{Delta0}
\Delta_0(\lambda)=\frac12[J_{14}-J_{23}-i(J_{13}+J_{24})]e^{i\pi\lambda}+[J_{12}+J_{34}]+ \frac12[J_{14}-J_{23}+i(J_{13}+J_{24})]e^{-i\pi\lambda}.
\end{equation}
╤фхырт яюфёЄрэютъє $e^{i\pi\lambda}=z$ т єЁртэхэшш $\Delta_0(\lambda)=0$, яюыєўшь ътрфЁрЄэюх єЁртэхэшх \eqref{sqeq} фы  яхЁхьхээющ $z$. ╤фхырт
юсЁрЄэє■ яюфёЄрэютъє, яюыєўшь єЄтхЁцфхэшх ЄхюЁхь√.
\end{proof}

  ╬яхЁрЄюЁ $L_{0,U}$ эрчютхь {\it ёшы№эю Ёхуєы Ёэ√ь}, хёыш
  фшёъЁшьшэрэЄ ътрфЁрЄэюую єЁртэхэш  \eqref{sqeq} юЄышўхэ юЄ эєы ,
  Є.х. ъюЁэш $z_,\, z_1$ Ёрчышўэ√.

\begin{Note}\label{no:reg}
┼ёыш ъЁрхтюх єёыютшх эх  ты хЄё  Ёхуєы Ёэ√ь (эряЁшьхЁ, хёыш ъю¤ЇЇшЎшхэЄ яЁш $e^{i\pi\lambda}$  т \eqref{Delta0}  Ёртхэ эєы■), Єю т ёыєўрх эхЁртхэёЄтр
эєы■ ъю¤ЇЇшЎшхэЄр яЁш  $e^{-i\pi\lambda}$ ёяхъЄЁ юяхЁрЄюЁр $L_{0,U}$ сєфхЄ ёюёЄю Є№ шч юфэющ ёхЁшш яЁюёЄ√ї ёюсёЄтхээ√ї чэрўхэшщ
$\lambda_n=\varkappa+2n$, $n\in\mathbb Z$. ┼ёыш яхЁт√щ ш ЄЁхЄшщ ъю¤ЇЇшЎшхэЄ√ Ёртэ√ эєы■ юср, Єю $\Delta_0(\lambda) = J_{12}+J_{34}$. ╥юуфр
ёюсёЄтхээ√ї чэрўхэшщ эхЄ, хёыш $J_{12}+J_{34}\ne 0 $ ш ёяхъЄЁ тё  ъюьяыхъёэр  яыюёъюёЄ№, хёыш $J_{12}+J_{34}=0$. ╦хуъю тшфхЄ№, ўЄю тёх ЄЁш ёшЄєрЎшш
Ёхрышчєхь√. ═ряЁшьхЁ, $\Delta_0(\lambda) \equiv 0$, хёыш
$$
\mathcal U =\begin{pmatrix} 1 &  1 & -1 & -1 \\
\beta & -\beta & \beta & -\beta\end{pmatrix} \quad 0\ne \beta \in
\mathbb C.
$$
\end{Note}

\begin{Theorem}\label{th:eigenfun}
═юЁьшЁютрээ√х ёюсёЄтхээ√х ЇєэъЎшш $\mathbf{y}_n$, $n\in\mathbb Z$ ёшы№эю Ёхуєы Ёэюую юяхЁрЄюЁр $L_{0,U}$ шьх■Є тшф
\begin{equation}\label{eigenfunc}
\left[\begin{array}{l}
\mathbf{y}_n=\gamma_{1,0}\begin{pmatrix}\cos(\lambda_nx)\\-\sin(\lambda_nx)\end{pmatrix}+\gamma_{2,0}\begin{pmatrix}\sin(\lambda_nx)\\\cos(\lambda_nx)\end{pmatrix},\quad n\in\mathbb Z,\quad\text{хёыш}\ n\ \text{ўхЄэю},\\\phantom{.}\\
\mathbf{y}_n=\gamma_{1,1}\begin{pmatrix}\cos(\lambda_nx)\\-\sin(\lambda_nx)\end{pmatrix}+\gamma_{2,1}\begin{pmatrix}\sin(\lambda_nx)\\\cos(\lambda_nx)\end{pmatrix},\quad
n\in\mathbb Z,\quad\text{хёыш}\ n\ \text{эхўхЄэю}.\end{array}\right.
\end{equation}
╫шёыр $\gamma_{i,j}$, уфх $i=1,\,2$, $j=0,\,1$ юяЁхфхы ■Єё 
ъЁрхт√ьш єёыютш ьш ш эх чртшё Є юЄ $n$.
\end{Theorem}
\begin{proof}
╤юсёЄтхээ√х ЇєэъЎшш, ттхфхээ√х т фюърчрЄхы№ёЄтх яЁхф√фє∙хщ ЄхюЁхь√ шьх■Є тшф
$$
\mathbf{y}_n=\gamma_1\mathbf{c}(x,\lambda_n)+\gamma_2\mathbf{s}(x,\lambda_n).
$$
╧Ёш ¤Єюь ўшёыр $\gamma_1$ ш $\gamma_2$ хёЄ№ Ёх°хэш  ёшёЄхь√ \eqref{eigeneq}, т ъюЄюЁющ $\lambda=\lambda_n$. ╙ўшЄ√тр , ўЄю ьрЄЁшЎр ¤Єющ ёшёЄхь√
$2$--яхЁшюфшўэр яю ярЁрьхЄЁє $\lambda$ яЁшїюфшь ъ Єюьє, ўЄю ўшёыр $\gamma_1$ ш $\gamma_2$ чртшё Є ыш°№ юЄ ўхЄэюёЄш шэфхъёр $n$. ╬ёЄрхЄё  яюфёЄртшЄ№
т√Ёрцхэш  фы  ЇєэъЎшщ $\mathbf{c}(x,\lambda_n)$ ш $\mathbf{s}(x,\lambda_n)$ ш ь√ яЁшфхь ъ ЇюЁьєырь \eqref{eigenfunc} фы  эхэюЁьшЁютрээ√ї ёюсёЄтхээ√ї
ЇєэъЎшщ. ╬ўхтшфэю, эюЁьр $\|\mathbf{y}_n\|$ Єръцх чртшёшЄ ыш°№ юЄ ўхЄэюёЄш шэфхъёр $n$. ╥ръшь юсЁрчюь, ЇюЁьєы√ \eqref{eigenfunc} ёюїЁрэ ■Єё  фы 
$\mathbf{y}_n/\|\mathbf{y}_n\|$.
\end{proof}

\begin{Theorem}\label{tm:Green}
╨хчюы№тхэЄр $R_0(\lambda)=(L_{0,U}-\lambda I)^{-1}$ Ёхуєы Ёэюую юяхЁрЄюЁр $L_{0,U}$  ты хЄё  шэЄхуЁры№э√ь юяхЁрЄюЁюь
\begin{equation}\label{Green}
R_0(\lambda)\mathbf{f}=\int_0^\pi G_0(x,t,\lambda)\mathbf{f}(t)dt.
\end{equation}
╘єэъЎш  $G_0(x,t,\lambda)$ эхяЁхЁ√тэр эр ътрфЁрЄх $(x,t)\in[0,\pi]^2$ чр шёъы■ўхэшхь фшруюэрыш $x=t$. ┬эх $\delta$- ъЁєцъют ё ЎхэЄЁрьш т эєы ї
$\lambda_n$  юяЁхфхышЄхы  $\Delta_0(\lambda)$ (т ўрёЄэюёЄш, тэх эхъюЄюЁющ яюыюё√ $|Im \lambda|>\alpha$) ЇєэъЎш  $G_0(x,t,\lambda)$ єфютыхЄтюЁ хЄ
юЎхэъх $|G_0(x,t,\lambda)|\le M$ ё эхъюЄюЁющ ъюэёЄрэЄющ $M$ чртшё ∙хщ юЄ ъЁрхт√ї єёыютшщ ш ўшёыр $\delta$  (шыш $\alpha$), эю эх чртшё ∙хщ юЄ $x$,
$t$ ш $\lambda$.
\end{Theorem}
\begin{proof}
╧Ёшьхэшь ьхЄюф  трЁшрЎшш яюёЄю ээ√ї ъ єЁртэхэш■ $l_0(\mathbf{y})=\lambda\mathbf{y}+\mathbf{f}$. ╥юуфр Ёх°хэшх яЁшьхЄ тшф
\begin{equation}\label{Varconst}
\begin{cases}
y_1(x)=\gamma_1\sin\lambda x+\gamma_2\cos\lambda x-\int_x^\pi\cos\lambda(x-t)f_1(t)dt+\int_x^\pi\sin\lambda(x-t)f_2(t)dt,\\
y_2(x)=\gamma_1\cos\lambda x-\gamma_2\sin\lambda x+\int_x^\pi\sin\lambda(x-t)f_1(t)dt+\int_x^\pi\cos\lambda(x-t)f_2(t)dt.\end{cases}
\end{equation}
╧юфёЄртшт ¤Єю Ёх°хэшх т ъЁрхт√х єёыютш , яюыєўшь ёшёЄхьє шч фтєї
єЁртэхэшщ фы  юяЁхфхыхэш  ўшёхы $\gamma_1$ ш $\gamma_2$:
$$
 M_0\begin{pmatrix}\gamma_1\\ \gamma_2\end{pmatrix}=
\begin{pmatrix}u_{11} & u_{12}\\ u_{21} & u_{22}\end{pmatrix}\hat{f}(\lambda),
$$
уфх ьрЄЁшЎр $M_0(\lambda)$ юяЁхфхыхэр т ╥хюЁхьх \ref{tm:spL0}, р
$$
\hat{f}(\lambda)=\begin{pmatrix}\int\limits_0^\pi\cos\lambda t
f_1(t)dt+\int\limits_0^\pi\sin\lambda t f_2(t)dt\\
\int\limits_0^\pi\sin\lambda t f_1(t)dt-\int\limits_0^\pi\cos\lambda t f_2(t)dt\end{pmatrix}.
$$
╥ръшь юсЁрчюь,
\begin{multline*}
\begin{pmatrix}\gamma_1\\ \gamma_2\end{pmatrix}=\frac1{\Delta_0(\lambda)}
\begin{pmatrix}u_{21}+u_{23}\cos\pi\lambda-u_{24}\sin\pi\lambda & -u_{11}-u_{13}\cos\pi\lambda+u_{14}\sin\pi\lambda\\
-u_{22}-u_{23}\sin\pi\lambda-u_{24}\cos\pi\lambda & u_{12}+u_{13}\sin\pi\lambda+u_{14}\cos\pi\lambda\end{pmatrix}\\
\times\begin{pmatrix}u_{11} & u_{12}\\ u_{21} & u_{22}\end{pmatrix}\hat{f}(\lambda)
\end{multline*}
ш ЇюЁьєыр \eqref{Green}  ёыхфєхЄ шч \eqref{Varconst}.

─юърцхь юЎхэъє ЇєэъЎшш $G_0(x,t,\lambda)$. ╧Ёютхфхь фюърчрЄхы№ёЄтю эр ЄЁхєуюы№эшъх $t<x$ --- эр тЄюЁюь ЄЁхєуюы№эшъх юЎхэъш яЁютюф Єё  рэрыюушўэю.
╤юуырёэю \eqref{Varconst}, ьрЄЁшЎ--ЇєэъЎш  $G_0(x,t,\lambda)$  яЁш $t < x$ шьххЄ тшф
\begin{multline*}
G_0(x,t,\lambda)=\frac1{\Delta_0(\lambda)}\begin{pmatrix}\sin\lambda x & \cos\lambda
x\\ \cos\lambda x & -\sin\lambda x\end{pmatrix}\times\\
\times\begin{pmatrix}u_{21}+u_{23}\cos\pi\lambda-u_{24}\sin\pi\lambda &
-u_{11}-u_{13}\cos\pi\lambda+u_{14}\sin\pi\lambda\\-u_{22}-u_{23}\sin\pi\lambda-u_{24}\cos\pi\lambda &
u_{12}+u_{13}\sin\pi\lambda+u_{14}\cos\pi\lambda\end{pmatrix}\\
\times\begin{pmatrix}u_{11} & u_{12}\\ u_{21} & u_{22}\end{pmatrix}\begin{pmatrix}\cos\lambda t & \sin\lambda t\\ \sin\lambda t & -\cos\lambda
t\end{pmatrix}.
\end{multline*}
╧хЁхьэюцр  ьрЄЁшЎ√ (ь√ юяєёърхь чфхё№ т√ўшёыхэш ), яюыєўрхь
\begin{multline*}
G_0(x,t,\lambda)=\frac1{\Delta_0(\lambda)}\left[\begin{pmatrix}0&J_{12}\\J_{12}&0\end{pmatrix}\sin\lambda(x-t)+
\begin{pmatrix}-J_{12}&0\\0&J_{12}\end{pmatrix}\cos\lambda(x-t)+
\right.\\+\begin{pmatrix}J_{23}&J_{13}\\J_{24}&J_{14}\end{pmatrix}\sin\lambda
t\sin\lambda(x-\pi)+\begin{pmatrix}-J_{24}&-J_{14}\\J_{23}&J_{13}\end{pmatrix}\sin\lambda t\cos\lambda(\pi-x)+\\ \left.
+\begin{pmatrix}J_{13}&-J_{23}\\J_{14}&-J_{24}\end{pmatrix}\cos\lambda
t\sin\lambda(x-\pi)+\begin{pmatrix}-J_{14}&J_{24}\\J_{13}&-J_{23}\end{pmatrix}\cos\lambda t\cos\lambda(\pi-x)\right].
\end{multline*}
╚ч яЁхфёЄртыхэш  \eqref{Delta0}, юўхтшфэю, ёыхфєхЄ юЎхэър $\Delta_0(\lambda) \geqslant Ce^{|\pi\, Im\, \lambda|}$   тэх $\delta$- ъЁєцъют ё ЎхэЄЁрьш
т эєы ї $\lambda_n$   ё ъюэёЄрэЄющ $C$, чртшё ∙хщ Єюы№ъю юЄ $\delta$ ш ўшёхы $J_{sj}$. ═ю тёх ёырурхь√х т яюёыхфэхь яЁхфёЄртыхэшш, чръы■ўхээ√х т
ътрфЁрЄэ√х ёъюсъш, юўхтшфэю, юЎхэштр■Єё  ётхЁїє Єръющ цх тхышўшэющ (эю ё фЁєующ ъюэёЄрэЄющ) яЁш $0\leqslant t\leqslant x\leqslant \pi$.  ▌Єю
фюърч√трхЄ ЄхюЁхьє
\end{proof}

\begin{Corollary}\label{co:complete}
╤шёЄхьр ёюсёЄтхээ√ї ЇєэъЎшщ Ёхуєы Ёэюую юяхЁрЄюЁр $L_{0,U}$
 ты хЄё  яюыэющ т яЁюёЄЁрэёЄтх $\mathbb H$.
\end{Corollary}
\begin{proof}
╬сючэрўшь $R_0(\lambda) = (L_{0,U} -\lambda)^{-1}$.  ┼ёыш   эрщфхЄё  тхъЄюЁ $\bold f\in \mathbb H$ юЁЄюуюэры№э√щ тёхь ёюсёЄтхээ√ь ЇєэъЎш ь юяхЁрЄюЁр
$L_{0,U}$, Єю тхъЄюЁ-ЇєэъЎш  $\bold F(\lambda) = R^*_0(\lambda)\bold f$  ты хЄё  Ўхыющ (єёыютшх юЁЄюуюэры№эюёЄш $\bold f$ ёюсёЄтхээ√ь ЇєэъЎш ь
юсхёяхўштрхЄ ЁртхэёЄтю эєы■ тёхї т√ўхЄют т яюы■ёрї Ёхчюы№тхэЄ√ ёюяЁ цхээюую юяхЁрЄюЁр). ┬эх $\delta$-ъЁєцъют ё ЎхэЄЁрьш т ёюсёЄтхээ√ї чэрўхэш ї
$\overline\lambda_n$ ёюяЁ цхээюую юяхЁрЄюЁр  ЇєэъЎш  $\bold F(\lambda) $ юЎхэштрхЄё  ъюэёЄрэЄющ. ╧Ёш фюёЄрЄюўэю ьрыюь $\delta>0$  ¤Єш ъЁєцъш эх
яхЁхёхър■Єё . ╚ч яЁшэЎшяр ьръёшьєьр ш ЄхюЁхь√ ╦шєтшыы  Єюуфр ёыхфєхЄ, ўЄю $\bold F(\lambda)= \bold g$ --- яюёЄю ээ√щ тхъЄюЁ. ═ю ёяЁртхфыштю ЁртхэёЄтю
$(L^*_{0,U}-\lambda)\bold f = \bold g$.  ╧Ёртр  ўрёЄ№ юЄ $\lambda$ эх чртшёшЄ, яю¤Єюьє $\bold f =0$.
\end{proof}

{\bf 3.2.}  ╚ч юс∙хщ ЄхюЁшш (ёь. \cite [├ы. 3, \S 6]{Na})  ёыхфєхЄ, ўЄю
 ёюяЁ цхээ√щ  ъ $L_{0,U}$  юяхЁрЄюЁ $L^*_{0,U} = L_{0,U^*}$
 яюЁюцфхэ Єхь цх фшЇЇхЁхэЎшры№э√ь т√Ёрцхэшхь \eqref{Dirac0} ш
ъЁрхт√ь єёыютшхь
$$
U^*(y) = \hat A \bold y(0) + \hat D \bold y(\pi) = 0,
$$
уфх $\hat A, \hat B$ ---  эхъюЄюЁ√х $2\times 2$ ьрЄЁшЎ√. ╧ютЄюЁ   Ёрёёєцфхэш  шч ЁрсюЄ√ ╩Ёхщэр \cite{Kr},  ьюцэю яюърчрЄ№,  ўЄю ¤Єш ьрЄЁшЎ√
юяЁхфхы ■Єё  ёююЄэю°хэшхь
$$
AJ\hat A = BJ\hat B,\qquad  \text{уфх}\ \ J=\begin{pmatrix} 0 & 1\\1 & 0\end{pmatrix}.
$$
╧юфЁюсэюх фюърчрЄхы№ёЄтю ¤Єюую ЇръЄр ь√ эх яЁштюфшь, Єръ ъръ фрыхх ¤Єю єЄтхЁцфхэшх эх шёяюы№чєхЄё . ┬рцэю чрьхЄшЄ№ ёыхфє■∙хх: {\sl ┼ёыш ъЁрхтюх
єёыютшх  $U(\bold y) =0$  ты хЄё  Ёхуєы Ёэ√ь, Єю Єхь цх ётющёЄтюь юсырфрхЄ  ёюяЁ цхээюх ъЁрхтюх єёыютшх $U(\bold y) =0$.} ─хщёЄтшЄхы№эю, ёюсёЄтхээ√х
чэрўхэш  $\lambda_n$  ш $\overline{\lambda_n}$ юяхЁрЄюЁют $L_{0,U}$ ш $L_{0,U^*}$ тчршьэю ёюяЁ цхэ√, яЁшўхь шї ъЁрЄэюёЄш ёютярфр■Є.  ┼ёыш ъЁрхтюх
єёыютшх $U^*(\bold y ) =0$ эхЁхуєы Ёэюх, Єю ЁхрышчєхЄё  юфэр шч тючьюцэюёЄхщ, юяшёрээ√ї т ╟рьхўрэшш \ref{no:reg}. ═ю ¤Єю эх ёюуырёєхЄё  ё ЁртхэёЄтрьш
$\lambda_n= \overline{\lambda_n}$  ё ёюїЁрэхэшхь ъЁрЄэюёЄхщ. ╥ръшь юсЁрчюь, ёюяЁ цхээ√щ юяхЁрЄюЁ  $L_{0,U^*}$ Ёхуєы Ёхэ юфэютЁхьхээю ё $ L_{0,U}.$

─ы  фюърчрЄхы№ёЄтр  ╥хюЁхь√ \ref{tm:Riesz0}  эрь яюэрфюсшЄё  ёыхфє■∙хх яЁюёЄюх єЄтхЁцфхэшх.
\begin{Lemma}\label{lem:scprod}
╧єёЄ№
\begin{gather*}
\mathbf{y}=\gamma_1\begin{pmatrix}\cos(n+\omega)x\\-\sin(n+\omega)x\end{pmatrix}+\gamma_2\begin{pmatrix}\sin(n+\omega)x\\
\cos(n+\omega)x\end{pmatrix},\\
\mathbf{z}=\gamma'_1\begin{pmatrix}\cos(n+\omega')x\\-\sin(n+\omega')x\end{pmatrix}+
\gamma'_2\begin{pmatrix}\sin(n+\omega')x\\
\cos(n+\omega')x\end{pmatrix}.
\end{gather*}
╥юуфр
$$
(\mathbf{y},\mathbf{z})_\mathbb
H=(\gamma_1\overline{\gamma_1}'+\gamma_2\overline{\gamma_2}')\frac{\sin\pi(\omega-\overline{\omega}')}{\omega-\overline{\omega}'}+
(\gamma_1\overline{\gamma_2}'-\gamma_2\overline{\gamma_1}')\frac{1-\cos\pi(\omega-\overline{\omega}')}{\omega-\overline{\omega}'},
$$
уфх яЁш ЁртхэёЄтх $\omega=\overline{\omega}'$ яхЁтр  фЁюс№
яюырурхЄё  Ёртэющ $\pi$, р тЄюЁр  --- эєы■.
\end{Lemma}
\begin{proof}
─юърчрЄхы№ёЄтю яюыєўрхЄё  эхяюёЁхфёЄтхээ√ь т√ўшёыхэшхь.
\end{proof}

═ряюьэшь, ўЄю ёшёЄхьр тхъЄюЁют т ушы№схЁЄютюь яЁюёЄЁрэёЄтх $H$ эрч√трхЄё  срчшёюь ╨шёёр, хёыш ёє∙хёЄтєхЄ юуЁрэшўхээ√щ ш юсЁрЄшь√щ юяхЁрЄюЁ т $H$,
ъюЄюЁ√щ яхЁхтюфшЄ ¤Єє ёшёЄхьє т юЁЄюэюЁьшЁютрээ√щ срчшё.  ╤шёЄхьр $\varphi_k$ т $H$  эрч√трхЄё  {\it схёёхыхтющ} ёшёЄхьющ,  хёыш фы  ы■сюую тхъЄюЁр
$f\in H$  Ё ф $\sum |\langle f, g\rangle|^2 $  ёїюфшЄё .  ╤шёЄхь√ $\{\varphi_k\} $  ш $\{\psi_k\}$ эрч√тр■Єё  \textit{сшюЁЄюуюэры№э√ьш}, хёыш
$\langle \varphi_k,\psi_j\rangle= \delta_{kj}$,  уфх  $\delta_{kj}$ --- ёшьтюы ╩ЁюэхъхЁр. ╤шёЄхьє  $\{\varphi_k\} $ эрч√тр■Є \textit{яюўЄш
эюЁьшЁютрээющ}, хёыш ўшёыр $\|\varphi_k\| $  ЁртэюьхЁэю юуЁрэшўхэ√ ш ЁртэюьхЁэю юЄфхыхэ√ юЄ эєы .

\begin{Theorem}\label{tm:Riesz0}
═юЁьшЁютрээр  ёшёЄхьр ёюсёЄтхээ√ї ЇєэъЎшщ ёшы№эю  Ёхуєы Ёэюую юяхЁрЄюЁр $L_{0,U}$  юсЁрчєхЄ срчшё ╨шёёр т яЁюёЄЁрэёЄтх $\mathbb H$.
\end{Theorem}
\begin{proof}\!\footnote{╤Ё. ё \cite[╦хььр 3.3]{DM3}).}
┬ёяюьэшь, ўЄю ёюсёЄтхээ√х ЇєэъЎшш $\mathbf{y}^0_n$  юяхЁрЄюЁр $L_{0, U}$ шьх■Є яЁхфёЄртыхэшх \eqref{eigenfunc}. ╥ръ ъръ юяхЁрЄюЁ $L_{0, U^*}$ Єюцх
ёшы№эю Ёхуєы Ёхэ, Єю хую ёшёЄхьр ёюсёЄтхээ√ї ЇєэъЎшщ, ёюуырёэю ╥хюЁхьх~\ref{th:eigenfun}, шьххЄ  рэрыюушўэюх яЁхфёЄртыхэшх
\begin{equation*}
\left[\begin{array}{l} \mathbf{z}_n=\gamma'_{1,0}\begin{pmatrix}\cos(\overline{\lambda_n}x)\\-\sin(\overline{\lambda_n}x)\end{pmatrix}+
\gamma'_{2,0}\begin{pmatrix}\sin(\overline{\lambda_n}x)\\\cos(\overline{\lambda_n}x)\end{pmatrix},\quad n\in\mathbb Z,\quad\text{хёыш}\ n\ \text{ўхЄэю},\\\phantom{.}\\
\mathbf{z}_n=\gamma'_{1,1}\begin{pmatrix}\cos(\overline{\lambda_n}x)\\-\sin(\overline{\lambda_n}x)\end{pmatrix}+
\gamma'_{2,1}\begin{pmatrix}\sin(\overline{\lambda_n}x)\\\cos(\overline{\lambda_n}x)\end{pmatrix},\quad n\in\mathbb Z,\quad\text{хёыш}\ n\
\text{эхўхЄэю},\end{array}\right.
\end{equation*}
╥ръ ъръ тёх ёюсёЄтхээ√х чэрўхэш  юяхЁрЄюЁр $L_{0, U}$ яЁюёЄ√х, Єю ёяЁртхфышт√ ёююЄэю°хэш  $\langle\bold y_n, \bold z_j\rangle =\alpha_n\delta_{nj}$,
уфх ўшёыр $\alpha_n$  чртшё Є юЄ эюЁьшЁютъш ёшёЄхь ЇєэъЎшщ   $\mathbf{y}^0_n$ ш  $\mathbf{z}^0_n$. ═ю, ёюуырёэю ╦хььх~\ref{lem:scprod},   ўшёыр
$\alpha_n$ чртшё Є Єюы№ъю юЄ ўхЄэюёЄш $n$,
 Є.х. $\alpha_{2n} =\alpha_0, \alpha_{2n+1} = \alpha_1$. ─юьэюцшт тхъЄюЁ√ ёшёЄхь√ $\{\bold z_n\}$ ё ўхЄэ√ьш
эюьхЁрьш эр ъюэёЄрэЄє $\alpha_0^{-1}$, р тхъЄюЁ√ ё эхўхЄэ√ьш эюьхЁрьш
--- эр ъюэёЄрэЄє $\alpha_1^{-1}$  яюыєўшь, ўЄю $\{\bold y_n\}$  ш
$\{\bold z_n\}$  тчршьэю сшюЁЄюуюэры№э√.  ╬сх ¤Єш ёшёЄхь√, юўхтшфэю, схёёхыхт√ ш юсх яюыэ√ т $\mathbb H$ ёюуырёэю ╤ыхфёЄтш■~\ref{co:complete}. ╥юуфр
т ёшыє ЄхюЁхь√ ┴рЁш--┴юрёр (ёь. \cite[├ы. 6]{GK})  юсх ёшёЄхь√ юсЁрчє■Є срчшё ╨шёёр т $\mathbb H$.
\end{proof}

╧ю ёэшь ЄхяхЁ№, ўЄю   яЁюшёїюфшЄ т ёыєўрх Ёхуєы Ёэ√ї, эю эх ёшы№эю
Ёхуєы Ёэ√ї ъЁрхт√ї єёыютш ї. ╧юыхчэю ЁрёёьюЄЁхЄ№ юяхЁрЄюЁ ─шЁрър ё
эєыхт√ь яюЄхэЎшрыюь ш
 ш ъЁрхт√ьш єёыютш ьш $y_1(0)=\alpha y_1(\pi)$, $y_2(0)=y_2(\pi)$. ╧Ёш
ы■сюь $\alpha\in\mathbb C$ ёяхъЄЁ ¤Єюую юяхЁрЄюЁр  ёюёЄюшЄ шч ёюсёЄтхээ√ї чэрўхэшщ $\lambda_n=2n$. ╩рцфюх Єръюх ёюсёЄтхээюх чэрўхэшх шьххЄ
рыухсЁршўхёъє■ ъЁрЄэюёЄ№ 2. ╧Ёш $\alpha=1$ ухюьхЄЁшўхёър  ъЁрЄэюёЄ№ ърцфюую ёюсёЄтхээюую чэрўхэш  Єръцх Ёртэр фтєь: хёЄ№ фтх юЁЄюуюэры№э√х фЁєу фЁєує
ёюсёЄтхээ√х ЇєэъЎшш $\mathbf{c}^0(x,\lambda_n)$ ш $\mathbf{s}^0(x,\lambda_n)$. ╧Ёш $\alpha\ne1$ ухюьхЄЁшўхёър  ъЁрЄэюёЄ№ Ёртэр хфшэшЎх: ёюсёЄтхээющ
ЇєэъЎшхщ  ты хЄё  Єюы№ъю ЇєэъЎш  $\mathbf{s}^0(x,\lambda_n)$ ш ъ эхщ шьххЄё  яЁшёюхфшэхээр . ▌Єю Єшяшўэ√щ яЁшьхЁ, ўЄю ёыхфєхЄ шч ёыхфє■∙хую
яЁхфыюцхэш .
\begin{Lemma}
╧єёЄ№ юяхЁрЄюЁ ─шЁрър $L_{0,U}$ Ёхуєы Ёхэ, эю эх ёшы№эю Ёхуєы Ёхэ. ╥юуфр ышсю тёх хую  ёюсёЄтхээ√х чэрўхэш  $\lambda_n= 2n +\kappa$ шьх■Є
ухюьхЄЁшўхёъє■ ъЁрЄэюёЄ№ 2 ш т ¤Єюь ёыєўрх шь юЄтхўр■Є ёюсёЄтхээ√х ЇєэъЎшш $\mathbf{c}^0(x,\lambda_n)$  ш $\mathbf{s}^0(x,\lambda_n)$,
 ышсю тёхь хую ёюсёЄтхээ√ь чэрўхэш ь юЄтхўр■Є цюЁфрэют√ Ўхяюўъш, Є.х. Єюы№ъю юфэр
 ёюсёЄтхээр  ЇєэъЎш  тшфр
$$
\mathbf{y}^0_n(x)=\gamma_1\mathbf{c}^0(x,\lambda_n)+\gamma_2\mathbf{s}^0(x,\lambda_n)
$$
ш яЁшёюхфшэхээр  ъ эхщ ЇєэъЎш  тшфр
\begin{equation}\label{assoc}
\mathbf{y}^1_n(x)=\gamma_2x\mathbf{c}^0(x,\lambda_n)-\gamma_1x\mathbf{s}^0(s,\lambda_n) +
\beta\gamma_1\mathbf{c}^0(x,\lambda_n)+\beta\gamma_2\mathbf{s}(s^0,\lambda_n).
\end{equation}
╟фхё№ ўшёыр $\gamma_1, \gamma_2$  юЄ $n$ эх чртшё Є.  ╫шёыю $\beta$ ьюцэю яюфюсЁрЄ№ Єръ, ўЄю ЇєэъЎшш $\mathbf{y}_n^0$ ш $\mathbf{y}_n^1$ сєфєЄ
тчршьэю юЁЄюуюэры№э√ьш. ╧Ёш ¤Єюь ўшёыю $\beta$ ш эюЁь√ $\|\mathbf{y}^0_n\|$, $\|\mathbf{y}^1_n\|$\ Єръцх эх чртшё Є юЄ $n$.
\end{Lemma}
\begin{proof}
╤юуырёэю юяЁхфхыхэш■, т єёыютш ї ыхьь√ шьххь, ўЄю ёюсёЄтхээ√х чэрўхэш  юяхЁрЄюЁр $L_{0,U}$ шьх■Є тшф $\lambda_n=2n+\omega$, $n\in\mathbb Z$.
 ═рышўшх шыш юЄёєЄёЄтшх тЄюЁюую ёюсёЄтхээюую тхъЄюЁр, юЄтхўр■∙хую
ёюсёЄтхээюьє чэрўхэш■ $\lambda_n$ юяЁхфхы хЄё  ъЁрхт√ьш єёыютш ьш, р Єюўэхх --- ьрЄЁшЎхщ $M_0(\lambda)$, ттхфхээющ т \eqref{eigeneq}. ▌Єр ьрЄЁшЎр
$2$--яхЁшюфшўэр,  яю¤Єюьє ёюсёЄтхээ√х ЇєэъЎшш   ты ■Єё  ышэхщэ√ьш ъюьсшэрЎш ьш $\mathbf{c}^0(x,\lambda_n)$  ш  $\mathbf{s}^0(x,\lambda_n)$  ё
яюёЄю ээ√ьш $\gamma_1, \gamma_2$,  эх чртшё ∙шьш юЄ $n$. ╧ю юяЁхфхыхэш■, яЁшёюхфшэхээр  ЇєэъЎш  $\mathbf{y}^1(x)$ фюыцэр єфютыхЄтюЁ Є№ ёшёЄхьх
$$
\begin{cases}-y_2'-\lambda_ny_1=\gamma_1\cos(\lambda_nx)+\gamma_2\sin(\lambda_nx),\\
y_1'-\lambda_ny_2=-\gamma_1\sin(\lambda_nx)+\gamma_2\cos(\lambda_nx),
\end{cases}
$$
юс∙шь Ёх°хэшхь  ъюЄюЁющ   ты хЄё  ЇєэъЎш 
$$
\gamma_2x\mathbf{c}^0(x,\lambda_n)-\gamma_1x\mathbf{s}^0(s,\lambda_n) +
\beta_1\mathbf{c}^0(x,\lambda_n)+\beta_2\mathbf{s}(s^0,\lambda_n)
$$
ё яЁюшчтюы№э√ьш яюёЄю ээ√ьш $\beta_1, \beta_2$.  ╧юфёЄрты   ¤Єє ЇєэъЎш■ т ъЁрхтюх
єёыютшх, эрщфхь ёююЄэю°хэшх ьхцфє $\beta_1$ ш $\beta_2$.
 ╧ЁютхЁър эхчртшёшьюёЄш эюЁь ЇєэъЎшщ $\mathbf{y}^0_n$,
$\mathbf{y}_n^1$  юЄ $n$  юёє∙хёЄты хЄё  ¤ыхьхэЄрЁэ√ь т√ўшёыхэшхь рэрыюушўэю ╦хььх \ref{lem:scprod}.
\end{proof}

┼ёыш ёюсёЄтхээ√ь чэрўхэш ь $\lambda_n$ юяхЁрЄюЁр $L_{0,U}$ юЄтхўр■Є цюЁфрэют√ Ўхяюўъш $\mathbf{y}^0_n, \mathbf{y}^1_n$, Єю ёюсёЄтхээ√ь чэрўхэш ь
$\overline \lambda_n$ ёюяЁ цхээюую юяхЁрЄюЁр $L^*_{0,U}$ Єръцх юЄтхўр■Є цюЁфрэют√ Ўхяюўъш $\mathbf{z}^0_n, \mathbf{z}^1_n$  Єюую цх тшфр. ╚ч юс∙шї
ёююЄэю°хэшщ сшюЁЄюуюэры№эюёЄш фы  тчршьэю ёюяЁ цхээ√ї юяхЁрЄюЁют (ёь., эряЁшьхЁ, \cite{Ke}) ёыхфєхЄ
$$
\langle\mathbf{z}^0_n, \mathbf{y}^0_n\rangle=0,\qquad \langle\mathbf{z}^0_n, \mathbf{y}^1_n\rangle=\alpha_n=\alpha\ne 0,\qquad \langle\mathbf{z}^1_n,
\mathbf{y}^0_n\rangle=\alpha_n'=\alpha'\ne0.
$$
╧Ёютюф , хёыш эрфю яхЁхэюЁьшЁютъє, ьюцхь ёўшЄрЄ№, ўЄю $\alpha=\alpha'=1$. ╟рьхэ  , хёыш эєцэю $\mathbf{z}^1_n$ эр $\mathbf{z}_n^1 +\beta
\mathbf{z}^0_n$ (Єрър  ЇєэъЎш  юёЄрхЄё  яЁшёюхфшэхээющ ъ $\mathbf{z}^0_n$),  ьюцэю фюсшЄ№ё  ЁртхэёЄтр $\langle\mathbf{y}^1_n, \mathbf{z}^1_n\rangle
=0$. ╥юуфр ьюцэю ёўшЄрЄ№, ўЄю {\sl ёшёЄхь√ $\{\mathbf{y}^0_n, \mathbf{y}^1_n\} $ ш $\{\mathbf{z}^0_n, \mathbf{z}^1_n\}$   ты ■Єё  тчршьэю
сшюЁЄюуюэры№э√ьш.} ═ю шч тшфр ¤Єшї ёшёЄхь ёыхфєхЄ, ўЄю юэш схёёхыхт√. ╬сх ёшёЄхь√ яюыэ√ т $\mathbb H$, р яюЄюьє ърцфр  шч эшї юсЁрчєхЄ срчшё ╨шёёр.
╥хь ёрь√ь фюърчрэ ёыхфє■∙шщ ЇръЄ (ёЁ. ё \cite[╦хьь√ 3.5 ш 3.6]{DM3}).

\begin{Theorem}\label{tm:3.9} ╤юсёЄтхээ√х ш яЁшёюхфшэхээ√х ЇєэъЎшш Ёхуєы Ёэюую, эю эх ёшы№эю Ёхуєы Ёэюую юяхЁрЄюЁр $L_{0,U}$
ьюцэю т√сЁрЄ№ Єръ, ўЄюс√ юэш юсЁрчют√трыш срчшё ╨шёёр т $\mathbb H$.
\end{Theorem}

{\bf 3.3}. ═ръюэхЎ фюърцхь трцэє■ т фры№эхщ°хь юЎхэъє  Ёхчюы№тхэЄ√ $R_0(\lambda) =(L_{0,U} -\lambda)^{-1}$  ъръ юяхЁрЄюЁр, фхщёЄтє■∙хую шч т
яЁюёЄЁрэёЄтрї $L_p[0,\pi]$. ─рыхх ўхЁхч $l_p$ юсючэрўр■Єё  юс√ўэ√х яЁюёЄЁрэёЄтр фтєёЄюЁюээшї ўшёыют√ї яюёыхфютрЄхы№эюёЄхщ $\{c_n\}$  ё эюЁьющ
$$
\|\{c_n\}\| = \left(\sum |c_n|^p\right)^{1/p}.
$$
┬ фюърчрЄхы№ёЄтх сєфхЄ шёяюы№чютрэ ёыхфє■∙шщ шчтхёЄэ√щ Ёхчєы№ЄрЄ (ёь., эряЁшьхЁ, \cite{BL}).

\noindent {\bf ╥хюЁхьр} (╒рєёфюЁЇр-▐эу). {\sl ╧єёЄ№ $\varphi_n(x) =e^{2\pi i nx}$, $c_n(f) = (f(x), \varphi_n(x))$.  ╥юуфр  яЁш тёхї $1\leqslant
p\leqslant 2$ юяхЁрЄюЁ
$$
Tf = \{c_n(f)\}_{-\infty}^\infty
$$
 ты хЄё  юуЁрэшўхээ√ь шч яЁюёЄЁрэёЄтр $L_p =L_p[0,\pi]$ т
яЁюёЄЁрэёЄтю $l_{p'}$,  уфх $1/p + 1/p' =1.$}

\begin{Theorem}\label{th:resolvet0}
╧єёЄ№ $R_0(\lambda)$ --- Ёхчюы№тхэЄр Ёхуєы Ёэюую юяхЁрЄюЁр $L_{0,U}$. ╧юыюцшь $\lambda = \mu +i\tau$  ш т√схЁхь  ўшёыю $\alpha>0$ Єръ, ўЄюс√ тёх
яюы■ёр Ёхчюы№тхэЄ√ ыхцрыш тэєЄЁш яюыюё√ $|\rm{Im}\,\lambda| <\alpha$.  ╥юуфр яЁш тёхї $\lambda$  тэх ¤Єющ яюыюё√ ш яЁш тёхї $1\le p\le 2\le q\le
\infty$  ёяЁртхфыштр юЎхэър
\begin{equation}\label{resolvent0}
\|R_0(\lambda)\|_{L_p\to L_q} \leqslant C_{p,q}|\tau|^{-1+1/p-1/q},
\end{equation}
ё ъюэёЄрэЄющ $C_{p,q}$,  чртшё ∙хщ Єюы№ъю юЄ $p$ ш $q$, эю эх юЄ $\lambda.$
\end{Theorem}

 \begin{proof}
╧єёЄ№  $L_{0,U}$ ёшы№эю Ёхуєы Ёхэ. ╨рёёьюЄЁшь ёэрўрыр ёыєўрщ $1<p<q<\infty$. ╬сючэрўшь ўхЁхч $p', q'$ ўшёыр, ёюяЁ цхээ√х яю ├хы№фхЁє ъ $p$ ш $q$. ┬
ёшыє ╥хюЁхь√ \ref{tm:Riesz0} шьххь яЁхфёЄртыхэшх фы  Ёхчюы№тхэЄ√
$$
R_0(\lambda) =\sum_{-\infty}^\infty \frac {\langle\cdot , \bold z_n\rangle\bold y_n}{\lambda -\lambda_n}.
$$

╬сючэрўшь $ f_k =\langle\bold f, \bold z_n\rangle, \ g_k =\langle\bold z_n, \bold g\rangle$.  ╥юуфр (фы  ъЁрЄъюёЄш ь√ яш°хь $\|\cdot\|_p $  тьхёЄю
$\|\cdot\|_{L_p}$ ш яюы№чєхьё  эхЁртхэёЄтюь ├хы№фхЁр ш ЄхюЁхьющ ╒рєёфюЁфр-▐эур)
\begin{multline}
\|R_0(\lambda)\|_{L_p\to L_q} = \sup_{\|f\|_p=1, \|g\|_{q'}=1} |\langle R_0(\lambda) \bold f, \bold g\rangle|\le
 \\
 \le \sum \frac{|f_n g_n|}{|\lambda -\lambda_n|} \le \left(\sum
|f_n|^{p'}\right)^{1/p'} \left(\sum\left( \frac {|g_n|}{|\lambda - \lambda_n|} \right)^p\right)^{1/p}
\\
\le\|f\|_p \left(\sum |g_n|^q\right)^{\frac 1p\, \frac pq} \left(\sum\left(\frac 1{|\lambda-\lambda_n|^p}\right)^{\left(\frac
qp\right)'}\right)^{\frac{\beta}p}\le
\\
\le \|f\|_p \|g\|_q \left(\sum \left(\frac 1{|\lambda-\lambda_n|}\right)^\beta\right)^{1/\beta},
\end{multline}
уфх $\beta = p(q/p)' =\frac{pq}{q-p}.$  ╚ч  ЇюЁьєы фы   $\lambda_n$ ыхуъю яюыєўрхь
$$
\frac 1{|\lambda -\lambda_n|} \le \frac C{|\tau| +|\mu- n|},
$$
яю¤Єюьє
$$
\left(\sum \left(\frac 1{|\lambda-\lambda_n|}\right)^\beta\right)^{1/\beta} \le C \left(\int_0^\infty \frac {dx}{(|\tau| +x)^\beta}\right)^{1/\beta}
\le C|\tau|^{1/\beta-1} =C|\tau|^{-1-1/q+1/p}.
$$
┼ёыш $p=2=q$,  Єю фюърчрЄхы№ёЄтю єяЁю∙рхЄё . ┬ ёыєўрх $p=1,$ ш $ 2\le q\le \infty$  ўшёыр $|f_n|$  юЎхэштр■Єё  эхъюЄюЁющ ъюэёЄрэЄющ $C$ ш
фюърчрЄхы№ёЄтю Єюцх єяЁю∙рхЄё . ═ръюэхЎ, т ёыєўрх $p=1$ ш $q=\infty$, эр°р юЎхэър т√ЄхърхЄ шч яЁхфёЄртыхэш  \eqref{Green} Ёхчюы№тхэЄ√ ъръ
шэЄхуЁры№эюую юяхЁрЄюЁр ш ЁртэюьхЁэющ юЎхэъш т эюЁьх $L_\infty$ ЇєэъЎшш ├Ёшэр (ёь ╥хюЁхьє \ref{tm:Green}).

┼ёыш юяхЁрЄюЁ $L_{0,U}$  Ёхуєы Ёхэ, эю эх ёшы№эю Ёхуєы Ёхэ, Єю т
ёыєўрх юЄёєЄёЄтш  яЁшёюхфшэхээ√ї ЇєэъЎшщ фюърчрЄхы№ёЄтю эх
ьхэ хЄё . ┼ёыш цх яЁшёюхфшэхээ√х ЇєэъЎшш хёЄ№, Єю эрфю
тюёяюы№чютрЄ№ё  ёыхфє■∙шь яЁхфёЄртыхэшхь фы  Ёхчюы№тхэЄ√ (ёь.
\cite{Ke})
$$
R_0(\lambda) =\sum_{-\infty}^\infty \frac {\langle\cdot , \bold z^0_n\rangle\bold y^1_n +  \langle\cdot , \bold z^1_n\rangle\bold y^0_n } {\lambda
-\lambda_n}+  \frac {\langle\cdot , \bold z^0_n\rangle\bold y^0_n}{(\lambda -\lambda_n)^2}.
$$
┬ЄюЁр  ёєььр юЎхэштрхЄё  ЄЁштшры№эю, р яхЁтр , ё єўхЄюь
схёёхыхтюёЄш ёшёЄхь ёюсёЄтхээ√ї ш яЁшёюхфшэхээ√ї ЇєэъЎшщ, Єюўэю
Єръцх ъръ т ёшы№эю Ёхуєы Ёэюь ёыєўрх.

\end{proof}

\section{╬яхЁрЄюЁ ─шЁрър ё Ёхуєы Ёэ√ьш ъЁрхт√ьш єёыютш ьш.}

{\bf 4.1} ─ы  шёёыхфютрэш  Ёхуєы Ёэюую юяхЁрЄюЁр $L_{Q,U}$ ь√ эрьхЁхэ√ тюёяюы№чютрЄ№ё  яюыєўхээ√ьш т ярЁруЁрЇх  2  рёшьяЄюЄшърьш. ╬фэръю ¤Єш
рёшьяЄюЄшъш яюыєўхэ√ Єюы№ъю т яЁюшчтюы№эющ яюыюёх $|\text{Im}\lambda|<\alpha,$  тэх яюыюё эрь рёшьяЄюЄшъш эх шчтхёЄэ√. ╧ю¤Єюьє трцэє■ Ёюы№ шуЁрхЄ
ёыхфє■∙шщ Ёхчєы№ЄрЄ, ъюЄюЁ√щ ь√ ёхщўрё яюыєўшь фы  ёыєўр , ъюуфр $L_{0,U}$  ёшы№эю Ёхуєы Ёхэ (яюЄюь ь√ яюърцхь, ўЄю Ёхчєы№ЄрЄ юёЄрхЄё  тхЁэ√ь фы 
яЁюшчтюы№э√ї Ёхуєы Ёэ√ї $L_{0,U}$).

\begin{Theorem}\label{tm:resmain} ╧єёЄ№  $L_{0,U}$ ёшы№эю Ёхуєы Ёхэ ш $Q\in L_1$.
╥юуфр эрщфхЄё   $\alpha >1$,  Єръюх, ўЄю яЁш  $|\text{Im}\lambda|>\alpha$ Ёхчюы№тхэЄр $R(\lambda) = (L_{Q,U}-\lambda)^{-1}$  ёє∙хёЄтєхЄ (ъръ юяхЁрЄюЁ
т $\mathbb H$) ш т√яюыэ хЄё  юЎхэър
$$\|R(\lambda)\| \le  C |\text{Im}\lambda|^{-1}
$$
ё яюёЄю ээющ  $C$, чртшё ∙хщ юЄ $\alpha$, эю эх чртшё ∙хщ юЄ $\lambda$ тэх яюыюё√ $\Pi_\alpha$. ╧Ёш ¤Єюь ёяхъЄЁ юяхЁрЄюЁр $L_{Q,U}$ фшёъЁхЄхэ.
\end{Theorem}
\begin{proof}
{\it ╪ру 1.} ╠√ єцх фюърчрыш, ўЄю Ёхчюы№тхэЄр $R(\lambda_0)=(L_{Q,U}-\lambda_0)^{-1}$  ёє∙хёЄтєхЄ  яЁш эхъюЄюЁюь $\lambda_0\in\mathbb C$, р чэрўшЄ
юсЁрч $R(\lambda_0)$ ёютярфрхЄ ё юсырёЄ№■ юяЁхфхыхэш   юяхЁрЄюЁр $L_{Q,U}$, р яюЄюьє ёюфхЁцшЄё  т яЁюёЄЁрэёЄтх $W_1^1$.  ═ю тыюцхэшх $W_1^1
\hookrightarrow L_2$ ъюьяръЄэю.  ёыхфютрЄхы№эю,   юяхЁрЄюЁ $R(\lambda_0)$ ъюьяръЄхэ,  р ёяхъЄЁ  $L_{Q,U}$ фшёъЁхЄхэ.

\noindent{\it ╪ру 2.} ╬сючэрўшь $R_0(\lambda) = (L_{0,U} -\lambda)^{-1}$  ш чряш°хь фы  Ёхчюы№тхэЄ√ $R(\lambda)= (L_{0,U} + Q -\lambda)^{-1}$
ЇюЁьры№э√щ Ё ф
\begin{equation}\label{form}
R(\lambda) = R_0(\lambda) +R_0(\lambda)QR_0(\lambda) +(R_0(\lambda)Q)^2R_0(\lambda)) +\dots .
\end{equation}
┬ёх ўыхэ√ ¤Єюую Ё фр ъюЁЁхъЄэю юяЁхфхыхэ√, Єръ ъръ юяхЁрЄюЁ $R_0(\lambda): L_1 \to L_\infty$  юуЁрэшўхэ. ═р°р Ўхы№ яюърчрЄ№, ўЄю эюЁь√ ёырурхь√ї т
¤Єюь Ё фх єс√тр■Є ёю ёъюЁюёЄ№■ ухюьхЄЁшўхёъющ яЁюуЁхёёшш ш Ё ф ёїюфшЄё  т эюЁьх $\mathbb H$.

╘шъёшЁєхь $\varepsilon>0$ (ъюЄюЁюх т√схЁхь яючцх) ш эрщфхь юуЁрэшўхээє■ ЇєэъЎш■ $V(x)$, Єръє■, ўЄю
$$
Q_\varepsilon(x) = Q(x) -V(x),\qquad \|Q_\varepsilon (x)\|_1 \le \varepsilon ,\qquad \|V(x)\|_\infty <C_\varepsilon.
$$
╥юуфр $(n+1)$-х ёырурхьюх $\left[R_0(\lambda)(Q_\varepsilon +V)\right]^nR_0(\lambda)$ т Ё фх \eqref{form} чряш°хЄё  т тшфх
$S_n(\lambda)R_0(\lambda)$,  уфх $S_n(\lambda)$ --- ёєььр, яюыєўхээр  яЁш ЁрёъЁ√Єшш ёъюсюъ т т√Ёрцхэшш $\left[R_0(Q_\varepsilon+V)\right]^n$. ▌Єр
ёєььр ёюёЄюшЄ шч $2^n$ ёырурхь√ї (юсючэрўшь шї ўхЁхч $T_{n,j}$, $j=1,\dots,2^n$), ърцфюх шч ъюЄюЁ√ї  ты хЄё  яЁюшчтхфхэшхь $n$ ьэюцшЄхыхщ тшфр
$R_0Q_\varepsilon$ шыш $R_0V$.

\noindent{\it ╪ру 3.} ╬Ўхэшь эюЁьє ърцфюую ёырурхьюую $\|T_{n,j}\|_{L_\infty\to L_2}$. ─ы  ¤Єюую тюёяюы№чєхьё  ╥хюЁхьющ \ref{th:resolvet0}, ёюуырёэю
ъюЄюЁющ эрщфхЄё  ўшёыю $\alpha_0$,  Єръюх, ўЄю яЁш $\tau = |\text{Im}\lambda|\ge \alpha_0$ ёяЁртхфышт√ юЎхэъш
\begin{gather*}
\|R_0(\lambda)\|_{L_2\to L_2} <C\tau^{-1}, \qquad \|R_0(\lambda)\|_{L_1\to
L_2}<C\tau^{-1/2}\\
\|R_0(\lambda)\|_{L_2\to L_\infty} <C\tau^{-1/2}, \qquad \|R_0(\lambda)\|_{L_1\to L_\infty} <C.
\end{gather*}
┬√схЁхь ЄхяхЁ№ ўшёыю $\varepsilon$,  Єръ, ўЄюс√ $\delta = \varepsilon C <1/4$,  р ўшёыю $\alpha>\alpha_0$  т√схЁхь Єръ, ўЄюс√
\begin{equation}\label{alpha}
(\pi)^{1/2} C_\varepsilon C \alpha^{-1/2} <\delta.
\end{equation}
┬ё■фє фрыхх юЎхэъш яЁютюфшь яЁш ЇшъёшЁютрээюь $\lambda$ ш $\tau = |\text{Im}\lambda|> \alpha$, уфх $\alpha$ єфютыхЄтюЁ хЄ \eqref{alpha}. ╤ єўхЄюь
эхЁртхэёЄт $\|Q_\varepsilon\|_1 <\varepsilon, \ \, \|V\|_\infty < C_\varepsilon $ яюыєўрхь
$$
\|Q_\varepsilon\|_{L_\infty\to L_1}<\varepsilon,\qquad\|V\|_{L_\infty\to L_2}\le\pi^{1/2}C_\varepsilon.
$$
╚ч яЁштхфхээ√ї юЎхэюъ фы  Ёхчюы№тхэЄ√ $R_0(\lambda)$ юЄё■фр ёыхфє■Є эхЁртхэёЄтр
\begin{gather*}
\|R_0Q_\varepsilon\|_{L_\infty \to L_\infty} \le \|Q_\varepsilon\|_{L_\infty \to L_1}\ \,
\|R_0\|_{L_1\to L_\infty} \le \varepsilon C\le\delta,\\
\|R_0Q_\varepsilon\|_{L_\infty \to L_2} \le \|Q_\varepsilon\|_{L_\infty \to L_1}\ \,
\|R_0\|_{L_1\to L_2} \le \varepsilon C\tau^{-1/2}\le\delta\tau^{-1/2}\\
\|R_0V\|_{L_\infty \to L_\infty} \le\|R_0\|_{L_2\to L_\infty}\|V\|_{L_\infty\to L_2}\le
C\tau^{-1/2}\pi^{1/2}C_\varepsilon\le\delta,\\
\|R_0V\|_{L_\infty \to L_2} \le\|R_0\|_{L_2\to L_2}\|V\|_{L_\infty\to L_2}\le C\tau^{-1}\pi^{1/2}C_\varepsilon\le\delta\tau^{-1/2}.
\end{gather*}
╥хяхЁ№ ь√ ьюцхь юЎхэшЄ№ $\|T_{n,j}\|_{L_\infty\to L_2}$. ▌Єю яЁюшчтхфхэшх ёюёЄртыхэю шч $n$ ьэюцшЄхыхщ $T_{n,j}=P_1P_2\dots P_n$, ърцф√щ шч ъюЄюЁ√ї
Ёртхэ ышсю $R_0Q_\varepsilon$, ышсю $R_0V$. ╥юуфр
$$
\|T_{n,j}\|_{L_\infty\to L_2}\le\|P_1\|_{L_\infty\to L_2}\|P_2\|_{L_\infty\to L_\infty}\dots\|P_n\|_{L_\infty\to L_\infty}\le
\delta\tau^{-1/2}\cdot\delta^{n-1}=\delta^n\tau^{-1/2}.
$$
{\it ╪ру 4.} ╤єььр $S_n$ ёюёЄртыхэр шч $2^n$ юяхЁрЄюЁют $T_{n,j}$, р чэрўшЄ
$$
\|S_n\|_{L_\infty\to L_2}\le2^n\delta^n\tau^{-1/2}<2^{-n}\tau^{-1/2}.
$$
╥юуфр ърцфюх ёырурхьюх Ё фр \eqref{form} юЎхэштрхЄё 
$$
\|S_nR_0(\lambda)\|_{L_2\to L_2}\le\|S_n\|_{L_\infty\to L_2}\|R_0\|_{L_2\to L_\infty}\le2^{-n}\tau^{-1/2}C\tau^{-1/2}=2^{-n}C\tau^{-1}.
$$
 ╥хь ёрь√ь, Ё ф \eqref{form} ёїюфшЄё  ш юЎхэштрхЄё  тхышўшэющ
 $2C\tau^{-1}.$  ╥хюЁхьр фюърчрэр.
\end{proof}

{\bf 4.2} ─рыхх ттхфхь юсючэрўхэшх
$$
E(x)=\frac 12
\exp\left\{\frac12\int_0^x(q_3(t)-q_2(t))dt\right\},\qquad\text{ш}\
E=E(\pi).
$$
 ▌Єр ЇєэъЎш  шуЁрхЄ трцэє■ Ёюы№. ┬ ёыхфє■∙хщ ЄхюЁхьх ь√ яюърцхь,
 ўЄю яЁш фюсртыхэшш ъ яЁюёЄхщ°хьє юяхЁрЄюЁє $L_{0,U}$  яюЄхэЎшрыр
 $ Q$ ёюсёЄтхээ√х чэрўхэш  тючьє∙хээюую юяхЁрЄюЁр $L_{Q,U}$
  рёшьяЄюЄшўхёъш ёютярфр■Є ё фтєь  (шыш юфэющ) ёхЁш ьш ёюсёЄтхээ√ї чэрўхэшщ
эхтючьє∙хээюую юяхЁрЄюЁр, хёыш  Єюы№ъю $E=1$,  т ўрёЄэюёЄш,  хёыш
$q_2=q_3$. ═ю яЁш $E\ne 1$  ¤Єю эх Єръ, ёюсёЄтхээ√х чэрўхэш 
юяхЁрЄюЁр $L_{Q,U}$  рёшьяЄюЄшўхёъш ёсышцр■Єё  ё ёюсёЄтхээ√ьш
чэрўхэш ьш фЁєуюую яЁюёЄхщ°хую юяхЁрЄюЁр $L_{0,U'}$, яюЁюцфхээюую
ъЁрхт√ь єёыютшхь $U'(\bold y) =0$,  ъюЄюЁюьє юЄтхўрхЄ Ёрё°шЁхээр 
ьрЄЁшЎр
 \begin{equation}\label{prime}
\mathcal U' = \begin{pmatrix} u_{11}& u_{12} & Eu_{13} & Eu_{14} \\
u_{21}& u_{22} & Eu_{23} &
 Eu_{24}
 \end{pmatrix}
 \end{equation}
┬ ўрёЄэюёЄш, ёюсёЄтхээ√х чэрўхэш  юяхЁрЄюЁр  $L_{Q,U}$ рёшьяЄюЄшўхёъш яЁюёЄ√х ш ЁртэюьхЁэю юЄфхыхэ√ фЁєу юЄ фЁєур, хёыш ш Єюы№ъю хёыш юяхЁрЄюЁ
$L_{0,U'}$  шьххЄ фтх эхёютярфр■∙шї ёхЁшш ёюсёЄтхээ√ї чэрўхэшщ. ▌Єю юяЁртф√трхЄ ёыхфє■∙хх юяЁхфхыхэшх: \textit{╬яхЁрЄюЁ  $L_{Q,U}$ эрчютхь ёшы№эю
Ёхуєы Ёэ√ь, хёыш ёшы№эю Ёхуєы Ёхэ юяхЁрЄюЁ  $L_{0,U'}$ яюЁюцфхээ√щ ъЁрхт√ь єёыютшхь $U'(\bold y) =0$,  ъюЄюЁюьє юЄтхўрхЄ Ёрё°шЁхээр  ьрЄЁшЎр
\eqref{prime}. }  ╚ч юяЁхфхыхэшщ ёыхфєхЄ, ўЄю юяхЁрЄюЁ√ $L_{Q,U}$ ш $L_{Q,U'}$ Ёхуєы Ёэ√ шыш эх Ёхуєы Ёэ√ юфэютЁхьхээю (яю¤Єюьє ьюцэю уютюЁшЄ№ ю
Ёхуєы ЁэюёЄш ъЁрхт√ї єёыютшщ фы  юяхЁрЄюЁр ─шЁрър ё яЁюшчтюы№э√ь яюЄхэЎшрыюь),  эю ёшы№эр  Ёхуєы ЁэюёЄ№ чртшёшЄ ъръ юЄ ъЁрхт√ї єёыютшщ, Єръ ш юЄ
яюЄхэЎшрыр.

 ╟рьхЄшь, хёыш  ъЁрхтюх єёыютшх яюЁюцфрхЄё  ьрЄЁшЎхщ
\eqref{prime},  Єю ёюуырёэю ╥хюЁхьх 3.1  юяхЁрЄюЁ $L_{0,U'}$
 шьххЄ фтх ёхЁшш яЁюёЄ√ї ёюсёЄтхээ√ї чэрўхэшщ (т ёыєўрх ёютярфхэш  юэш юсЁрчє■Є фтєъЁрЄэ√х ёюсёЄтхээ√х чэрўхэш ), шьх■∙шї рёшьяЄюЄшъє
\begin{equation}\label{eigen'}
\lambda^0_n=\kappa_j+n,
\end{equation}
уфх $\kappa_j =\kappa_{j(n)}=j -i\pi^{-1}\ln z_j$, уфх $j=0$, хёыш $n$ ўхЄэю ш $j=1$, хёыш $n$ эхўхЄэю,  р ўшёыр $z_0$, $z_1$ ёєЄ№ ъюЁэш ътрфЁрЄэюую
єЁртэхэш 
\begin{equation}\label{sqeq'}
E[J_{14}-J_{23}-i(J_{13}+J_{24})]z^2+[J_{12}+E^2
J_{34}]z+E[J_{14}-J_{23}+i(J_{13}+J_{24})]=0.
\end{equation}
╤ўшЄрхь, ўЄю тхЄтш фюурЁшЇьр т юяЁхфхыхэшш ўшёхы $\kappa_j$ т√сЁрэ√ Єръ, ўЄю  $\text{Re}\kappa_j\in(-1,1]$. ┬√схЁхь $\varepsilon >0$ Єръ, ўЄюс√
чрьъэєЄ√х ъЁєцъш Ёрфшєёр $\varepsilon$  ё ЎхэЄЁрьш т ёюсёЄтхээ√ї чэрўхэш ї \eqref{eigen'} юяхЁрЄюЁр $L_{0,U'}$  эх яхЁхёхърышё№ (т ёыєўрх Ёхуєы Ёэ√ї,
эю эх ёшы№эю Ёхуєы Ёэ√ї єёыютшщ ъЁєуш ё эюьхЁрьш $2n$ ш $2n+1$, $n\in\mathbb Z$, ёютярфр■Є). ┬ёяюьэшь юяЁхфхыхэшх ЇєэъЎшш $\Upsilon(x,\lambda)$  шч
ярЁруЁрЇр 2 ш фы  $Q\in L_p$, $p\in[1,2]$, юяЁхфхышь ўшёыр
\begin{equation}\label{s-n}
s_n(r) = \max_{|\lambda -\lambda_n^0|\le r} \left(\Upsilon(\pi,\lambda)+\Upsilon_{p'}(\lambda)\right),\quad 0\le r\le\varepsilon,
\end{equation}
уфх $1/p'+1/p=1$. ╚ч юяЁхфхыхэш  ш ╦хьь√ \ref{lem:Ups1} ёыхфєхЄ
\begin{Lemma}\label{lem:snto0}
╩ръют с√ эш с√ы ёєььшЁєхь√щ
 яюЄхэЎшры $Q$ ўшёыр $s_n(\varepsilon)\to 0$  яЁш $n\to \pm\infty$
\end{Lemma}
╠√ уюЄютшьё  фюърчрЄ№ яхЁт√х Ёхчєы№ЄрЄ√ юс рёшьяЄюЄшўхёъюь яютхфхэшш ёюсёЄтхээ√ї чэрўхэшщ юяхЁрЄюЁр $L_{Q,U}$ ё Ёхуєы Ёэ√ьш єёыютш ьш. ─ы  яЁютхфхэш 
юЎхэюъ эрь сєфхЄ єфюсэю ттхёЄш ёыхфє■∙шх юсючэрўхэш . ╧єёЄ№
$$
\mathbf{c}(x,\lambda) =(c_1(x,\lambda),c_2(x,\lambda))^t,\qquad \mathbf{s}(x,\lambda)= (s_1(x,\lambda),s_2(x,\lambda))^t
$$
--- Ёх°хэш  єЁртэхэш  \(l(y)=\lambda y\), юяЁхфхыхээ√х т ярЁруЁрЇх 2.  ╬с∙хх
Ёх°хэшх ¤Єюую єЁртэхэш  шьххЄ тшф \(\mathbf{y}=\gamma_1\mathbf{c}+\gamma_2\mathbf{s}\), уфх \(\gamma_1\), \(\gamma_2\) --- яюёЄю ээ√х. ╧юфёЄрты   ¤Єю
Ёх°хэшх т ъЁрхтюх єёыютшх  \(U(\mathbf{y})=0\), яюыєўрхь
\begin{equation}\label{eigeneq'}
\begin{cases}[u_{12}+u_{13} s_1(\pi, \lambda) +u_{14}s_2(\pi,\lambda)]\gamma_1+[u_{11}+u_{13}c_1(\pi,\lambda)-u_{14}c_2(\pi,\lambda)]\gamma_2=0,\\
[u_{22}+u_{23}s_1(\pi,\lambda)+u_{24}s_2(\pi,\lambda)]\gamma_1+[u_{21}+u_{23}s_2(\pi,\lambda)-u_{24}s_1(\pi,\lambda)]\gamma_2=0.\end{cases}
\end{equation}
╬сючэрўшь  $2\times 2$ ьрЄЁшЎє, яюЁюцфр■∙є■ ¤Єє юфэюЁюфэє■ ёшёЄхьє ышэхщэ√ї єЁртэхэшщ ўхЁхч $\mathcal M(\lambda)$. ╫шёыю $\lambda\in\mathbb C$
 ты хЄё  ёюсёЄтхээ√ь чэрўхэшхь юяхЁрЄюЁр $L_{Q,U}$ Єюуфр ш Єюы№ъю Єюуфр, ъюуфр юяЁхфхышЄхы№ $\Delta(\lambda): =\det \mathcal M(\lambda)=0$. ╧юёых
эхёыюцэ√ї т√ўшёыхэшщ яЁшїюфшь ъ ёыхфє■∙хьє т√Ёрцхэш■ фы  їрЁръЄхЁшёЄшўхёъюую юяЁхфхышЄхы 
\begin{equation}\label{eq:1.5A}
    \Delta(\lambda)=J_{12}+J_{34}+J_{13}s_1(\pi,\lambda)+J_{14}s_2(\pi,\lambda)+J_{32}c_1(\pi,\lambda)+J_{42}c_2(\pi,\lambda).
\end{equation}
╧юфёЄрты   ё■фр рёшьяЄюЄшўхёъшх ЇюЁьєы√ \eqref{cosas} ш \eqref{sinas}, яюыєўшь яЁхфёЄртыхэшх ЇєэъЎшш $\Delta(\lambda)$  т яюыюёх $\Pi_\alpha$
\begin{gather}
\Delta(\lambda)=\Delta_0(\lambda)+\delta(\lambda),\quad\text{уфх}\label{Deltaas}\\
\Delta^0(\lambda)=J_{12}+J_{34}+\frac{E(\pi)e^{i\pi\lambda}}2\left[J_{14}+J_{32}-i(J_{42}-J_{13})\right]+
\frac{E(\pi)e^{-i\pi\lambda}}2\left[J_{14}+J_{32}+i(J_{42}-J_{13})\right].\notag
\end{gather}
╠√ тшфшь, ўЄю юяЁхфхышЄхы№ $\Delta^0(\lambda)$  ёютярфрхЄ  ё їрЁръЄхЁшёЄшўхёъшь юяЁхфхышЄхыхь рёёюЎшшЁютрээюую юяхЁрЄюЁр $L_{0,U'}$. ─рыхх, шч
\eqref{cosas} ш \eqref{sinas} ёыхфєхЄ, ўЄю яЁш тёхї $\lambda\in\Pi_\alpha$, фы  ъюЄюЁ√ї
$$
\Upsilon(\lambda)<1/(8k^4),\qquad \text{уфх}\ \ k=2+12R\ch(2\pi\alpha+1),\ R=\|Q\|_{L_1},
$$
т√яюыэхэр юЎхэър
\begin{equation}\label{Deltaas1}
|\delta(\lambda)|\le M(R,\alpha)\Big(\Upsilon(\pi,\lambda)+\Upsilon_{p'}(\lambda)\Big)\le (1+\pi^{1/p'})M(R,\alpha)\Upsilon(\lambda)
\end{equation}
ё эхъюЄюЁющ ъюэёЄрэЄющ $M$, юяЁхфхы хьющ яю $R$ ш $\alpha$. ╠√ єцх юЄьхўрыш, ўЄю ЇєэъЎш  $\Upsilon(\lambda)\to0$ яЁш $\text{Re}\lambda\to\pm\infty$,
р чэрўшЄ эрщфхЄё  Єръюх ўшёыю $\beta>0$, ўЄю т юсырёЄш $P_{\alpha,\beta}=\{\lambda\in\mathbb C\vert\ |\text{Im}\lambda|<\alpha,\
|\text{Re}\lambda|>\beta\}$ т√яюыэхэю яЁхфёЄртыхэшх \eqref{Deltaas} ё юЎхэъющ \eqref{Deltaas1}.

╬яЁхфхышЄхы№ $\Delta^0$ Ёртхэ эєы■ яЁш $\lambda=\lambda_n^0$. ┼ёыш ъЁрхт√х єёыютш  єёшыхээю Ёхуєы Ёэ√, Єю тёх хую эєыш яЁюёЄ√ ш Єюуфр т ъЁєух
$|\lambda-\lambda_n^0|<\varepsilon$, уфх ўшёыю $\varepsilon$ юяЁхфхыхэю т√°х, т√яюыэхэр юЎхэър $|\Delta^0(\lambda)|>c_1|\lambda-\lambda_n^0|$. ╧Ёш
¤Єюь, яюёъюы№ъє ЇєэъЎш  $\Delta^0(\lambda)$ яхЁшюфшўэр, ўшёыю $c_1$ ьюцэю т√сЁрЄ№ эх чртшё ∙шь юЄ $n$. ╬ўхтшфэю, ўЄю ЇєэъЎш  $s_n(r)$ эхяЁхЁ√тэр ш
ьюэюЄюээю тючЁрёЄрхЄ яю $r\in[0,\varepsilon]$. ┬√схЁхь ўшёыю $N=N(Q,U)$ Єръшь, ўЄю яЁш тёхї $|n|\ge N$ ъЁєуш $|\lambda-\lambda_n^0|\le\varepsilon$
ыхцрыш т юсырёЄш $\Pi_{\alpha,\beta}$, р ўшёыр $s_n(\varepsilon)$ эх яЁхтюёїюфшыш $c_1\varepsilon/M$ (юэш ёЄЁхь Єё  ъ эєы■, ёюуырёэю ╦хььх
\ref{lem:snto0}). ╥юуфр фы  ърцфюую Єръюую $n$ эрщфхЄё  ўшёыю $r=r_n$,  ты ■∙ххё  ъюЁэхь єЁртэхэш  $s_n(r)=c_1r/M$. ╧Ёш эхъюЄюЁ√ї $n$ ¤Єю єЁртэхэшх
ьюцхЄ шьхЄ№ эхёъюы№ъю Ёх°хэшщ --- ь√ т√схЁхь эршьхэ№°хх шч эшї (юэю, тючьюцэю, Ёртэю эєы■). ╥ръшь юсЁрчюь
$$
\max_{|\lambda-\lambda_n^0|\le r_n}\Big(\Upsilon(\pi,\lambda)+\Upsilon_{p'}(\lambda)\Big)=s_n(r_n)=\frac{c_1r_n}{M}.
$$
╫шёыр $s_n(r_n)$ ь√ сєфхь юсючэрўрЄ№ яЁюёЄю $s_n$.

┼ёыш цх ъЁрхт√х єёыютш  Ёхуєы Ёэ√, эю эх єёшыхээю Ёхуєы Ёэ√, Єю тёх эєыш $\lambda_n^0$ фтєъЁрЄэ√ ш т ъЁєух $|\lambda-\lambda_n^0|<\varepsilon$
т√яюыэхэр юЎхэър $|\Delta^0(\lambda)|>c_2|\lambda-\lambda_n^0|^2$, уфх $c_2$ тэют№ эх чртшёшЄ юЄ $n$. ┬ ¤Єюь ёыєўрх ь√ Єръцх ЁрёёьюЄЁшь ЇєэъЎшш
$s_n(r)$, эю ўшёыю $r_n$ юяЁхфхышь ъръ эршьхэ№°шщ ъюЁхэ№ єЁртэхэш  $s_n(r)=c_2r^2/M$, чртхфюью ёє∙хёЄтє■∙шщ эр юЄЁхчъх $0\le r\le\varepsilon$, хёыш
$s_n(\varepsilon)\le c_2\varepsilon^2/M$, ўЄю ёяЁртхфыштю яЁш фюёЄрЄюўэю сюы№°шї $n$.

\begin{Theorem}\label{tm:4.1}
╧єёЄ№ $Q\in L_1$ ш $L_{Q,U}$ ---  Ёхуєы Ёэ√щ юяхЁрЄюЁ ─шЁрър. ╬сючэрўшь  ўхЁхч $\{\lambda_n^0\}$  ёюсёЄтхээ√х чэрўхэш  \eqref{eigen'}
рёёюЎшшЁютрээюую юяхЁрЄюЁр $L_{0,U'}$  ш  ўхЁхч $\lambda_n$ ёюсёЄтхээ√х чэрўхэш  юяхЁрЄюЁр $L_{Q,U}$ ё єўхЄюь рыухсЁршўхёъющ ъЁрЄэюёЄш.  ╥юуфр яЁш
яюфїюф ∙хщ эєьхЁрЎшш яюёыхфютрЄхы№эюёЄш $\{\lambda_n\}$ (ш Єрър  эєьхЁрЎш  тючьюцэр) т ёыєўрх ёшы№эющ Ёхуєы ЁэюёЄш юяхЁрЄюЁр $L_{Q,U}$   ёяЁртхфыштр
рёшьяЄюЄшър
\begin{equation}\label{asymreg} |\lambda_n -\lambda_n^0| \le r_n\le \frac{M(R,\alpha)}{c_1(U)}s_n(\varepsilon),
\end{equation}
р т ёыєўрх юс√ўэющ Ёхуєы ЁэюёЄш
\begin{equation}\label{asymnonreg}
|\lambda_n -\lambda_n^0|^2\le r_n\le\sqrt{\frac{M(R,\alpha)}{c_2(U)}s_n(\varepsilon)},
\end{equation}
уфх эхЁртхэёЄтр ёяЁртхфышт√, эрўшэр  ё эхъюЄюЁюую эюьхЁр $N$, чртшё ∙хую юЄ яюЄхэЎшрыр $Q$ ш ъЁрхт√ї єёыютшщ $U$.
\end{Theorem}
\begin{proof}
╤эрўрыр ЁрёёьюЄЁшь ёыєўрщ, ъюуфр ъюЁэш $\Delta^0(\lambda)$ яЁюёЄ√х. ═р юъЁєцэюёЄ ї $\gamma_n$ Ёрфшєёр $r_n$ ё ЎхэЄЁрьш т эєы ї $\lambda_n^0$ шьххь
юЎхэъш
$$
|\Delta^0(\lambda)| > c_1r_n,\qquad |\delta(\lambda)| < M\Big(\Upsilon(\pi,\lambda)+\Upsilon_{p'}(\lambda)\Big)\le c_1r_n.
$$
╥юуфр эр юъЁєцэюёЄ ї $\gamma_n$ сєфхь шьхЄ№ юЎхэъє $|\Delta^0(\lambda)|>|\delta(\lambda)|$. ╤юуырёэю ЄхюЁхьх ╨є°х, тэєЄЁш юъЁєцэюёЄш $\gamma_n$
эрїюфшЄё  Ёютэю юфшэ эєы№ $\lambda_n$ ЇєэъЎшш $\Delta(\lambda)$, Є.х. яЁш сюы№°шї $|n|$ т√яюыэ ■Єё  юЎхэъш $|\lambda_n-\lambda_n^0| < r_n$. ┬ ёыєўрх
фтєъЁрЄэ√ї ъюЁэхщ  эр юъЁєцэюёЄ ї $\gamma_n$ ёяЁртхфышт√ юЎхэъш
$$
|\Delta^0(\lambda)| >c_2r_n^2,\qquad |\delta(\lambda)| < M\Big(\Upsilon(\pi,\lambda)+\Upsilon_{p'}(\lambda)\Big)\le c_2r_n^2,
$$
ш шч ЄхюЁхь√ ╨є°х ёыхфєхЄ, ўЄю тэєЄЁш ърцфющ юъЁєцэюёЄш $\gamma_n$ ыхцрЄ Ёютэю фтр эєы  ЇєэъЎшш $\Delta(\lambda)$.

─ы  чртхЁ°хэш  фюърчрЄхы№ёЄтр ЄхюЁхь√ юёЄрхЄё  яюърчрЄ№, ўЄю ўшёыю ъюЁэхщ ЇєэъЎшщ $\Delta^0(\lambda) $ ш $\Delta(\lambda)$ т ъЁєух фюёЄрЄюўэю
сюы№°юую Ёрфшєёр ёютярфрхЄ. ╠√ єцх чэрхь, ўЄю  яЁш фюёЄрЄюўэю сюы№°юь $\alpha_0$ тэх яюыюё√ $|\text{Im}\lambda|<\alpha_0$ юсх ЇєэъЎшш ъюЁэхщ эх
шьх■Є. ▌Єю ётющёЄтю ёюїЁрэ хЄё  фы  ёхьхщёЄтр юяхЁрЄюЁют $L(t) =L_{tQ,U}$, яюЁюцфхээ√ї яюЄхэЎшрырьш $tQ$,  уфх $0\le t\le 1$. ─хщёЄтшЄхы№эю, ёюуырёэю
╥хюЁхьх \ref{tm:resmain}, °шЁшэр яюыюё√ --- ўшёыю $\alpha_0$, юяЁхфхы хЄё  эхЁртхэёЄтюь
$$
 \alpha_0^{1/2}>\frac{C\|V\|_\infty}{\|Q_\varepsilon\|_1},\qquad\text{уфх}\ \ Q=Q_\varepsilon+V, \qquad\|Q_\varepsilon\|_1<\varepsilon.
$$
╧єёЄ№ ЇєэъЎшш $Q_\varepsilon$ ш $V$ єцх яюёЄЁюхэ√ фы  яюЄхэЎшрыр $Q$. ╧ЁхфёЄртшь яюЄхэЎшры $tQ$ т тшфх $tQ=tQ_\varepsilon+tV$ --- Єюуфр
$\|tQ_\varepsilon\|_1<\varepsilon$, р ўшёыю $\alpha_0$ эх чртшёшЄ юЄ $t$. ╘єэъЎш  $\Upsilon(\lambda)$, юяЁхфхыхээр  \eqref{gamma} Єръцх чртшёшЄ юЄ
яюЄхэЎшрыр $Q$, эю $\Upsilon(tQ,\lambda)= t\Upsilon(Q,\lambda)$. ╧єёЄ№ $\Gamma$
 --- яЁ ьюєуюы№эшъ, уюЁшчюэЄры№э√х ёЄюЁюэ√ ъюЄюЁюую ыхцрЄ эр яЁ ь√ї $|\text{Im}\lambda|=\pm \alpha_0$,
р тхЁЄшъры№э√х ёЄюЁюэ√ т√сЁрэ√ ёЄюы№ фрыхъшьш, ўЄю эр эшї т√яюыэ хЄё  юЎхэър $|\delta(\lambda)|< M\Upsilon(\lambda) < |\Delta^0(\lambda)|$. ╥юуфр яЁш
$t\in [0,1]$ їрЁръЄхЁшёЄшўхёъшщ юяЁхфхышЄхы№ $\Delta_t(\lambda)$, юЄтхўр■∙шщ яюЄхэЎшрыє $tQ$,  эх шьххЄ эєыхщ эр ёЄюЁюэрї $\Gamma$, яюёъюы№ъє
$|\delta_t(\lambda)| < M t \Upsilon(\lambda) <|\Delta^0(\lambda)|$.  ╤ыхфютрЄхы№эю, эєыш $\lambda_n(t)$ Ўхыющ ЇєэъЎшш $\Delta_t(\lambda)$,   ты  ё№
эхяЁхЁ√тэ√ьш ЇєэъЎш ьш юЄ $t\in[0,1]$, эх яхЁхёхър■Є уЁрэшЎ√ $\Gamma$,  р яюЄюьє шї ўшёыю яЁш $t=0$  ш $t=1$  юфшэръютю (т ёыхфє■∙хщ ЄхюЁхьх ь√
ёЄЁюую юсюёэєхь эхяЁхЁ√тэюёЄ№ ёюсёЄтхээ√ї чэрўхэшщ юЄ ярЁрьхЄЁр $t$). ╤ыхфютрЄхы№эю, "<сышчъшх">  ёюсёЄтхээ√х чэрўхэш  $\lambda_n$ ьюцэю чрэєьхЁютрЄ№
Єръ, ўЄюс√ ЇюЁьєы√ \eqref{asymreg}  ш \eqref{asymnonreg}  шьхыш ёь√ёы схч "<ёсю "> т эєьхЁрЎшш. ╥хюЁхьр фюърчрэр.
\end{proof}

╩юэхўэю, трцэр шэЇюЁьрЎш  ю Єюь, ъръ єс√тр■Є ўшёыр $r_n$, їрЁръЄхЁшчє■∙шх ёъюЁюёЄ№ яЁшсышцхэш  ёюсёЄтхээ√ї чэрўхэшщ $\lambda_n$ ъ ёюсёЄтхээ√ь
чэрўхэш ь рёёюЎшшЁютрээюую юяхЁрЄюЁр ё эєыхт√ь яюЄхэЎшрыюь. ┬ ёыєўрх $Q\in L_1$  ь√, ёюуырёэю ╦хььх \ref{lem:snto0}, ьюцхь єЄтхЁцфрЄ№, ўЄю $r_n\to
0$. ╬фэръю єцх т ёыєўрх $Q\in L_p$, $p\in(1,2]$,   ь√ ьюцхь яюыєўшЄ№ сюыхх Єюўэє■ шэЇюЁьрЎш■. ─ы  ¤Єюую эрь яюЄЁхсєхЄё  тёяюьюурЄхы№э√щ Ёхчєы№ЄрЄ,
ъюЄюЁ√щ  ты хЄё  ёыхфёЄтшхь ЄхюЁхь√ ╩рфхЎр (ёь. \cite{Kad},  р Єръцх \cite{Hry10} фы  сюыхх яюыэющ шэЇюЁьрЎшш).

\begin{Lemma}\label{Kadets} ╧єёЄ№ $|\lambda_n-2n| <\varepsilon < 1/(2p)$  яЁш тёхї $n\in \mathbb Z$ ш эхъюЄюЁюь $p\in[1,2]$.
╥юуфр фы  тёхї $f\in L_p[0,\pi]$ ёяЁртхфыштр юЎхэър
$$
\left(\sum_{n\in\mathbb Z} \left|\int\limits_0^\pi f(x) e^{i\lambda_n x}\,dx\right|^{p'}\right)^{1/p'} \le C \|f\|_{L_p},
$$
уфх яюёЄю ээр  $C$  чртшёшЄ  юЄ $\varepsilon$, эю эх чртшёшЄ юЄ яюёыхфютрЄхы№эюёЄш $\{\lambda_n\}$.
\end{Lemma}
\begin{proof}
╨рёёьюЄЁшь ышэхщэ√щ юяхЁрЄюЁ $T:\,L_p\to l_{p'}$,
$$
Tf=\left\{\int\limits_0^\pi f(x)e^{i\lambda_n x}\,dx\right\}_{n\in\,\mathbb{Z}}.
$$
╟рьхЄшь, ўЄю яЁш $p=2$ ш $p=\infty$ ¤ЄюЄ юяхЁрЄюЁ юуЁрэшўхэ, яЁшўхь хую эюЁьр юуЁрэшўхэр ъюэёЄрэЄрьш $M_2$ ш $M_\infty$, эх чртшё ∙шьш юЄ
яюёыхфютрЄхы№эюёЄш $\{\lambda_n\}$ (чфхё№ ь√ шёяюы№чєхь ЄхюЁхьє ╩рфхЎр). ╧Ёшьхэшт шэЄхЁяюы Ўшюээє■ ЄхюЁхьє ╨шёёр--╥юЁшэр (ёь., эряЁшьхЁ, \cite{BL}),
яюыєўшь єЄтхЁцфхэшх ыхьь√.
\end{proof}
╥хяхЁ№ яюыєўшь юЎхэъє фы  ўшёхы $\{r_n\}$. ╬ЄьхЄшь, ўЄю $r_n=Ms_n/c_1\le Ms_n(\varepsilon)/c_1$ (фы  єёшыхээю Ёхуєы Ёэюую юяхЁрЄюЁр) ш
$r_n=\sqrt{Ms_n/c_2}\le\sqrt{Ms_n/c_2}$ (фы  Ёхуєы Ёэюую, эю эх єёшыхээю Ёхуєы Ёэюую юяхЁрЄюЁр). ╥ръшь юсЁрчюь, юЎхэъш, яюыєўхээ√х фы  ўшёхы
$s_n(\varepsilon)$, тыхъєЄ юЎхэъш фы  $r_n$.
\begin{Theorem}  ─ы  ы■сюую $Q\in L_1$ ўшёыр $s_n\to0$ яЁш $n\to\infty$. ┼ёыш цх яюЄхэЎшры $Q$  ыхцшЄ т $L_p$  яЁш $p\in(1,2]$, Єю
$\{s_n\}\in l_{p'}$,  уфх $1/p +1/p' =1$.
\end{Theorem}
\begin{proof}
─ы  $p=1$ ЄхюЁхьр єцх фюърчрэр (ёь. ╦хььє \ref{lem:snto0}). ╧єёЄ№ ЄхяхЁ№ $p\in(1,2]$. ╬сючэрўшь ўхЁхч $\mu_n$ --- Єюўъш, фы  ъюЄюЁ√ї
$$
\Upsilon(\pi,\mu_n)+\Upsilon_{p'}(\mu_n)=\max\limits_{|\lambda-\lambda_n^0|\le\varepsilon}\left(\Upsilon(\pi,\lambda)+\Upsilon_{p'}(\lambda)\right)
$$
╥юуфр $\mu_n=n+\kappa_j+o(1)$. ╟рьхЄшь, хёыш т ╦хььх \ref{Kadets} ўшёыр $2n$ чрьхэшЄ№ эр $2n+\kappa$, Єю юЎхэър ёюїЁрэшЄё , тючьюцэю ё эютющ
ъюэёЄрэЄющ $C$, чртшё ∙хщ Єюы№ъю юЄ $\varepsilon$ ш $\kappa$. ╧юёъюы№ъє тёх ёырурхь√х, тїюф ∙шх т юяЁхфхыхэшх ЇєэъЎшш $\Upsilon(x,\lambda)$,
юЎхэштр■Єё  юфшэръютю, яЁютхфхь Ёрёёєцфхэшх фы  ЇєэъЎшш
$$
\upsilon_0(x,\lambda)=\int_0^x q(t)e^{2i\lambda t}\,dt,
$$
уфх $q\in L_p[0,\pi]$. ┬ ёшыє ЄхюЁхь√ ╩рфхЎр, яЁшьхэхээющ ъ ЇєэъЎшш $q(t)\chi_{[0,x]}(t)$, яЁш фюёЄрЄюўэю сюы№°юь $N$ яюыєўшь
$$
\sum_{|n|>N} |\upsilon_0(x,\mu_n)|^{p'} \le (C\|Q\|_{L_p})^{p'}.
$$
╥юуфр, т ёшыє ЄхюЁхь√ ╦хтш ю яЁхфхы№эюь яхЁхїюфх,
$$
\sum_{|n|>N}\int_0^\pi|\upsilon_0(x,\mu_n)|^{p'}\,dx\le\pi(C\|Q\|_{L_p})^{p'}.
$$
╧Ёшьхэ   х∙х Ёрч ЄхюЁхьє ╩рфхЎр ъ яюёыхфютрЄхы№эюёЄш $\{\upsilon_0(\pi,\mu_n)\}_{|n|>N}$, яюыєўшь, ўЄю
$$
\|\{s_n(\varepsilon)\}_{|n|>N}\|_{l_{p'}}\le (1+\pi^{1/p'})C\|Q\|_{L_p}.
$$
▌Єр юЎхэър тыхўхЄ єЄтхЁцфхэшх ЄхюЁхь√ фы  $p\in(1,2]$.
\end{proof}

{\bf 4.3}. ╥хяхЁ№ ь√ фюърцхь Ёхчєы№ЄрЄ ю эхяЁхЁ√тэющ чртшёшьюёЄш Ёхчюы№тхэЄ√ $(L_{Q,U}-\lambda I)^{-1}$ юЄ яюЄхэЎшрыр $Q\in L_1$ т яЁхфяюыюцхэшш
Ёхуєы ЁэюёЄш юяхЁрЄюЁр $L_{Q,U}$. ╧єёЄ№ $Q_\varepsilon$ --- ёхьхщёЄтю  ёєььшЁєхь√ї ЇєэъЎшщ, єфютыхЄтюЁ ■∙шї єёыютш■
\(\|Q_{\varepsilon}(x)-Q(x)\|_{L_1}\to 0\) яЁш \(\varepsilon\to 0\). ─рыхх яюырурхь,  ўЄю $Q$  ш $U$  ЇшъёшЁютрээ√ьш ш яюырурхь $L= L_{Q,U},
L_\varepsilon = L_{Q_\varepsilon, U}$. ─юърцхь ёыхфє■∙шщ Ёхчєы№ЄрЄ.

\begin{Theorem}\label{tm:1.6} ╧єёЄ№ $K$
ъюьяръЄ т ъюьяыхъёэющ яыюёъюёЄш \(\mathbb C\),  яЁшэрфыхцр∙шщ Ёхчюы№тхэЄэюьє ьэюцхёЄтє  Ёхуєы Ёэюую юяхЁрЄюЁр $L= L_{Q,U}$. ╧єёЄ№
\(\|Q_{\varepsilon}(x)-Q(x)\|_{L_1}\to 0\) яЁш \(\varepsilon\to 0\) ш  $L_\varepsilon = L_{Q_\varepsilon, U}$. ╥юуфр  эрщфхЄё  $\varepsilon_0>0$,
Єръюх, ўЄю яЁш тёхї $\varepsilon<\varepsilon_0$ ъюьяръЄ $K$  ыхцшЄ т Ёхчюы№тхэЄэ√ї ьэюцхёЄтрї юяхЁрЄюЁют $L_\varepsilon$ ш ЁртэюьхЁэю яю $\lambda\in
K$ ёхьхщёЄтю юяхЁрЄюЁют \((L_{\varepsilon}-\lambda)^{-1}\) ёїюфшЄё  яЁш \(\varepsilon\to 0\) т ЁртэюьхЁэющ юяхЁрЄюЁэющ Єюяюыюушш, Є.х.
\begin{equation*}
   \sup_{\lambda\in K} \|(L_{\varepsilon}-\lambda)^{-1}-(L-\lambda)^{-1}\|\to 0
\end{equation*}
\end{Theorem}
\begin{proof}  ─юёЄрЄюўэю фюърчрЄ№ ЄхюЁхьє фы  ЇшъёшЁютрээющ Єюўъш
$\lambda$,  ыхцр∙хщ т Ёхчюы№тхэЄэюь ьэюцхёЄтх $L$. ┬ ёшыє эхяЁхЁ√тэюёЄш   Ёхчюы№тхэЄ√ юЄ $\lambda$  ёїюфшьюёЄ№ ёюїЁрэшЄё  т эхъюЄюЁющ ьрыющ юЄъЁ√Єющ
юъЁхёЄэюёЄш Єюўъш $\lambda$,  р чрЄхь фюърчрЄхы№ёЄтю ьюцэю чртхЁ°шЄ№ ё яюью∙№■ ёЄрэфрЁЄэющ Єхїэшъш т√сюЁр ъюэхўэюую яюфяюъЁ√Єш  ъюьяръЄр $K$.

┬ фюърчрЄхы№ёЄтх сєфхь ёє∙хёЄтхээю шёяюы№чютрЄ№ ╥хюЁхьє~\ref{tm:0}. ╬сючэрўшь ўхЁхч \(\mathbf{c}_{\varepsilon}(x)\), \(\mathbf{s}_{\varepsilon}(x)\)
ярЁє Ёх°хэшщ ёшёЄхь√
\begin{equation}\label{eq:T.1.7}
    B\mathbf{y}'+Q_{\varepsilon}(x)\mathbf{y}=\lambda \mathbf{y},
\end{equation}
єфютыхЄтюЁ ■∙шї єёыютш ь
\[
    \mathbf{c}_{1,\varepsilon}(0)=1,\quad\mathbf{c}_{2,\varepsilon}(0)=0;\quad
    \mathbf{s}_{1,\varepsilon}(0)=0,\quad\mathbf{s}_{2,\varepsilon}(0)=1.
\]
╤юуырёэю ╥хюЁхьх~\ref{tm:0}, шьххь
\begin{equation}\label{eq:bas}
    \|\mathbf{c}_{j,\varepsilon}(x)-\mathbf{c}_{j}(x)\|_{W_1^1}+
    \|\mathbf{s}_{j,\varepsilon}(x)-\mathbf{s}_{j}(x)\|_{W_1^1}\le
    C\|Q_{\varepsilon}(x)-Q(x)\|_{L_1},
    \quad j=1,2.
\end{equation}
╧Ёшьхэшь ьхЄюф  трЁшрЎшш яюёЄю ээ√ї ъ єЁртэхэш■ $l_\varepsilon(\mathbf{y})=\lambda\mathbf{y}+\mathbf{f}$. ╥юуфр Ёх°хэшх яЁшьхЄ тшф
$$
\mathbf{y}_\varepsilon=\gamma_{1,\varepsilon}\mathbf{c}_\varepsilon(x,\lambda)+\gamma_{2,\varepsilon}\mathbf{s}_\varepsilon(x,\lambda)-\int_x^\pi{\mathbf
U}_\varepsilon(x-t,\lambda)\mathbf{f}(t)dt,
$$
уфх
$$
{\mathbf U}_\varepsilon(x,\lambda)=\begin{pmatrix}c_{1,\varepsilon}(x,\lambda)&
s_{1,\varepsilon}(x,\lambda)\\c_{2,\varepsilon}(x,\lambda)&s_{2,\varepsilon}(x,\lambda)\end{pmatrix}
$$
╧ютюЄюЁ   Ёрёёєцфхэш  ╥хюЁхь√ \ref{tm:Green}, яюыєўшь
\begin{multline*}
\begin{pmatrix}\gamma_{1,\varepsilon}\\ \gamma_{2,\varepsilon}\end{pmatrix}=\frac1{\Delta_\varepsilon(\lambda)}
\begin{pmatrix}u_{21}+u_{23}s_{2,\varepsilon}(\pi,\lambda)-u_{24}c_{2,\varepsilon}(\pi,\lambda) & -u_{11}-u_{13}c_{1,\varepsilon}(\pi,\lambda)+u_{14}s_{1,\varepsilon}(\pi,\lambda)\\
-u_{22}+u_{23}c_{2,\varepsilon}(\pi,\lambda)-u_{24}s_{2,\varepsilon}(\pi,\lambda) & u_{12}+u_{13}s_{1,\varepsilon}(\pi,\lambda)+u_{14}c_{1,\varepsilon}(\pi,\lambda)\end{pmatrix}\times\\
\times\begin{pmatrix}u_{11} & u_{12}\\ u_{21} & u_{22}\end{pmatrix}\int_0^\pi{\mathbf U}_\varepsilon(-t,\lambda)\mathbf{f}(t).
\end{multline*}
╧юёъюы№ъє ъЁрхт√х єёыютш  Ёхуєы Ёэ√, Єю
\(\Delta=\Delta_{\varepsilon}(\lambda) \not\equiv 0\) яЁш тёхї
\(\varepsilon\ge 0\). ╚ч юЎхэъш~\eqref{eq:bas} шьххь
\begin{align*}
    \|U(\mathbf{c}_\varepsilon)-U(\mathbf{c}_{\delta})|+|U(\mathbf{s}_{\varepsilon})-U(\mathbf{s}_{\delta})|\le C\|Q_{\varepsilon}-Q_{\delta}\|_{L_1},
\end{align*}
р яюЄюьє \(|\Delta_{\varepsilon}-\Delta_{\delta}|\le C\|Q_{\varepsilon}- Q_{\delta}\|_{L_1}\). ┬√схЁхь ЄхяхЁ№ ўшёыю \(\lambda\) Єръшь, ўЄюс√
\(\Delta_0(\lambda)\neq 0\). ╥юуфр \(|\Delta_{\varepsilon}(\lambda)|>c=c(\lambda)\) яЁш тёхї фюёЄрЄюўэю ьры√ї \(\varepsilon>0\). ╤ыхфютрЄхы№эю,
яюёЄю ээ√х \(\gamma_{1,\varepsilon}\), \(\gamma_{2,\varepsilon}\) Єръют√, ўЄю
\[
    |\gamma_{1,\varepsilon}-\gamma_{1,\delta}|+|\gamma_{2,\varepsilon}-\gamma_{2,\delta}|\le
    C\|Q_{\varepsilon}-Q_{\delta}\|_{L_1}\|f\|_{L_2},
\]
уфх \(C\) чртшёшЄ Єюы№ъю юЄ т√сЁрээюую ўшёыр \(\lambda\).
╧юыєўхээ√х юЎхэъш яюърч√тр■Є, ўЄю фы  Ёх°хэшщ
\[
    \mathbf{y}_{\varepsilon}=(L_{\varepsilon}-\lambda)^{-1}\mathbf{f}
\]
ёяЁртхфышт√ эхЁртхэёЄтр
\[
    \|\mathbf{y}_{\varepsilon}-\mathbf{y}_0\|_{L_2}\le C\|\mathbf{y}_{\varepsilon}-\mathbf{y}_0\|_{W_1^1}\le C\|Q_{\varepsilon}-Q_0\|_{L_1}
    \|f\|_{L_2}.
\]
╥хь ёрь√ь фюърчрэю ёююЄэю°хэшх~\eqref{eq:T.1.7}, шыш ЁртэюьхЁэр 
Ёхчюы№тхэЄэр  ёїюфшьюёЄ№ юяхЁрЄюЁют \(L_{\varepsilon}\) яЁш
\(\varepsilon\to 0\).
\end{proof}

\begin{Corollary}
╚ч ╥хюЁхь√ \ref{tm:1.6}  ёыхфєхЄ, ўЄю  ёюсёЄтхээ√х чэрўхэш 
юяхЁрЄюЁр $L_{Q,U}$ эхяЁхЁ√тэю чртшё Є юЄ яюЄхэЎшрыр $Q$, т эюЁьх
яЁюёЄЁрэёЄтр $L_1$. ─юърчрЄхы№ёЄтю ёЁрчє яюыєўрхЄё  ё яюью∙№■
ёЄрэфрЁЄэющ Єхїэшъш ё шёяюы№чютрэшхь  яЁюхъЄюЁют ╨шёёр.
\end{Corollary}

{\bf 4.4.} ╥хяхЁ№ фюърцхь ЄхюЁхьє юс рёшьяЄюЄшъх ёюсёЄтхээ√ї
ЇєэъЎшщ  ёшы№эю Ёхуєы Ёэюую юяхЁрЄюЁр.  ─ы  ёыєўр  юс√ўэющ
Ёхуєы ЁэюёЄш эєцэю  ЁрёёьрЄЁштрЄ№ чрфрўє юс рёшьяЄюЄшўхёъюь
яютхфхэшш фтєьхЁэ√ї яюфяЁюёЄЁрэёЄт, юЄтхўр■∙шї ёсышцр■∙шьё 
ёюсёЄтхээ√ь чэрўхэш ь, эю хх  ь√ т ¤Єющ ЁрсюЄх эх ЁрёёьрЄЁштрхь.

\begin{Theorem}\label{tm:2.6}\phantom{.}\!\!\!\footnote{┬ юЁшушэры№эющ ёЄрЄ№х Math. Notes \textbf{96} (5) т ЇюЁьєышЁютъх ╥хюЁхь√ 4.8 шьххЄё  юяхўрЄър: тьхёЄю
$L_{p'}$--эюЁь√ т юяЁхфхыхэшш ўшёхы $b_n$ чряшёрэр $C$--эюЁьр, Є.х. эюЁьр т яЁюёЄЁрэёЄтх эхяЁхЁ√тэ√ї   ЇєэъЎшщ.   ╩юэхўэю, яЁхф°хёЄтє■∙шх Ёхчєы№ЄрЄ√
(ёь. фюърчрЄхы№ёЄтю ╥хюЁхь√ 4.5) уютюЁ Є ю Єюь, ўЄю яЁртшы№эр  юЎхэър фрхЄё  $L_{p'}$-эюЁьющ.} ╧єёЄ№ $Q\in L_p$  яЁш эхъюЄюЁюь $p\ge 1$ ш яєёЄ№
юяхЁрЄюЁ $L_{Q,U}$ ёшы№эю Ёхуєы Ёхэ. ╬сючэрўшь ўхЁхч $\{\mathbf{y}_n(x)\}$ ёюсёЄтхээ√х ЇєэъЎшш ¤Єюую юяхЁрЄюЁр, юЄтхўр■∙шх ёюсёЄтхээ√ь чэрўхэш ь
$\{\lambda_n\}$,  р ўхЁхч $\{\mathbf{y}_n^0(x)\}$   ёюсёЄтхээ√х ЇєэъЎшш (ё эюЁьющ =1) рёёюЎшшЁютрээюую юяхЁрЄюЁр $L_{0,U'}$,  юЄтхўр■∙шх ёюсёЄтхээ√ь
чэрўхэш ь $\{\lambda^0_n\}$. ╥юуфр ёяЁртхфышт√ ЁртхэёЄтр
\begin{equation}\label{efass}
\mathbf{y}_n(x)=E(x)\mathbf{y}_n^0(x)+\bold R_n(x),\qquad
\end{equation}
уфх  ўшёыр $b_n =\|\bold R_n\|_{L_{p'}}$ Єръют√, ўЄю яюёыхфютрЄхы№эюёЄ№ $\{b_n\} \in l_{p'}, \ 1/p +1/p' =1,$  р яЁш $p=1$ яюёыхфютрЄхы№эюёЄ№
$b_n=\|\bold R_n\|_C$ ёЄЁхьшЄё  є эєы■. ┴юыхх Єюую, ёяЁртхфыштю яЁхфёЄртыхэшх
$$
\bold R_n(x)=\mathcal R_n(x)\begin{pmatrix}\cos(\lambda_n^0x),&\sin(\lambda_n^0x)^t\end{pmatrix},
$$
уфх ьрЄЁшЎр  $\mathcal R_n$ Єръютр, ўЄю яЁюшчтюфэ√х хх ¤ыхьхэЄют $r_{jk, n}(x)$   яюфўшэхэ√ юЎхэъх
$$
|r'_{jk, n}(x)|\le C (|q_1(x)| +|q_2(x)+q_3(x)|)
$$
ё яюёЄю ээющ $C$,   эх чртшё ∙хщ юЄ $n$. ╤юсёЄтхээ√х ЇєэъЎшш ёюяЁ цхээюую юяхЁрЄюЁр, юЄтхўр■∙шх ёюсёЄтхээ√ь чэрўхэш ь $\{\overline\lambda_n\}$
фюяєёър■Є яЁхфёЄртыхэшх
\begin{equation}\label{efass'}
\bold z_n(x)=E^{-1}(x)\bold z_n^0(x)+\bold R^*_n(x),\qquad
\end{equation}
уфх $\bold z_n^0(x)$ --- ёюсёЄтхээ√х ЇєэъЎшш юяхЁрЄюЁр, ёюяЁ цхээюую ъ $L_{0,U'}$, эюЁьшЁютрээ√х  ( яЁш сюы№°шї $n$) єёыютшхь
$\langle\mathbf{y}^0_n, \bold z^0_n\rangle =1$. ╧Ёш ¤Єюь ЇєэъЎшш $\bold R^*_n(x)$  єфютыхЄтюЁ ■Є Єхь цх єёыютш ь, ўЄю ш $\bold R_n(x).$
\end{Theorem}
\begin{proof}
─юърчрЄхы№ёЄтю ыхуъю яюыєўрхЄё  шч ╥хюЁхь \ref{tm:2.4} ш \ref{tm:2.3}. ╧ю ёэшь, ъръ ¤Єю ёфхырЄ№.

╧єёЄ№ $\lambda_n^0$ --- юфэю шч ёюсёЄтхээ√ї чэрўхэшщ юяхЁрЄюЁр $L_0$, т юъЁхёЄэюёЄш ъюЄюЁюую эрїюфшЄё  ёюсёЄтхээю чэрўхэшх $\lambda_n$ юяхЁрЄюЁр $L$.
╥юуфр $\Delta_0(\lambda_n^0)=\Delta(\lambda_n)=0$. ╬сючэрўшь ўхЁхч $\bold w^0$ тхъЄюЁ ё фтєь  ъюьяюэхэЄрьш $w_1^0$ ш $w_2^0$, Єръшьш, ўЄю ЇєэъЎш 
$w_1^0\mathbf{c}^0(x,\lambda_n^0)+w_2^0\mathbf{s}^0(x,\lambda_n^0)$  ты хЄё  ёюсёЄтхээющ фы  юяхЁрЄюЁр $L_{0,U'}$. ╥ръ ъръ $\lambda_n^0$ --- яЁюёЄющ
эюы№ ЇєэъЎшш $\Delta_0(\lambda)$, Єю тхъЄюЁ $w^0$ юяЁхфхыхэ ё ЄюўэюёЄ№■ фю ъюэёЄрэЄ√. ┬√схЁхь хх Єръ, ўЄюс√ $\|\bold w^0\|=1$. └эрыюушўэю юяЁхфхышь
тхъЄюЁ $\bold w$, ўЄюс√  ЇєэъЎш  $w_1\mathbf{c}(x,\lambda_n)+w_2\mathbf{s}(x,\lambda_n)$ с√ыр ёюсёЄтхээющ фы  юяхЁрЄюЁр $L_{Q,U}$. ╚ч юяЁхфхыхэш 
ёыхфє■Є ЁртхэёЄтр
$$
M_0(\lambda_n^0)\bold w^0=M(\lambda_n)\bold w=0,
$$
уфх
$$
M_0(\lambda)=\begin{pmatrix}u_{11}&u_{12}\\u_{21}&u_{22}\end{pmatrix}+E\begin{pmatrix}u_{13}&u_{14}\\u_{23}&u_{24}\end{pmatrix}
\begin{pmatrix}c^0_1(\pi,\lambda)&s^0_1(\pi,\lambda)\\c^0_2(\pi,\lambda)&s^0_2(\pi,\lambda)\end{pmatrix},
$$
р ьрЄЁшЎ-ЇєэъЎш  $M(\lambda)$ юяЁхфхыхэр рэрыюушўэю (ёь. ╦хььє 4.2). ─рыхх сєфхь шёяюы№чютрЄ№ юсючэрўхэшх
$$
V_0(\lambda)=\begin{pmatrix}c^0_1(\pi,\lambda)&s^0_1(\pi,\lambda)\\c^0_2(\pi,\lambda)&s^0_2(\pi,\lambda)\end{pmatrix},\qquad
V(\lambda)=\begin{pmatrix}c_1(\pi,\lambda)&s_1(\pi,\lambda)\\c_2(\pi,\lambda)&s_2(\pi,\lambda)\end{pmatrix},
$$
╚ч рёшьяЄюЄшўхёъшї ЁртхэёЄт ёыхфєхЄ $\{\|V(\lambda_n)-V_0(\lambda_n^0)\|\}\in l_{p'}$ (ёь. фюърчрЄхы№ёЄтю ╥хюЁхь√ 4.5)  ╧юёъюы№ъє ЇєэъЎш 
$M_0(\lambda)$ яхЁшюфшўэр, р $\lambda_n^0$
--- яЁюёЄющ эюы№ ¤Єющ ЇєэъЎшш, Єю
$$
\|w^0-w\|\le C_1\|M_0(\lambda_n^0)-M(\lambda_n)\|\le C_2E\|V(\lambda_n)-V_0(\lambda_n^0)\|,
$$
уфх $C_1$ ш $C_2$ чртшё Є Єюы№ъю юЄ ъЁрхт√ї єёыютшщ. ╥юуфр
$$
\mathbf{y}_n(x)=w_1^0\mathbf{c}(x,\lambda_n)+w_2^0\mathbf{s}(x,\lambda_n)+(w_1-w_1^0)\mathbf{c}(x,\lambda_n)+(w_2-w_2^0)\mathbf{s}(x,\lambda_n),
$$
р шёяюы№чє  рёшьяЄюЄшўхёъшх ЇюЁьєы√ \eqref{sinas} ш \eqref{cosas}, яЁшїюфшь ъ ЁртхэёЄтє
$$
\mathbf{y}_n(x)=E(x)(w_1^0\mathbf{c}^0(x,\lambda_n)+w_2^0\mathbf{s}^0(x,\lambda_n))+\bold R^1_n(x),
$$
уфх $\bold R^1(x)$  єфютыхЄтюЁ хЄ ёЇюЁьєышЁютрээ√ь т ╥хюЁхьх єёыютш ь. ─рыхх, яюы№чє ё№ ЄЁшуюэюьхЄЁшўхёъшьш ЇюЁьєырьш ш рёшьяЄюЄшъющ $\{\lambda_n\}$
яюыєўрхь,  ўЄю ЇєэъЎшш
$$
\bold R^2(x) = \mathbf{c}^0(x,\lambda_n) - \mathbf{c}^0(x,\lambda^0_n), \quad \bold R^3(x) = \mathbf{s}^0(x,\lambda_n) - \mathbf{s}^0(x,\lambda^0_n)
$$
єфютыхЄтюЁ ■Є Єхь цх єёыютш ь. ╧Ёш ¤Єюь, хёыш $ b_n^j=\|\bold R_n^j(x)\|_{L_{p'}}$,  Єю яюёыхфютрЄхы№эюёЄш $\{b_n^j\}$  яЁшэрфыхцрЄ $l_{p'}$. ╬Ўхэъш
фы  яЁюшчтюфэ√ї $r'_{jk,n}$ ёЁрчє цх ёыхфє■Є шч юЎхэюъ фы  яЁюшчтюфэ√ї $\rho'_{j,k}$, яюыєўхээ√ї т ЄхюЁхьрї \ref{tm:2.4} ш \ref{tm:2.3}.

╥ръ ъръ ёюяЁ цхээ√щ юяхЁрЄюЁ $L^*_{Q,U}$  Єюцх ёшы№эю Ёхуєы Ёхэ, Єю єЄтхЁцфхэшх ЄхюЁхь√ фы  эхую ёюїЁрэ хЄё , эю эєцэю єўхёЄ№, ўЄю ЇєэъЎш  $E(x)$
фы  ёюяЁ цхээюую яюЄхэЎшрыр $\overline Q$ яЁхюсЁрчєхЄё  т ЇєэъЎш■ $E^{-1}(x)$. ╥хюЁхьр фюърчрэр.
\end{proof}

 {\bf 4.5.}
┬ чръы■ўхэшх фюърцхь ЄхюЁхьє ю срчшёэюёЄш ╨шёёр ёюсёЄтхээ√ї ш яЁшёюхфшэхээ√ї ЇєэъЎшщ ёшы№эю Ёхуєы Ёэюую юяхЁрЄюЁр $L_{Q,U}$ яЁш єёыютшш $Q\in L_1$.
╟рЄхь, шёяюы№чє  юфшэ яЁшхь шч \cite{S1} (яЁшьхэхээ√щ т фЁєующ ёшЄєрЎшш) ш ЄхюЁхьє ╩рЎэхы№ёюэр--╠рЁъєёр--╠рЎрхтр, ь√ яюърцхь, ўЄю ёяЁртхфышт√ь
юёЄрхЄё  єЄтхЁцфхэшх ю срчшёэюёЄш ╨ьёёр шч фтєьхЁэ√ї яюфяЁюёЄЁрэёЄт, юЄтхўр■∙шї ёсышцр■∙шьё  ёюсёЄтхээ√ь чэрўхэш ь. ╬яЁхфхыхэшх срчшёр ╨шёёр шч
яюфяЁюёЄЁрэёЄт (шыш схчєёыютэющ срчшёэюёЄш шч яюфяЁюёЄЁрэёЄт) ёь. т ьюэюуЁрЇшш \cite[уы. 6]{GK}.

\begin{Theorem}\label{tm:2.9}
╧єёЄ№ $Q\in L_1$ ш яєёЄ№ юяхЁрЄюЁ $L_{Q,U}$ ёшы№эю Ёхуєы Ёхэ.
╥юуфр ёшёЄхьр хую ёюсёЄтхээ√ї ш яЁшёюхфшэхээ√ї ЇєэъЎшщ (яЁш
єёыютшш, ўЄю хую ёюсёЄтхээ√х ЇєэъЎшш эюЁьшЁютрэ√ эр 1) юсЁрчє■Є
срчшё ╨шёёр. ┼ёыш  юяхЁрЄюЁ $L_{Q,U}$  Ёхуєы Ёхэ, эю эх ёшы№эю
Ёхуєы Ёхэ, Єю  фтєьхЁэ√х яюфяЁюёЄЁрэёЄтр, юЄтхўр■∙шх ёсышцр■∙шьё 
ёюсёЄтхээ√ь чэрўхэш ь юсЁрчє■Є схчєёыютэ√щ срчшё яЁюёЄЁрэёЄтр
$\mathbb H$.
\end{Theorem}
\begin{proof}
╤эрўрыр чрьхЄшь, ўЄю ёюсёЄтхээ√х ш яЁшёюхфшэхээ√х ЇєэъЎшш $\{\mathbf{y}_n(x)\}$  ш $\{\mathbf{z}_n(x)\}$  юяхЁрЄюЁр $L_{Q,U}$  ш хую ёюяЁ цхээюую
$L^*_{Q,U}$ юсЁрчє■Є яюыэ√х ёшёЄхь√ т яЁюёЄЁрэёЄтх $\mathbb H$.  ▌Єю ёыхфєхЄ шч юЎхэъш Ёхчюы№тхэЄ√ ¤Єшї юяхЁрЄюЁют, тэх эхъюЄюЁющ яюыюё√
$|\text{Im}\lambda|>\alpha$, яюыєўхээющ т ╥хюЁхьх 4.1 ш Єюую ЇръЄр, ўЄю Ёхчюы№тхэЄ√ ¤Єшї юяхЁрЄюЁют  ты ■Єё  шэЄхуЁры№э√ьш юяхЁрЄюЁрьш, ЇєэъЎшш ├Ёшэр
ъюЄюЁ√ї ЁртэюьхЁэю юуЁрэшўхэ√ тэх ъЁєцъют ьрыюую Ёрфшєёр $\varepsilon$ ё ЎхэЄЁрьш т ёюсёЄтхээ√ї чэрўхэш ї $\lambda_n$. ╬уЁрэшўхээюёЄ№ ЇєэъЎшщ ├Ёшэр т
яюыюёх (Єрь, уфх шьх■Єё  рёшьяЄюЄшъш ЇєэфрьхэЄры№эющ ёшёЄхь√ Ёх°хэшщ) яюыєўрхЄё  Єюўэю Єръцх, ъръ ¤Єю ёфхырэю т \S 3   фы  ЇєэъЎшш $G_0(x,\xi,
\lambda)$ (яюфёЄрэютър ЇєэъЎшщ $\mathbf{c}(x,\lambda$ ш $\mathbf{s}(x,\lambda)$  тьхёЄю $\mathbf{c}^0(x,\lambda$ ш $\mathbf{s}^0(x,\lambda)$ яЁш
яюёЄЁюхэшш ЇєэъЎшш ├Ёшэр эшўхую эх ьхэ хЄ т фюърчрЄхы№ёЄтх яЁш эрышўшш ёююЄтхЄёЄтє■∙шї рёшьяЄюЄшъ). ╧ю¤Єюьє яюыэюЄр ёшёЄхь $\{\mathbf{y}_n(x)\}$ ш
$\mathbf{z}_n(x)$ яюыєўрхЄё  Єюўэю Єръцх, ъръ яЁш $Q=0.$

╥хяхЁ№ яЁхфяюыюцшь, ўЄю юяхЁрЄюЁ  $L_{Q,U}$ ёшы№эю Ёхуєы Ёхэ. ─юърцхь, ўЄю юсх  ёшёЄхь√ $\{\mathbf{y}_n(x)\}$  ш $\{\mathbf{z}_n(x)\}$ схёёхыхт√.
┬юёяюы№чєхьё  ╥хюЁхьющ \ref{tm:2.6}, ъюЄюЁр  фрхЄ яЁхфёЄртыхэш  фы  ЇєэъЎшщ ¤Єшї ёшёЄхь. ┬ ¤Єшї яЁхфёЄртыхэш ї ь√ ёўшЄрхь, ўЄю
$\|\mathbf{y}_n^0(x)\|=1$ ш $\langle \mathbf{y}_n^0, \mathbf{z}^0_n\rangle=1$. ╙цх  фюърчрэю, ўЄю ёшёЄхь√ $\{\mathbf{y}_n^0(x)\}$ ш
$\{\mathbf{z}^0_n(x)\}$, яюфўшэхээ√х єърчрээющ эюЁьшЁютъх,  ты ■Єё  схёёхыхт√ьш. ═ю Єюуфр
  ёшёЄхь√, ёюёЄртыхээ√х шч  яхЁт√ї ёырурхь√ї  $\{E(x)\mathbf{y}_n^0(x)\}$ ш
$\{E^{-1}(x)\mathbf{z}^0_n(x)\}$ т яЁхфёЄртыхэш ї \eqref{efass}  ш \eqref{efass'},  Єръцх  ты ■Єё  схёёхыхт√ьш (єьэюцхэшх тёхї ЇєэъЎшщ ёшёЄхь√  эр
юуЁрэшўхээє■ ЇєэъЎш■ юёЄрты хЄ ёшёЄхьє схёёхыхтющ). ╬ёЄрхЄё  яюърчрЄ№, ўЄю ёшёЄхь√, ёюёЄртыхээ√х шч тЄюЁ√ї ёырурхь√ї $\{\bold R_n(x)\}$ ш $\{\bold
R^*_n(x)\}$  ты ■Єё  схёёхыхт√ьш. ─юърцхь ¤Єю ётющёЄтю фы  яхЁтющ ёшёЄхь√, ётющёЄтр тЄюЁющ ёшёЄхь√ Єх цх ёрь√х.

╚ч яЁхфёЄртыхэшщ фы  ЇєэъЎшщ $\bold R_n(x)$ ёыхфєхЄ, ўЄю
фюёЄрЄюўэю фюърчрЄ№ ёїюфшьюёЄ№ Ё фр
\begin{equation}\label{sum1}
\sum_n\left|\int\limits_0^\pi f(x)r_n(x)\sigma(\lambda_n^0 x)dx\right|^2,
\end{equation}
уфх $f\in L_2[0,\pi], \ \, \sigma(\lambda_n^0)=\cos(\lambda_n^0 x)$ шыш $\sin(\lambda_n^0 x)$, р ЇєэъЎш  $r_n(x) \in W_1^1[0,\pi]$
--- юфэр шч ъюьяюэхэЄ ьрЄЁшЎ√ $\mathcal R_n(x)$,  яюфўшэхэр єёыютш■
$$
|r'_n(x)|\le
C\left(|q_1(x)|+\frac12|q_2(x)|+\frac12|q_3(x)|\right),
$$
уфх яюёЄю ээр  $C$  эх чртшёшЄ юЄ $n$. ▌Єю чрьхўрэшх ёыхфєхЄ шч Єюую, ўЄю Ё ф $\sum|\langle\mathbf{f}(x), \bold R_n(x)\rangle|^2 $ яЁхфёЄрты хЄё  т
тшфх ўхЄ√Ёхї Ё фют тшфр \eqref{sum1}.

╬сючэрўшь ёєььє $|q_1(x)|+\frac12|q_2(x)|+\frac12|q_3(x)|$ ўхЁхч
$q(x)$. ─рыхх сєфхь шёяюы№чютрЄ№ юЎхэъє $|r'_n(x)| \le Cq(x)$.
╤єььє\eqref{sum1} яхЁхяш°хь т  тшфх
\begin{gather*}
\sum_n\left|\int\limits_0^\pi f(x)\left[r_n(0) +\int\limits_0^x r'_n(t)dt\right]\sigma(\lambda_n^0 x)dx\right|^2\le\\ \le 2\sum_n
\left|r_n(0)\int\limits_0^\pi f(x)\sigma(\lambda_n^0 x)dx\right|^2 + 2 \sum_n\left|\int\limits_0^\pi r'_n(t)\int\limits_t^\pi f(x)\sigma(\lambda_n^0
x)\, dx\, dt\right|^2.
\end{gather*}
╟фхё№ яхЁт√щ Ё ф  т яЁртющ ўрёЄш ёїюфшЄё , яюёъюы№ъє ёшёЄхьр $\{\sigma(\lambda_n^0 x)\}$ схёёхыхтр, р ўшёыр $r_n(0)$ юуЁрэшўхэ√ (фрцх ёЄЁхь Єё  ъ
эєы■). ╧Ёюфюыцшь юЎхэъє тЄюЁюую Ё фр. ╧хЁхяш°хь хую т тшфх
\begin{gather*}
\sum_n\left|\int\limits_0^\pi\int\limits_0^\pi r'_n(t) r'_n(s)\int\limits_t^\pi f(x)\sigma(\lambda_n^0 x)dx\int\limits_s^\pi f(y)\sigma(\lambda_n^0
y)dy\,dsdt\right|\le\\
\le C^2\int\limits_0^\pi\int\limits_0^\pi q(s)q(t)\sum_n\left|\int\limits_t^\pi f(x)\sigma(\lambda_n^0 x)dx\int\limits_s^\pi f(y)\sigma(\lambda_n^0
y)dy\right|dsdt \le C^2 C_1\left(\int\limits_0^\pi q(x)\,dx\right)^2 \|f\|^2,
\end{gather*}
яюёъюы№ъє   ёшёЄхьр $\{\sigma(\lambda_n^0 x)\}$ схёёхыхтр ш ёєььр яюф шэЄхуЁрырьш ЁртэюьхЁэю яю $ s,t\in [0,\pi]\times [0,\pi]$  юЎхэштрхЄё 
тхышўшэющ (ўхЁхч $\chi_{[t,\pi]}$  юсючэрўрхь їрЁръЄхЁшёЄшўхёъє■ ЇєэъЎш■ юЄЁхчър $[t,\pi]$)
\begin{gather*}
\sum_n\left|\int\limits_0^\pi f(x)\chi_{[t,\pi]}\sigma(\lambda_n^0
x)dx\int\limits_0^\pi f(y)\chi_{[s,\pi]}\sigma(\lambda_n^0 y)dy\right|\le \\
\le \sum_n\left|\int\limits_0^\pi f(x)\chi_{[t,\pi]}\sigma(\lambda_n^0 x)dx\right|^2 + \left|\int\limits_0^\pi f(y)\chi_{[s,\pi]}\sigma(\lambda_n^0
y)dy\right|^2 \le  C_1 \|f\|^2.
\end{gather*}
╥хь ёрь√ь фюърчрэр схёёхыхтюёЄ№ юсхшї ёшёЄхь $\{\mathbf{y}_n(x)\}$  ш $\{\mathbf{z}_n(x)\}$. ╚ч яЁхфёЄртыхэшщ \eqref{efass} ш \eqref{efass'}
 яЁш сюы№°шї $n$ яюыєўрхь
$$
\alpha_n: =\langle \mathbf{y}_n(x),\mathbf{z}_n(x)\rangle=\langle \mathbf{y}^0_n(x),\mathbf{z}^0_n(x)\rangle +o(1) =  1+o(1).
$$
╤ыхфютрЄхы№эю, ёшёЄхьр  $\{\alpha_n^{-1}\mathbf{z}_n(x)\}$ (ъюэхўэю, юэр Єюцх схёёхыхтр)  ты хЄё  сшюЁЄюуюэры№эющ  ъ  $\{\mathbf{y}_n(x)\}$. ╥хяхЁ№
тюёяюы№чєхьё  ЄхюЁхьющ ┴рЁш-┴юрёр \cite[уы. 6]{GK}: {\sl ┼ёыш ёшёЄхь√ $\{\varphi_n\}$  ш  $\{\psi_n\}$  тчршьэю сшюЁЄюуюэры№э√, яюыэ√ ш схёёхыхт√, Єю
юэш юсх  ты ■Єё  срчшёрьш ╨шёёр.}

╟рьхЄшь, ўЄю ёююЄэю°хэш  сшюЁЄюуюэры№эюёЄш ь√ т√яшё√трыш яЁш сюы№°шї $n$, ъюуфр ёюсёЄтхээ√ї чэрўхэш  яЁюёЄ√х ш яЁшёюхфшэхээ√ї эхЄ. ╬фэръю т ёшы№эю
Ёхуєы Ёэюь ёыєўрх яЁшёюхфшэхээ√ї ЇєэъЎшщ ьюцхЄ с√Є№ Єюы№ъю ъюэхўэюх ўшёыю, р яюЄюьє  эшъръюую тыш эш  эр їюф фюърчрЄхы№ёЄтр юэш эх юърч√тр■Є.

╥хяхЁ№ яЁхфяюыюцшь, ўЄю юяхЁрЄюЁ $L=L_{Q,U}$  Ёхуєы Ёхэ, эю эх ёшы№эю Ёхуєы Ёхэ.  ─юсртшь ъ $Q$ юуЁрэшўхээ√щ яюЄхэЎшры $Q_0$, Єръющ, ўЄюс√  $L_{Q
+Q_0, \, U}$ с√ы ёшы№эю Ёхуєы Ёхэ. ▌Єюую ьюцэю фюсшЄ№ё  яюфсшЁр  тэхфшруюэры№э√х ЇєэъЎшш Єръ, ўЄюс√ $E(\pi)\ne 1$. ╥юуфр ёюсёЄтхээ√х ш яЁшёюхфшэхээ√х
ЇєэъЎшш юяхЁрЄюЁр $L+Q_0$ юсЁрчє■Є срчшё ╨шёёр. ╥ръ ъръ яЁшёюхфшэхээ√ї ЇєэъЎшщ ъюэхўэюх ўшёыю, Єю эрщфхЄё  ъюэхўэюьхЁэ√щ юяхЁрЄюЁ $K$,  Єръющ, ўЄю
$L+Q_0+K$  шьххЄ Єюы№ъю ёюсёЄтхээ√х ЇєэъЎшш. ╥юуфр т ёъры Ёэюь яЁюшчтхфхэшш, Єюяюыюушўхёъш ¤ътштрыхэЄэюьє шёїюфэюьє, юэ  ты хЄё  эюЁьры№э√ь
юяхЁрЄюЁюь, р ё єўхЄюь, ўЄю хую ёюсёЄтхээ√х чэрўхэш  ыхцрЄ т яюыюёх, юэ яЁхфёЄртшь т тшфх $T +B$, уфх $T$ ёрьюёюяЁ цхэ, р $B$  юуЁрэшўхэ.
╤ыхфютрЄхы№эю, $L = T+B -Q_0 -K$   ты хЄё  юуЁрэшўхээ√ь тючьє∙хэшхь ёрьюёюяЁ цхээюую юяхЁрЄюЁр $T$, ёюсёЄтхээ√х чэрўхэш  ъюЄюЁюую   шьх■Є рёшьяЄюЄшъє
$\mu_n = n +O(1)$.  ═ю хёыш ёюсёЄтхээ√х чэрўхэш  $\mu_n$  ёрьюёюяЁ цхээюую юяхЁрЄюЁр $T$   єфютыхЄтюЁ ■Є эхЁртхэёЄтрь $|\mu_n| \ge  c |n|$, Єю
ёюуырёэю ЄхюЁхьх ╩рЎэхы№ёюэр-╠рЁъєёр-╠рЎрхтр ( ёь. \cite {Kaz, MM} {\sl ёшёЄхьр ёюсёЄтхээ√ї ш яЁшёюхфшэхээ√ї ЇєэъЎшщ юуЁрэшўхээюую тючьє∙хэш 
юяхЁрЄюЁр $T$ юсЁрчєхЄ срчшё ╨шёёр ёю ёъюсърьш (шыш срчшё ╨шёёр шч яюфяЁюёЄЁрэёЄт)}. ╬ЄьхЄшь, ўЄю т єърчрээ√ї ЁрсюЄрї эх юуютюЁхэю ёяхЎшры№эю, ъръ
ЁрёёЄрты Є№ ёъюсъш шыш юс·хфшэ Є№ тхъЄюЁ√ т яюфяЁюёЄЁрэёЄтр, эю шч їюфр фюърчрЄхы№ёЄтр (ёь. \cite{M1}) ыхуъю шчтыхў№, ўЄю т эр°хщ ёшЄєрЎшш юс·хфшэ Є№
эрфю ЇєэъЎшш, юЄтхўр■∙шх ёсышцр■∙шьё  ёюсёЄтхээ√ь чэрўхэш ь. ╥хюЁхьр фюърчрэр.
\end{proof}

\end{document}